\documentclass[10pt]{amsart}
\usepackage{amsmath,amssymb,amsfonts}
\usepackage{enumerate}
\usepackage[pagebackref=true]{hyperref}
\hypersetup{
   colorlinks=true,
    linkcolor=blue,
   filecolor=magenta,      
    urlcolor=cyan,
    pdftitle={Overleaf Example},
    pdfpagemode=FullScreen,
    }
\usepackage{color}
\usepackage[dvipsnames]{xcolor}
\definecolor{sz}{RGB}{115,45,2}

\newcommand{\tb}[1]{\textcolor{sz}{\bf #1}}
\usepackage{geometry}
\usepackage{tikz}
\usepackage{tikz-cd}
\usepackage{mathtools}

\usepackage{epsfig}  		
\usepackage{epic,eepic}  
\usepackage{graphicx}
\usepackage{verbatim}
\usepackage{mathdots}
\usepackage{leftidx}

\newcommand{\w}{\mathrm{w}}
\newcommand{\rightarrowdbl}{\rightarrow\mathrel{\mkern-14mu}\rightarrow}

\theoremstyle{plain}
\newtheorem{thm}{Theorem}[section]

\newtheorem{lem}[thm]{Lemma}
\newtheorem{prop}[thm]{Proposition}

\theoremstyle{remark}

\newtheorem{remark}[thm]{Remark}

\DeclareMathOperator{\Hom}{Hom}

\DeclareMathOperator{\Gal}{Gal}

\begin{document}
\title{Stability of Rankin-Selberg local $\gamma$-factors for split classical groups: The symplectic case}
\author{Taiwang Deng, Dongming She}
\address{        
Morningside Center of Mathematics, Chinese Academy of Sciences, No.55, Zhongguancun East Road, Haidian District, Beijing 100190, China.}
\email{shed@amss.ac.cn}
\address{ 
Yanqi Lake Beijing Institute of Mathematical Sciences and Applications (BIMSA), Huairou District, 100084, Beijing, China.}
\email{dengtaiw@bimsa.cn}
\maketitle

\begin{abstract}
    Given a split classical group of symplectic type and a split general linear group over a local field $F$, we use Langlands-Shahidi method to construct their Rankin-Selberg local $\gamma$-factors and prove the corresponding analytic stability for generic representations. The idea generalizes the work of J. Cogdell, F. Shahidi, T.-L. Tsai in 2017 and D. She in 2023 in the study of asymptotic behaviors of partial Bessel functions. Different from the known cases, suppose $P=MN$ is the maximal parabolic subgroup with Levi component $M\simeq \mathrm{GL}_r\times\mathrm{Sp}_{2m}$ that defines the local factors, the action of the maximal unipotent subgroup of $M$ on $N$ have non-trivial stabilizers, and the space of integration for the corresponding local coefficient is no longer isomorphic to a torus. We will separate its toric part out in our cases and show that it plays the same role as the torus over which the integral representing the local coefficient is taken in the known cases. This is a new phenomenon with sufficient generality and we believe that it may provide us with a possible direction towards a uniform proof of stability of Langlands-Shahidi $\gamma$-factors in our future work.    
\end{abstract}
\tableofcontents
\section{Introduction}

Let $M$ be a connected reductive group defined over a local field $F$, $\pi$ an irreducible admissible represetation of $M(F)$. Fix a non-trivial additive character $\psi$ of $F$. In the study of Langlands functoriality, or more precisely the L-function theory, one needs to consider a local $\gamma$-factor $\gamma(s,\pi, r,\psi)$ for any finite dimensional complex representation $r$ of the Langlands L-group ${^LM}=\hat{M}\rtimes \Gamma_F$, where $\Gamma_F=\Gal(\overline{F}/F)$, or $W_F$, the local Weil group of $F$. The general definition of $\gamma(s,\pi,r,\psi)$ is unknown so far, but in special cases there are many methods of constructing them, a very important one among which is the Langlands-Shahidi method \cite{sha90}, \cite{sha10}. It defines the local $\gamma$-factors $\gamma^{sh}(s, \pi, r,\psi)$ when $\pi$ is $\psi$-generic, namely, $\pi$ admits a non-zero Whittaker model, $M$ appears as a Levi subgroup of a larger reductive group $G$, and $r$ an irreducible constituent of the adjoint action of ${^LM}$ on the Lie algebra of ${^LN}$, the L-group of $N$ where $N$ is the unipotent radical of the parabolic $P=MN$.

The local $\gamma$-factors $\gamma(s,\pi,r,\psi)$, once defined, are related to the local $\epsilon$- and L-factors in the following way:
$$\gamma(s, \pi,r,\psi)=\epsilon(s,\pi,r,\psi)\frac{L(1-s,\tilde{\pi}, r)}{L(s,\pi, r)}$$ The local arithmetic and analytic $\gamma$-factors are expected to satisfy some stable equality under highly ramified twists, called the arithmetic and analytic stability respectively. The arithmetic stability is fully proved by P. Deligne \cite{Deligne73}, but the analytic stability is only known for certain cases. To be precise, if $\pi_1$ and $\pi_2$ are irreducible admissible representations of $M(F)$ sharing the same central character, then for a highly ramified character $\chi$ of $F^\times$, regarded as a character of $M(F)$ via $m\mapsto \chi(\det\mathrm{Ad}_\mathfrak{n}(m))$, where $\mathrm{Ad}: M\rightarrow \mathrm{GL}(\mathfrak{n})$ is the adjoint representation and $\mathfrak{n}:=\mathrm{Lie}(N)$, we expect to have
$$\gamma(s,\pi_1\otimes \chi, r, \psi)=\gamma(s, \pi_2\otimes\chi, \psi).$$

We care about the analytic stability since it serves as an important intermediate step to prove many crucial results in the Langlands program. For example, local converse theorems, functoriality conjecture, and the local Langlands correspondence. Many results are known in these areas. We apologize not to be able to list all of them due to the large number of authors and their work in these directions. 

The analytic stability of local $\gamma$-factors is known in many cases: for $M=\mathrm{GL}_n\times \mathrm{GL}_m$ by Jacquet \& Shalika \cite{Jac85}; for $M$ a $F$-split classical group, either of symplectic or orthogonal type, by Rallis \& 
 Soury \cite{Ral05} via the doubling method; for $M=\mathrm{GL}_n$, $r=\mathrm{Sym}^2$ or $\wedge^2$, by Cogdell, Shahidi \& Tsai \cite{CST17}; for the twisted symmetric and exterior square local factors of $\mathrm{GL}_n$, by She \cite{She23}; for Asai local factors, by Shankman \cite{Shan18}; for $M=\mathrm{GL}_2$, $r=\mathrm{Sym}^3$, by Shankman \& She \cite{Shan19}; for $M$ a general quasi-split reductive group, by Cogdell, Piatetski-Shapiro \& Shahidi \cite{Cog08}, but under the assumptions that $\dim (U_M\backslash N)=2$ and $\mathrm{rank}\{Z_G\backslash T_{w}\}=2$, where $T_{w}$ is the subtorus defined in (3.6) of \cite{Cog08}.

 In this paper, we establish the stability for the Rankin-Selberg product local $\gamma$-factor attached to a classical group of symplectic type, and a general linear group. Namely, we prove the following result:

\begin{thm}\label{thm1.1} Let $M=\mathrm{GL}_r\times \mathrm{Sp}_{2m}$ be split over a p-adic field $F$, $\sigma_{i}$ and $\tau_i, (i=1,2)$ are irreducible $\psi$-generic representations of $\mathrm{GL}_r(F)$ and $\mathrm{Sp}_{2m}(F)$ respectively. Let $\pi_i=\sigma_i\boxtimes \tau_i$($i=1,2$), assume $\omega_1=\omega_2$, where $\omega_i$ is the central character of $\pi_i$. Take a continuous character $\chi: F^\times \rightarrow \mathbb{C}^\times$, regarded as a character of $M(F)$ via $(m_1, m_2)\mapsto \chi(\det(m_1))$ for $m_1\in \mathrm{GL}_r(F)$ and $m_2\in \mathrm{Sp}_{2m}(F)$. Assume $\chi$ is sufficiently ramified, then we have
$$\gamma(s, (\sigma_1\times \tau_1)\otimes \chi, \psi )=\gamma(s, (\sigma_2\times \tau_2)\otimes \chi, \psi))$$
\end{thm}

The strategy is to use Langlands-Shahidi method to define the local $\gamma$-factors, and reduce the stability to the stability of the corresponding local coefficient. Similar to \cite{CST17}, the local coefficient admits an integral representation as the Mellin transform of certain partial Bessel functions, whose asymptotic expansion via relevant Bruhat cells breaks into a sum of two parts, the first part depends only on the central character, the second part is a uniform smooth function on certain subtorus. Hence when twisted by a highly ramified character, the second part becomes zero. This gives the stability of the local coefficient.

We point out two main differences in our cases compared to the known cases using this method. Firstly, the action of $U_M$ on $N$ have non-trivial stabilizers, and consequently the geometry of the orbit space $U_M\backslash N$ becomes more subtle. Secondly, the orbit space $U_M\backslash N$, over which the integral of the partial Bessel function is taken, is no longer isomorphic to a torus in our cases. Therefore a similar argument to \cite{CST17} can not be directly applied. Instead, we separate the toric part out of the orbit space $U_M\backslash N$, and show that only the toric part accounts for the proof of stability. We also observe that the map $n\mapsto m$ via the Bruhat decomposition $\dot{w}_0^{-1}n=mn'\bar{n}$ is finite $\acute{e}tale$ onto its image with covering group isomorphic to finitely many copies of $\mathbb{Z}/2\mathbb{Z}$. This is a new phenomenon with sufficient generality and we believe it may show us one possible direction towards a uniform proof of stability of Langlands-Shahidi $\gamma$-factors in our future work. We also generalize the asymptotic analysis for partial Bessel integrals in \cite{CST17} and \cite{She23}. Lastly, we remark that the cases where we replace the symplectic groups by odd orthogonal groups are very similar. A slight modification of our arguments would also give a proof of the same result for odd orthogonal cases.

Finally, many of our computations rely on the computational softwares Sagemath, Mathematica and GP/PARI.

\par \vskip 1pc
{\bf Acknowledgement.}
The work was started when the first author was a postdoc at the Yau Mathematical Sciences Center (YMSC) and was finished at Yanqi Lake Beijing Institute of Mathematical Sciences and Applications (BIMSA) where he became an assistant fellow. He wants to thank both institutes for their hospitality. At the same time, the second author was a postdoc at the Morningside Center of Mathematics(MCM), he would like to also thank the institute for its hospitality as well, and for providing him a first-rate research environment. 

Many ideas of the work are based on the second author's Ph.D theis under the supervision of Prof. Freydoon Shahidi. The second author sincerely expresses his deep sense of gratitude to Prof. Freydoon Shahidi for his invaluable guidance. 
 
\section{Construction of the Rankin-Selberg local $\gamma$-factors} Consider the split reductive group $G=G(n)=\mathrm{Sp}_{2n}$ defined over a p-adic field $F$ which is a finite extension $\mathbb{Q}_p$. We realize $G$ as
$$\mathrm{Sp}_{2n}=\{h\in\mathrm{GL}_{2n}: {^th}J'h=J'\}$$
where $$J'=J'_{2n}:=\begin{bmatrix}
    & J_n \\
    -{^tJ}_n & \\ 
\end{bmatrix} \mathrm{\  and \ } J_n=\begin{bmatrix}
     & &  & 1 \\
     & & -1 & \\
     & \iddots & \\
     (-1)^{n-1} & & & \\
\end{bmatrix}.$$ It is a semi-simple group of type $C_n$. Fix a Borel subgroup $B=TU$ consisting of the upper triangular matrices in $\mathrm{Sp}_{2n}$, then
$$T=\{t=\mathrm{diag}(t_1, \cdots, t_n, t_n^{-1},\cdots, t_1^{-1}): t_i\in \mathbb{G}_m\}.$$

Following Bourbaki's labeling \cite{Bou08} of the root systems, the set of positive roots is given by $$\Phi^+=\{\alpha_1,\cdots, \alpha_n, \sum_{i\le k< j}\alpha_k (1\leq i<j\leq n), \sum_{i\leq k<j}\alpha_i+2\sum_{j\leq k<n}\alpha_k+\alpha_n (1\leq i<j\le n), 2\sum_{i\leq k <n}\alpha_k+\alpha_n (1\leq i\le n)\}.$$ Note that all the standard maximal parabolic subgroups of $\mathrm{Sp}_{2n}$ are self associate. Indeed, let $\Delta_r=\Delta-\{\alpha_r\}$, $ P_{\Delta_r}=M_{\Delta_r}N_{\Delta_r}$, and $w_0=\w_Gw_{M_{\Delta_r}}^{-1}$ where $w_G$ and $w_{M_{\Theta_r}}$ are the long Weyl group element of $G$ and $M_{\Delta_r}$ respectively, then $w_0(\Delta_r)=\Delta_r$, and $w_0(\alpha_r)<0$. For simplicity, we write $P=MN$ from now on. The Levi components of the maximal parabolic subgroups of $G(n)$ are of the form $\mathrm{GL}_r\times \mathrm{Sp}_{2m}$, with $r+m=n$, obtained by removing the simple root $\alpha_r(1\leq r\leq n)$ from the set of simple roots. (When $r=n$, the corresponding maximal Levi is isomorphic to $\mathrm{GL}_n$, by $\mathrm{Sp}_0$ we mean the trivial group 1.) We realize $M\simeq \mathrm{GL}_r\times\mathrm{Sp}_{2m}$ via
$$m=\begin{bmatrix}
    m_1 &  & \\ 
        & m_2 & \\
        &     &  J_r {^tm_1^{-1}} J_r^{-1}\\
\end{bmatrix}\mapsto (m_1,m_2)$$ where $m_1\in \mathrm{GL}_r$, $m_2\in \mathrm{Sp}_{2m}$. We also have

$$N=\{n=\begin{bmatrix}
    I_r &  X  & Y \\
        & I_{2m} & Z \\
        &     & I_r \\
\end{bmatrix}\in \mathrm{Sp}_{2n}: X\in \mathrm{Mat}_{r\times 2m}, Y\in \mathrm{Mat}_{r\times r}\}$$
$$=\{n(X,Y):=\begin{bmatrix}
    I_r & X & Y\\
    & I_{2m} & (-1)^rJ'_{2m} {^tX}J_r \\
     & & I_r\\
\end{bmatrix}: J_r{^tY}-Y{^tJ_r}+(-1)^rX J'_{2m} {^tX}=0\}.$$

 We will use the local theory of Langlands-Shahidi method to construct the local Rankin-Selberg $\gamma$-factors in our cases. Fix a non-trivial additive character $\psi$ of $F$, and suppose $\sigma$ and $\tau$ are irreducible $\psi$-generic representations of $\mathrm{GL}_r(F)$ and $\mathrm{Sp}_{2m}(F)$ respectively. Set $\pi\simeq \sigma\boxtimes \tau$, then its central character $\omega_\pi\simeq \omega_\sigma\boxtimes\omega_{\tau}$. Fix a non-zero Whittaker functional $\lambda\in \Hom_{U_M(F)}(\pi, \psi)$, i.e., 
 $$\lambda(\pi(u)v)=\chi(u)\lambda(v)$$ for any $v\in \pi$ and $u\in U_M(F)$. Suppose that $\dot{w}_0\in G(F)$ is a representative of $\omega_0=\omega_G\cdot \omega_M^{-1}$, such that $\dot{w}_0=\dot{w}_G\dot{w}_M^{-1}$ is compatible with $\psi$, in the sense that $\psi(\dot{w}_0u\dot{w}_0^{-1})=\psi(u)$ for all $u\in U_M(F)$.  For $\nu\in \mathfrak{a}_{P,\mathbb{C}}^*$, where $\mathfrak{a}_{P}^*=\Hom(X^*(M)_F,\mathbb{R})$, and $\mathfrak{a}_{P,\mathbb{C}}^*=\mathfrak{a}_{P}^*\otimes_{\mathbb{R}}\mathbb{C}$, let $$I(\nu,\pi)=\mathrm{Ind}_{M(F)N(F)}^{G(F)}\pi\otimes q_F^{\langle \nu+\rho_P, H_P(\cdot)\rangle}\otimes\textbf{1}_{N(F)}$$ be the normalized parabolic induction. Then the functional $$\lambda_{\psi}(\nu,\pi):I(\nu,\pi)\longrightarrow \mathbb{C}$$ 
$$\lambda_\psi(\nu, \pi)f=\int_{N'(F)}\lambda(f(\dot{w}_0^{-1} n'))\overline{\chi}(n')dn'$$ where $N'=w_0 N^-w_0^{-1}$, defines a non-zero Whittaker functional on $I(\nu,\pi)$. Suppose that $w\in W(G,T)$, such that $\Theta'=w(\Theta)\subset \Delta$ with $\dot{w}\in G(F)$ is chosen to be compatible with $\psi$, then $$A(\nu, \pi,\dot{w}):I(\nu,\pi)\longrightarrow I(\dot{w}\nu, \dot{w}(\pi))$$
 $$A(\nu,\pi,\dot{w})f=(g\mapsto \int_{N_{\dot{w}(F)}}f(\dot{w}^{-1}ng)dn)$$ where $N_{\dot{w}}=U \cap \dot{w}N^-\dot{w}^{-1}$, defines an intertwining operator between the induced representations. Consequently, $\lambda_\psi(\dot{w}(\nu),\dot{w}(\pi))\circ A(\nu,\pi,\dot{w})$ is another non-zero Whittaker functional on $I(\nu,\pi)$. By uniqueness of local Whittaker functionals, there exists a non-zero constant $C_\psi(\nu,\sigma,\dot{w})$, called the Shahidi local coefficient, such that
$$\lambda_\psi(\nu,\pi)=C_\psi(\nu,\pi,\dot{w})\lambda_\psi(\dot{w}(\nu),\dot{w}(\pi))\circ A(\nu,\pi,\dot{w}).$$

In our cases, $P=MN=P_{\Delta_r}=M_{\Delta_r} N_{\Delta_r}$ is self-associate for all $1\leq r\leq n$. Let $\nu=s\tilde{\alpha}$, and we fix $$\dot{w}_G=\begin{bmatrix}
    &  &  J_r \\
    & (-1)^rJ'_{2m} & \\
    -{^tJ_r} &  & \\
\end{bmatrix}, \dot{w}_M=\begin{bmatrix}
    J_r & & \\
     & (-1)^rJ'_{2m} & \\
      &  & J_r {^tJ_r^{-1}} J_r^{-1} \\
\end{bmatrix}=\begin{bmatrix}
    J_r & & \\
    & (-1)^rJ'_{2m} & \\
    &  & J_r
\end{bmatrix}$$ since ${^tJ_r^{-1}}=J_r$. Then we pick the representative $\dot{w}_0$ to be
$$\dot{w}_0=\dot{w}_G\dot{w}_M^{-1}=\begin{bmatrix}
     & & I_r \\
     & I_{2m} & \\
     (-1)^rI_r & & \\
\end{bmatrix}$$ A straightforward computation shows that $\psi(\dot{w}_0 u\dot{w}_0^{-1})=\psi(u)$ for all $u\in U_M$, i.e., our choice of $\dot{w}_0$ is compatible with the generic character $\psi$. Denote by $A(s,\pi):=A(s\tilde{\alpha},\pi, \dot{w}_0)$ and $I(s,\pi):=I(s\tilde{\alpha}, \pi)$
then
$$A(s,\pi):I(s,\pi)\longrightarrow I(-s,\dot{w}_0(\pi))$$ and the local coefficient $C_\psi(s,\pi):=C_\psi(s\tilde{\alpha},\pi,\dot{w}_0)$ is a product of two $\gamma$-factors, namely,
$$C_\psi(s,\pi)=\gamma(s,\sigma\times\tau,\psi)\cdot \gamma(s, \sigma,\wedge^2,\psi)$$ where $\gamma(s,\sigma, \wedge^2,\psi)$ is the exterior square local factor of $\mathrm{GL}_r$ attached to $\sigma$, which is defined in \textcolor{blue}{\cite{CST17}} also by Langlands-Shahidi method. Therefore we have defined the Rankin-Selberg local factor $\gamma(s,\sigma\times\tau, \psi)$ for the reductive groups $\mathrm{GL}_r$ and $\mathrm{Sp}_{2m}$. Moreover, since the analytic stability of $\gamma(s,\sigma,\wedge^2,\psi)$ is established in \cite{CST17}, it reduces the stability of $\gamma(s,\sigma\times \tau,\psi)$ to the stability of the corresponding local coefficient $C_\psi(s,\pi)$. 

\section{A Bruhat decomposition} From the construction of the local Whittaker functionals on the space of induced representations, we would like to write $\dot{w}_0^{-1}n\in P\overline{N}$ where $\overline{N}=\dot{w}_0N\dot{w}_0^{-1}$ is the opposite of $N$. This decomposition does not hold for all $n\in N$, but holds for a Zariski open dense subset $N'$ of $N$. By our choice of the representative $\dot{w}_0$ of $w_0$, for $t=\mathrm{diag}(t_1, t_2,\cdots, t_n, t_n^{-1}, \cdots, t_2^{-1}, t_1^{-1})\in T$, we have
$$\dot{w}_0 t\dot{w}_0^{-1}=\mathrm{diag}(t_r^{-1},\cdots, t_1^{-1}, t_{r+1},\cdots, t_n, t_n^{-1},\cdots, t_{r+1}^{-1}, t_1,\cdots, t_r).$$ 
From this we observe that
$$w_0(\alpha_i)=\alpha_{r-i}, \mathrm{\ for \ } 1\leq i\leq r-1,$$
$$w_0(\alpha_r)=-((\alpha_1+\alpha_2+\cdots+ \alpha_r)+2(\alpha_{r+1}+\cdots+\alpha_{n-1})+\alpha_n).$$ Therefore $w_0(\Delta_r)=\Delta_r$ and $w_0(\alpha_r)<0$. In particular, this shows that the parabolic subgroup $P=P_{\Delta_r}$ is self-associate.

We will show that there exists a Zariski open dense subset $N'\subset N$ such that the following Bruhat decomposition holds
$$\dot{w}_0^{-1}n=mn' \overline{n}.$$
A typical element in $ \overline{N}=\dot{w}_0N\dot{w}_0^{-1}$ is of the form $\overline{n}=\overline{n}(X_1, Y_1)=\begin{bmatrix}
    I_r & & \\
    (-1)J_{2m}' {^tX_1}J_r & I_{2m} & \\
    (-1)^rY_1 & (-1)^r X_1 & I_r \\ 
\end{bmatrix}$ for some $X_1\in\mathrm{Mat}_{r\times 2m}, Y_1\in \mathrm{Mat}_{r\times r}$. Let $m=(m_1, m_2)=\begin{bmatrix}
    m_1 & & \\
     & m_2 & \\
      & & J_r {^tm_1^{-1}} J_r\\
\end{bmatrix}\in M $ and $n'=n(X', Y')=\begin{bmatrix}
    I_r & X' & Y' \\
     & I_{2m} & (-1)^rJ_{2m}' {^tX'} J_r\\
     & & I_r\\
\end{bmatrix}\in N$.  One computes that
$$mn'\overline{n}=\begin{bmatrix}
    m_1 & & \\
    & m_2 & \\
    & & J_r {^tm_1^{-1}}J_r^{-1}
\end{bmatrix}\begin{bmatrix}
    I_r & X' & Y' \\
     & I_{2m} & (-1)^rJ_{2m}' {^tX'} J_r\\
     & & I_r\\
\end{bmatrix}\begin{bmatrix}
    I_r & & \\
   (-1)^r J_{2m}' {^tX_1}J_r & I_{2m} & \\
    (-1)^rY_1 & (-1)^r X_1 & I_r \\ 
\end{bmatrix}$$
$$=\begin{bmatrix}
    m_1+(-1)^rm_1X'J_{2m}'{^tX_1}J_r+(-1)^rm_1Y'Y_1 & m_1X'+(-1)^rm_1Y'X_1 & m_1Y' \\
(-1)^rm_2J_{2m}'{^tX_1}J_r+m_2J_{2m}' {^tX'}J_r Y_1 & m_2+m_2J_{2m}' {^tX'}J_rX_1 & (-1)^rm_2J_{2m}'{^tX'}J_r \\
    (-1)^rJ_r{^tm_1^{-1}}J_r^{-1}Y_1 & (-1)^r J_r{^tm_1^{-1}} J_r^{-1} X_1 & J_r{^tm_1^{-1}}J_r^{-1}\\
\end{bmatrix}.$$
On the other hand, let $n=n(X,Y)\in N$, we have
$$\dot{w}_0^{-1}n=\begin{bmatrix}
    & & (-1)^rI_r\\
    & I_{2m} & (-1)^rJ_{2m}'{^tX}J_r\\
    I_r & X & Y \\
\end{bmatrix}.$$ Compare both sides we obtain the following equalities:
\begin{enumerate}
    \item $$m_1+(-1)^rm_1X'J_{2m}'{^tX_1}J_r+(-1)^rm_1Y'Y_1=m_1X'+(-1)^rm_1Y'X_1=(-1)^rm_2J_{2m}'{^tX_1}J_r+m_2J_{2m}' {^tX'}J_r Y_1 =0,$$
    \item $$m_1Y'=(-1)^rI_r=J_r{^tm_1^{-1}}J_r^{-1}Y_1,$$
    \item $$m_2+m_2J_{2m}' {^tX'}J_rX_1=I_{2m},$$
    \item 
$$m_2J_{2m}'{^tX'}J_r= J_{2m}' {^tX}J_r,$$
    \item $$(-1)^r J_r{^tm_1^{-1}} J_r^{-1} X_1=X,$$
    \item $$J_r{^tm_1^{-1}}J_r^{-1}=Y.$$
\end{enumerate}
The outer automorphism $\theta_r(g):=J_r{^tg^{-1}}J_r^{-1}, g\in \mathrm{GL}_r$ defines a involution of $\mathrm{GL}_r$. Assume $\det(Y)\neq 0$, then  $$(6)\Leftrightarrow m_1=J_r{^tY^{-1}}J_r^{-1}=\theta_r(Y),$$ $$(5)\Leftrightarrow X_1=(-1)^r\theta_r(m_1)^{-1}X=(-1)^rY^{-1}X,$$
$$(2)\Leftrightarrow Y_1=(-1)^rY^{-1}, Y'=(-1)^r\theta_r(Y)^{-1}.$$ 
The second equality of (1) is equivalent to $$X'=(-1)^{r-1}\theta_r(Y^{-1})Y^{-1}X.$$
Plug this in (3), we have
$$m_2(I_{2m}+J'_{2m}{^tX'}J_rX_1)=m_2(I_{2m}-J'_{2m}{^tX}{^tY^{-1}}\theta_r({^tY^{-1}})J_r Y^{-1}X),$$
$$=m_2(I_{2m}-J'_{2m}{^tX}{^tY}^{-1}J_r X)=m_2(I_{2m}+(-1)^{r-1}\theta_{r,m}(Y^{-1}X)X)=I_{2m},$$
where $$\theta_{r,m}: \mathrm{Mat}_{r\times 2m}\longrightarrow \mathrm{Mat_{2m\times r}},$$
$$X\mapsto {^t(J_rXJ_{2m}')}={^tJ'_{2m}}{^tX}{^tJ_r}=(-1)^rJ'_{2m}{^tX}J_r.$$ 
If we further assume that $I_{2m}-J'_{2m}{^tX}{^tY}^{-1}J_r X\in\mathrm{GL}_{2m}$, which is a Zariski open dense condition on $N$, then
$$m_2=(I_{2m}-J'_{2m}{^tX}{^tY}^{-1}J_r X)^{-1}=(I_{2m}+(-1)^{r-1}\theta_{r,m}(Y^{-1}X)X)^{-1}.$$
One easily checks that other equalities in (1) are automatically satisfied for our solutions of $m_1,m_2, X', Y', X_1, Y_1.$ To check that (4) holds, since $m_2\in \mathrm{Sp}_{2m}$, i.e., ${^tm_2}J'_{2m}m_2=J'_{2m}$, (4) is equivalent to $X'm_2^{-1}=X$. So we need to check that $$(-1)^{r-1}\theta_r(Y^{-1})Y^{-1}X(I_{2m}-J'_{2m}{^tX}{^tY^{-1}}J_rX)=X,$$ hence it suffices to check that $$(-1)^{r-1}J_r{^tY}J_r^{-1}Y^{-1}+(-1)^rJ_r{^tY}J_r^{-1}Y^{-1}X J'_{2m}{^tX}{^tY^{-1}}J_r=I_r.$$ Multiplying $YJ_r{^tY^{-1}}J_r^{-1}$ on the left and $J_r^{-1}{^tY}$ on the right, one simplifies to get that $J_r{^tY}-Y{^tJ_r}+(-1)^rXJ'_{2m}{^tX}=0$, which holds automatically by the structure of $\mathrm{Sp}_{2n}$.

\section{The orbit space $Z_M^0U_M\backslash N$ and its invariant measure.} Recall that the Levi component $M$ of the standard parabolic subgroup of $G=\mathrm{Sp}_{2n}$ obtained by removing $\alpha_r$ in its Dynkin diagram is isomorphic to $\mathrm{GL}_r\times \mathrm{Sp}_{2m}$ with $r+m=n$. Let $U_M=U\cap M$, then $U_M$ acts on $N$ by conjugation. The local coefficient in these cases can be represented as the Mellin transform of certain partial Bessel functions over the orbit space $U_M\backslash N$. In order to proceed with the tools in harmonic analysis so that we could obtain a nice asymptotic expansion formula for partial Bessel functions, and finally prove the analytic stability, we need to find a Zariski open dense subset $N'\subset N$ such that the geometric quotient $U_M\backslash N'$ exists, and construct an invariant measure on the corresponding p-adic manifold over a non-Archimedean local field $F$ of characteristic 0. This is the main goal of this section.

We will use the same notations as in section 3. A simple calculation shows that $U_M$ acts on $N$ by
$$U_M\times N\longrightarrow N$$
$$(u, n(X,Y))\mapsto n(u_1Xu_2^{-1}, u_1YJ_r{^tu_1}J_r^{-1})=n(u_1Xu_2^{-1}, u_1Y \theta_r(u_1^{-1}))$$ where we identify $u=\textrm{diag}(u_1, u_2, J_r{^tu_1^{-1}}J_r^{-1})\in U_M$ with $(u_1, u_2)\in U_{\mathrm{GL}_r}\times U_{\mathrm{Sp}_{2m}}$, where $U_{\mathrm{GL}_r}$ and $U_{\mathrm{Sp}_{2m}}$ are the maximal unipotent subgroups of $\mathrm{GL}_r$ and $\mathrm{Sp}_{2m}$ respectively. 

Since $n(X,Y)\in N$, $X$ and $Y$ are related by $ J_r{^tY}-Y{^tJ_r}+(-1)^rX J'_{2m} {^tX}=0$. Set $Z=J_r{^tY}+\frac{(-1)^rXJ'_{2m}{^tX}}{2}$, thus 
$$ J_r{^tY}-Y{^tJ_r}+(-1)^rX J'_{2m} {^tX}=0\Leftrightarrow Z={^tZ}.$$
Then the action $(X,Y)\mapsto (u_1Xu_2^{-1}, u_1Y\theta_r(u_1^{-1}))$
is equivalent to 
$X\mapsto u_1Xu_2^{-1}, Z\mapsto u_1Z{^tu_1}$. The advantage of this change of variable is that now $X$ and $Z$ are independent. Denote the space of $k\times k$ symmetric matrices as $\mathrm{Sym}^k$. We have the following description on the orbit space representatives and measures for $U_M\backslash N$:

\begin{prop} \label{prop4.1}There exists a Zariski open dense subset $N'\subset N$, such that the Bruhat decomposition $\dot{w}_0^{-1}n=mn'\overline{n}$ holds, and the orbit space $U_M\backslash N'$ admits a set of representatives of the form
$$R=(R_X, R_Z)$$ where
$$R_X=\begin{cases}
\{x_{r-k, 1+k}(0\le k\le r-1), x_{i,j}(2m-r+1\le j-i, j\le 2m)\}, &  r\le m,\\
\{x_{i,j}: i+j=r+1 (j\le m+1), \mathrm{\ or\ }r-m-l\le i\le r-m+l,j=m+l+1 ( 1\le l\le r-m-1), & \\
\mathrm{\ or \ } 1\le i\le 2(r-m)-1+k, j=r+k(0\le k\le 2m-r)\}, & m<r<2m,\\
    \{x_{i,j}: i+j=r+1(j\le m+1), \mathrm{\ or \ }r-m-l\le i\le r-m+l, j=m+l+1(1\le l \le m-1)\}, & r\ge 2m,
\end{cases}$$
$R_Z=
    \mathrm{Sym}^{r}$ if $r<2m$, and $R_Z=
    \{(z_{ij})\in\mathrm{Sym}^r| z_{ij}=0 \text{ if } i, j\leq r-2m \text{ and }i\neq j\}$ if $r\ge 2m$.
    
Set $d_X=\dim R_X$, and $d_Z=\dim R_Z$.
Then 
$$d_X=\begin{cases}
      \frac{r(r+1)}{2}, & \mathrm{\ if \ } r\le m,\\
      2rm-\frac{r(r-1)}{2}-m^2, & \mathrm{ \ if \ } m<r<2m,\\
      m(m+1), & \mathrm{\ if \ } r\ge 2m,
  \end{cases}, 
\mathrm{\ and\ } 
d_Z=\begin{cases}
      \frac{r(r+1)}{2} & \mathrm{\ if \ } r< 2m,\\
      (2m+1)(r-m) & \mathrm{\ if \ } r\ge 2m.
  \end{cases}$$ Note that when $r<m$, the action has a non-trivial but base point free stabilizer isomorphic to $U_{\mathrm{Sp}_{2(m-r)}}$ and when $r\geq m$ the stabilizer is always trivial. Moreover, the corresponding invariant measure $d\mu$ on $R$ is given by $d\mu=d\mu_X\wedge d\mu_Z$, where
  $$d\mu_X=\begin{cases}

\vert x_{r,1}^{r+2m-2}x_{r-1, 2}^{r+2m-5}x_{r-2, 3}^{r+2m-8}\cdots x_{1,r}^{2m-2r+1}\vert \prod dx_{ij}, & r\le m,\\

\vert x_{r,1}^{r+2m-2}x_{r-1, 2}^{r+2m-5}x_{r-2, 3}^{r+2m-8}\cdots x_{r-m+1,m}^{r-m+1}x_{r-m, m+1}^{r-m-1} x_{r-m-1,m+2}^{r-m-2}\cdots x_{2,r-1}\vert \prod dx_{ij}, & m<r<2m,\\

\vert x_{r,1}^{r+2m-2}x_{r-1, 2}^{r+2m-5}x_{r-2, 3}^{r+2m-8}\cdots x_{r-m+1,m}^{r-m+1}x_{r-m, m+1}^{r-m-1} x_{r-m-1,m+2}^{r-m-2}\cdots x_{r-2m+1, 2m}^{r-2m}\vert \prod dx_{ij}, & r\ge 2m,\\
\end{cases}$$
 and $$d\mu_Z=\begin{cases}
     \prod_{i,j}dz_{i,j} & \mathrm{\ if \ } r\leq 2m,\\
     \vert z_{r-2m, r-2m}^{r-2m-1}z_{r-2m-1, r-2m-1}^{r-2m-2}\cdots z_{2,2}\vert\prod dz_{ij} & \mathrm{\ if \ } r>2m,
 \end{cases}$$ where the product runs over all $(i,j)$'s such that $x_{i.j}$ and $z_{i,j}$ are non-zero in each case. In addition, if we use $(X,Y)$ instead of $(X,Z)$ to parameterize the orbit space, the corresponding invariant measures are related by $$d\mu_X\wedge d\mu_Y=d\mu_X\wedge (d\mu_Z\cdot J_r).$$ 
\end{prop}

\begin{proof}
 Based on the arguments before the statement of the proposition, it suffices to study the action $X\mapsto u_1Xu_2^{-1}, Z\mapsto u_1Z{^tu_1}$, by induction of the size of $X$ and $Z$. If $r=1$, then $u_1=1$, the action degenerates as $(X,Z)\mapsto (Xu_2^{-1}, Z)$. Assume $x_{1,1}\neq 0$, then it is easy to see that the orbit space representative can be given as $((x_{1,1},0,\cdots, 0),Z)\in \mathbb{A}^{2m}\times \mathrm{Sym}^r$. 

Suppose $r>1$. We study the action on $X$ first. Write $$u_1=\begin{bmatrix}
    v_1 & {^t\delta_1} \\
    & 1 \\
\end{bmatrix}, u_2=\begin{bmatrix}
    1 & \delta_2\\
     & v'_2\\
\end{bmatrix}, X=\begin{bmatrix}
    {^t\alpha} & X'_1 \\
    x_{r,1} & \beta\\
\end{bmatrix}$$ with $\alpha, \delta_1\in \mathbb{A}^{r-1}$, $\beta, \delta_2\in \mathbb{A}^{2m-1}$, $v_1\in U_{\mathrm{GL}_{r-1}},v'_2\in U_{\mathrm{GL}_{2m-1}}$ , and $X'_1\in \mathrm{Mat}_{(r-1)\times(2m-1)}$. Then
$$u_1 Xu_2^{-1}=\begin{bmatrix}
    v_1 & {^t\delta_1} \\
    & 1 \\
\end{bmatrix}\begin{bmatrix}
    {^t\alpha} & X'_1 \\
    x_{r,1} & \beta\\
\end{bmatrix}\begin{bmatrix}
    1 & -\delta_2 (v'_2)^{-1}\\
     & (v'_2)^{-1}\\
\end{bmatrix}=\begin{bmatrix}
    v_1{^t\alpha}+{^t\delta_1}x_{r,r+1} & v_1 X'_1+{^t\delta_1}\beta \\
    x_{r,1} & \beta\\
\end{bmatrix}\begin{bmatrix}
    1  &  -\delta_2 (v'_2)^{-1}\\
     & (v'_2)^{-1}\\
\end{bmatrix}$$
$$=\begin{bmatrix}
    v_1{^t\alpha}+{^t}\delta_1x_{r,1} & -v_1{^t\alpha}\delta_2 (v'_2)^{-1}-x_{r,1}{^t\delta_1}\delta_2 (v'_2)^{-1}+v_1X'_1(v'_2)^{-1}+{^t\delta_1}\beta (v'_2)^{-1}\\
    x_{r,1} & -x_{r,1}\delta_2
(v'_2)^{-1}+\beta (v'_2)^{-1}\end{bmatrix}.$$
Assume $x_{r,1}\neq 0$, choose $\delta_1$ and $\delta_2$ such that $v_1{^t\alpha}+{^t}\delta_1x_{r,1}=-x_{r,1}\delta_2
(v'_2)^{-1}+\beta (v'_2)^{-1}=0$, i.e. $\delta_1=-\frac{\alpha{^tv_1}}{x_{r,1}}$, $\delta_2=\frac{\beta}{x_{r,1}}$, then
$-v_1{^t\alpha}\delta_2 (v'_2)^{-1}-x_{r,1}{^t\delta_1}\delta_2 (v'_2)^{-1}+v_1X'_1(v'_2)^{-1}+{^t\delta_1}\beta (v'_2)^{-1}=v_1X'_1(v'_2)^{-1}-\frac{v_1{^t\alpha}\beta (v'_2)^{-1}}{x_{r,1}}$.
Let $X''_1=X'_1-\frac{{^t\alpha}\beta}{x_{r,1}}$, then $$v_1X'_1(v'_2)^{-1}-\frac{v_1{^t\alpha}\beta (v'_2)^{-1}}{x_{r,1}}=v_1X''_1 (v'_2)^{-1}.$$

We have $X''_1\in \mathrm{Mat}_{(r-1)\times (2m-1)}$, $v_1\in U_{\mathrm{GL}_r}$ and $v'_2\in U_{\mathrm{GL}_{2m-1}}$. Now if $m=1$ then $v_2'=1$, and the action degenerates as $X''_1\mapsto v_1 X''_1$ with $X''_1\in \mathrm{Mat}_{(r-1)\times 1}$. So it is clear that if we assume that the $(r-1,2)$-th entry of $X_1''$ is non-zero, we obtain an orbit space representative of the action on $X$ given by 
$R_X=\begin{bmatrix}
    0 & 0\\
    \vdots & \vdots\\
    0 & x_{r-1,2}\\
    x_{r,1} & 0\\
\end{bmatrix}$. Therefore, we assume that $m>1$ from now on.

Write $u_2=\begin{bmatrix}
    1 & \delta_2 \\
     & v_2'
\end{bmatrix}=\begin{bmatrix}
    1 & \gamma_2 & x\\
    & v_2 & {^t\gamma_2'} \\
    & & 1\\
\end{bmatrix},$ where $\delta_2=(\gamma_2, x)$ and $v'_2=\begin{bmatrix}
    v_2 & {^t\gamma_2'} \\
    & 1\\
\end{bmatrix}$, with $\gamma_2,\gamma_2'\in \mathbb{A}^{2m-2}$. Moreover, by the structure of $\mathrm{Sp}_{2m}$, a simple calculation shows that $\gamma_2'=-\gamma_2v_2^{-1}J_{2m-2}'$, $\gamma_2'{^tJ'_{2m-2}}{^t\gamma'_2}=0$, $x$ is free, and ${^tv_2}J_{2m-2}'v_2=J_{2m-2}'$, hence $v_2\in U_{\mathrm{Sp}_{2m-2}}$. Note that since $\delta_2$ is determined by the above process, so is $\gamma_2$.
We also write $X''_1=\begin{bmatrix}
    X_1 & {^t\gamma}\\
\end{bmatrix}$ with $X_1\in \mathrm{Mat}_{(r-1)\times (2m-2)}$ and $\gamma\in \mathbb{A}^{r-1}$.
Therefore we can write
$$v_1 X_1'(v_2')^{-1}=
    v_1 \begin{bmatrix}
        X_1 & {^t\gamma}\\
        \end{bmatrix}\begin{bmatrix}
            v_2^{-1} & -v_2^{-1}{^t\gamma_2'}\\
            & 1\\
        \end{bmatrix}
=\begin{bmatrix}
    v_1 X_1 v_2^{-1} & v_1 X_1 v_2^{-1}{^t\gamma'_2}+v_1{^t\gamma}\\
\end{bmatrix}$$
$$=\begin{bmatrix}
    v_1X_1v_2^{-1} & v_1X_1 v_2^{-1}(-{^tJ'_{2m-2}}{^tv_2^{-1}}{^t\gamma_2})+v_1{^t\gamma}\\
\end{bmatrix}=\begin{bmatrix}
    v_1X_1v_2^{-1} & -v_1X_1({^tv_2}J'_{2m-2}v_2)^{-1}{^t\gamma_2}+ v_1{^t\gamma}\\
\end{bmatrix}$$
$$=\begin{bmatrix}
    v_1X_1v_2^{-1} & v_1X_1 J'_{2m-2} {^t\gamma_2} +v_1{^t\gamma}\\
\end{bmatrix}=\begin{bmatrix}
     v_1X_1v_2^{-1} & (v_1X_1 v_2^{-1})v_2J'_{2m-2} {^t\gamma_2} +v_1{^t\gamma}\\
\end{bmatrix}.$$
From this observation we see that it suffices to study the action $$(U_{\mathrm{GL}_{r-1}}\times U_{\mathrm{Sp}_{2m-2}})\times (\mathrm{Mat}_{(r-1)\times (2m-2)}\times \mathbb{A}^{r-1})\longrightarrow \mathrm{Mat}_{(r-1)\times(2m-2)}\times \mathbb{A}^{r-1}$$
$$((v_1,v_2), (X_1,\gamma))\mapsto (v_1X_1v_2^{-1}, v_1{^t\gamma}).$$

Set $X_0:=X$, $r'=r-1$ and $m'=m-1$, replace $X$ by $X_1$, $u_1$ by $v_1$, $u_2$ by $v_2$, $r$ by $r'$, and $m$ by $m'$. Continue with the above process, we can construct $X_i$'s ($i\ge 0$) inductively. Since the relative size of $r$ and $2m$ will affect the inductive process, we will have to discuss three cases separately:

\textbf{Case (1):} $r\le m$. Then $r'\le m'$ for all $r'$ and $m'$ in the inductive process. Note that in this case all the entries of $v_1$, hence of $u_1$, are chosen in the inductive process. Consequently we do not have free variables in $v_1$ when considering the action $\gamma\mapsto v_1{^t\gamma}$. We further assume that all the entries of $v_1{^t\gamma}$ are non-zero. We eventually obtain that $X_{r-1}=(\tilde{x}_{1,r}, *, \cdots, *)$ is a vector in $\mathbb{A}^{2(m-r)+2}$ and we are left to consider the right-action of $U_{\mathrm{Sp}_{2(m-r)+2}}$. Write the action as $X_{r-1}\mapsto X_{r-1}u^{-1}$ with $u=\begin{bmatrix}
    1 &  \delta \\
     & v'\\
\end{bmatrix}=\begin{bmatrix}
    1 & \gamma & a\\
    & v & {^t\gamma'}\\
    & & 1\\
\end{bmatrix}\in U_{Sp_{2(m-r)+2}}$, where $\delta=(\gamma,a)\in \mathbb{A}^{2(m-r)+1}$, $\gamma,\gamma'\in \mathbb{A}^{2(m-r)}$, and $v$ is a unipotent matrix of size $2(m-r)$. Then a similar calculation as above shows that $\gamma'=-\gamma v^{-1}J'_{2(r-m)}$, $a$ is free, and $v\in U_{\mathrm{Sp}_{2(m-r)}}$. Write $X_{r-1}=(\tilde{x}_{r,1}, \alpha')$ with $\alpha'\in \mathbb{A}^{2(m-r)+1}$. By assuming $\tilde{x}_{1,r}\neq 0$, we choose $\delta=(\gamma, a)$ such that
$\tilde{x}_{r,1}\delta+\alpha'v=0$, then $X_{r-1}u^{-1}=(\tilde{x}_{r,1},0,\cdots, 0)$. So we obtain our orbit space representative in this case as
$$R_X=\begin{bmatrix}
    0 & 0 & \cdots & 0 & x_{1,r} & 0 \ \cdots \ 0 & x_{1,r+2(m-r)+2} & x_{1,r+2(m-r)+3} & \cdots & x_{1,2m}\\
    0 & 0 & \cdots & x_{2,r-1} & 0 & 0 \ \cdots \ 0 & 0 & x_{2, r+2(m-r)+3} & \cdots  & x_{2, 2m} \\
    \vdots & \vdots & \iddots &\vdots & \vdots & 0\ \cdots \ 0 & \vdots & \vdots & \ddots &\vdots \\
    0 & x_{r-1,2} &  \cdots & 0 & 0 & 0\ \cdots \ 0 & 0 & 0 & \cdots & x_{r-1, 2m}\\
    x_{r,1} & 0 & \cdots & 0 & 0 & 0\ \cdots \ 0 & 0 & 0 & \cdots & 0
\end{bmatrix},$$
where there are $2(m-r)+1$ zero columns in the middle. From the last step above, we also see that in this case the action $X\mapsto u_1Xu_2^{-1}$ has stablizer $U_{\mathrm{Sp}_{2(m-r)}}$ given by $v$, provided that $m>r$. Moreover, one computes easily that 
$\dim R_X=r+ \frac{r(r-1)}{2}=\frac{r(r+1)}{2}$.

\textbf{Case (2):} $r\ge 2m$. Then $r'\ge 2m'$ for all $r'$ and $m'$ in the inductive process. Perform the inductive process to the $m$-th step, we obtain that $X_m=\begin{bmatrix}
    * & *\\
    \vdots & \vdots \\
    * & x_{r-m,m+1}\\
    x_{r-m+1.m} & *\\
\end{bmatrix}$, now $v_1=\begin{bmatrix}
    v_1' & {^t\delta_1}\\
    & 1\\
\end{bmatrix}\in U_{\mathrm{GL}_{r-m}}$, and $v_2=\begin{bmatrix}
    1 & x \\
    & 1\\
\end{bmatrix}$. Choose $\delta_1$ and $x$ accordingly we can make $X_m$ to be of the form $\begin{bmatrix}
    0 & *\\
    \vdots &\vdots \\
    0 & x_{r-m,m+1}\\
    x_{r-m+1,m} & 0\\
\end{bmatrix}$. From the $(m+1)$-th step on the action degenerates as the left action of $U_{\mathrm{GL}_{r-m-1}}$ on $(X, {^t\gamma})$. By determining the last column of $v_1$ each time to make the first column but the last entry of $X$ zero, we eventually obtain a vector of the form ${^t(*,\cdots, * , \tilde{x}_{r-2m+1, 2m})}\in\mathbb{A}^{r-2m+1}$. During this process, all the entries of $u_2$ together with all the last $2m-1$ columns of $u_1$ are determined, and we are left to consider the left-action of $U_{\mathrm{GL}_{r-2m+1}}$ on $\mathbb{A}^{2m-r+1}$. Assume $\tilde{x}_{r-2m+1,2m}\neq 0$, we pick the orbit representative in the last step as $^t(0,\cdots, 0, \tilde{x}_{r-2m+1, 2m})$. As a result, we obtain our orbit space representative in this case as

$$R_X=\begin{bmatrix}
    0 & 0 & \cdots & 0 & 0 & 0 & 0 & 0 & 0\\
      & & \vdots & \vdots & \vdots & \vdots & \vdots & \vdots & \vdots \\
    0 & 0 & \cdots & 0 & 0 & 0 & 0 & 0 & 0\\
    0 & 0 & \cdots & 0 & 0 & 0 & 0 & 0 & x_{r-2m+1, 2m}\\
    0 & 0 & \cdots & 0 & 0 & 0 & 0 & x_{r-2m+2,2m-1} & *\\
     & & \vdots & 0 & 0 & 0 & \iddots & \vdots & \vdots\\
     0 & 0 & \cdots & 0 & 0 & x_{r-m-1, m+2}& \cdots & * & *\\ 
     0 & 0 & \cdots & 0 & x_{r-m,m+1} & x_{r-m, m+2} & \cdots & * & *\\ 
     0 & 0 & \cdots & x_{r-m+1,m} & 0 & x_{r-m+1,m+2} & \cdots & * & *\\ 
     & & \iddots & \vdots & \vdots & \vdots & \vdots & \vdots & \vdots \\
     0 & 0 & \cdots & 0 & 0 & 0 & 0 & x_{r-2,2m-1} & *\\
     0 & x_{r-1,2} &\cdots & 0 & 0 & 0 & 0 & 0 & x_{r-1,2m}\\
     x_{r,1} & 0 & \cdots & 0 & 0 & 0 & 0 &  0 & 0,\\ 
\end{bmatrix}$$ where there are $r-2m$ zero rows on the top.
Observe that in this case the action $X\mapsto u_1X u_2^{-1}$ has a stablizer $U_{\mathrm{GL}_{r-2m}}$ given by the upper-left unipotent submatrix of $u_1$ of size $r-2m$, if $r>2m$. We also have $\dim R_X=m+(m+\frac{m(m-1)}{2}\cdot 2)=m(m+1)$.

\textbf{Case (3)}: $m<r<2m$. In this case it is possible that there exists some $k\ge 0$ such that $r'=r-k\ge 2(m-k)=2m'$, i.e., $2m-r\le k\le m$. If this happens, continue with our algorithm in case (2). Specifically, note that when $k=2m-r>0$, $r'=2m'=2(r-m)$, we can conclude as in the last step of case (2) that the action $X\mapsto u_1Xu_2^{-1}$ has trivial stablizer.  If for all $r'$ and $ m'$ in the inductive process, we have $r'<2m'$, we perform a similar inductive argument as in case (1). Since $r'<2m'$ for all $r',m'$ in the inductive process, by a similar argument as in case (1), we have exhausted all the possible choices of entries in $u_1$. On the other hand, since $r'>m'$, when $r'$ decreases to $1$, so does $m'$, hence we also exhaust all the possible choices of entries in $u_2$. Consequently, the action $X\mapsto u_1Xu_2^{-1}$ has trivial stabilizer. Eventually one obtains the orbit space representative in this case as:
$$R_X=\begin{bmatrix}
    0 & \cdots & 0 & 0 & 0 & \cdots\ \ \ 0 & x_{1, r} & x_{1,r+1} & \cdots & x_{1,2m}\\
    0  & \cdots & 0 & 0 & 0 & \cdots\ \ \ x_{2,r-1} & x_{2,r} & x_{2, r+1} &\cdots & x_{2,2m}\\
    \vdots & & \vdots & \vdots & \vdots & \iddots\ \ \ \vdots & \vdots & \vdots & \cdots & \vdots\\
     0  & \cdots & 0 & 0 & x_{r-m-1,m+2} & \cdots \ \ \ * & * & * & \cdots & *\\ 
     0 & \cdots & 0 & x_{r-m,m+1} & x_{r-m,m+2} & \cdots \ \ \ * & * & * &\cdots & *\\ 
     0 & \cdots & x_{r-m+1,m} & 0 & x_{r-m+1,m+2} & \cdots \ \ \ * & * & * & \cdots & *\\
      \vdots &  & \vdots & \vdots & \vdots & \cdots\ \ \ \vdots & \vdots & \vdots & \cdots & \vdots\\
      0 &\cdots & 0 & 0 & 0 & \cdots \ \ \ 0 & x_{2(r-m)-1,r} & * & \cdots & * \\
     0 &\cdots & 0 & 0 & 0 & \cdots \ \ \ 0 & 0 & x_{2(r-m),r+1} & \cdots & * \\
     0 & \cdots & 0  & 0 & 0 & \cdots\ \ \ 0 & 0 & 0 & \cdots & *\\
     \vdots & \iddots & \vdots & \vdots & \vdots & \cdots\ \ \  \vdots & \vdots & \vdots & \ddots & \vdots\\
     0 & \cdots & 0 & 0 & 0 & \cdots\ \ \  0 & 0 & 0 & \cdots & x_{r-1,2m}\\
    x_{r,1} & \cdots & 0 & 0 & 0 & \cdots\ \ \  0 & 0 &  0 & \cdots & 0\\ 
\end{bmatrix}.$$
We have $\dim R_X=m+(r-m)^2+\frac{(3r-2m-1)(2m-r)}{2}=2rm-\frac{r(r-1)}{2}-m^2$.

Note that only in case (2) we also have to consider the action of $U_{\mathrm{GL}_{r-2m}}$ on $Z\in \mathrm{Sym}^r$. Write $Z=\begin{bmatrix}
    Z_1 & h \\
    {^th} & Z_2\\
\end{bmatrix}$, with $Z_1\in \mathrm{Sym}^{r-2m}, Z_2\in\mathrm{Sym}^{2m}$, $h\in\mathrm{Mat}_{(r-2m)\times 2m}$, and $u_1=\begin{bmatrix}
    u'_1 & w \\
    & u_1''\\
\end{bmatrix}$ with $u_1'\in U_{\mathrm{GL}_{r-2m}}$, $u_1''\in U_{\mathrm{GL}_{2m}}$ and $w\in \mathrm{Mat}_{(r-2m)\times 2m}$. A simple calculation gives
$$u_1Z{^tu_1}=\begin{bmatrix}
    u_1'Z_1{^tu'_1}+w{^th}{^tu_1'}+(u'_1h+wZ_2){^tw} & u_1'h{^tu_1''}+wZ_2 {^tu_1''}\\
    u_1''{^th}{^tu'_1}+u_1''Z_2{^tw} & u_1''Z_2 {^tu_1''}\\
\end{bmatrix}$$ Note that both $w$ and $u_1''$ are fixed during the inductive process of the action $X\mapsto u_1Xu_2^{-1}$, so it suffices to consider the action $(u_1', (Z_1, h))\mapsto (u_1'Z_1{^tu_1'}, u_1'h)$. We further write $u_1'=\begin{bmatrix}
    \tilde{u}_1 & {^t\alpha_1}\\
    & 1
\end{bmatrix}$, and $Z_1=\begin{bmatrix}
    Z_1' & {^t\delta'_1}\\
    \delta'_1 & z\\
\end{bmatrix}$. Then
$$u_1'Z_1{^tu'_1}=\begin{bmatrix}
    \tilde{u}_1 Z_1'{^t\tilde{u}_1}+{^t\alpha_1}\delta'_1{^t\tilde{u}_1}+(\tilde{u}_1{^t\delta'_1}+{^t\alpha_1z})\alpha_1 & \tilde{u}_1{^t\delta'_1}+{^t\alpha_1}z \\
    \delta'_1{^t\tilde{u}_1}+z\alpha_1 & z\\
\end{bmatrix}$$ Assume $z\neq 0$ we can choose $\alpha_1$ such that $\tilde{u}_1{^t\delta'_1}+{^t\alpha_1}z=0$, then the matrix becomes $\begin{bmatrix}
    \tilde{u}_1 (Z_1'-\frac{^t\delta'_1\delta'_1}{2}){^t\tilde{u}_1} & 0 \\
    0 & z\\
\end{bmatrix}$. Replace $Z_1$ by $Z_1'-\frac{^t\delta'_1\delta'_1}{2}$, $u_1'$ by $\tilde{u}_1$ and continue with the above process, we exhaust all possible choices of entries in $u_1$. Hence there is no need to consider $h\mapsto u_1' h$. Assume that after our choices of all entries in $u_1'$, we have $u_1'h\neq 0$. We obtain a orbit space representative given by $$R_Z=\begin{bmatrix}
    z_{11} & & &z_{1,r-2m+1} & * & \cdots & * & z_{1,r}\\
    & \ddots & & \vdots & \vdots &\cdots & \vdots & \vdots\\
    & & z_{r-2m,r-2m} & z_{r-2m, r-2m+1} & * & \cdots & * & *\\ 
    z_{1,r-2m+1} & \cdots & z_{r-2m, r-2m+1} & z_{r-2m+1,r-2m+1} & * & \cdots & * & *\\
     * &\cdots & *  & * & * & \cdots & * & * \\ 
     \vdots &\cdots &\vdots & \vdots & \vdots & \ddots & \vdots & \vdots\\
     * &\cdots & *  & * & * & \cdots & z_{r-1, r-1} & *\\
     z_{1,r} &\cdots & *& *  & * &\cdots  &  * & z_{r,r}\\
\end{bmatrix}.$$
It follows that $\dim R_Z=r-2m+(r-2m)2m+ \frac{(2m+1)2m}{2}=(2m+1)(r-m)$.

Next, we study the invariant measure on our chosen orbit space representatives. We will first work on $R_X$. Denote by $U_0$ the stabilizer in $U_M\simeq U_{\mathrm{GL}_r}\times U_{\mathrm{Sp}_{2m}}$, which is independent of the base point as can be observed from our previous discussion. The action is given by
$$(U_{\mathrm{GL}_r}\times U_{\mathrm{Sp}_{2m}})/ U_0\times R_X\longrightarrow \mathrm{Mat}_{r\times 2m}$$
$$((u_1, u_2), X_0)\mapsto u_1X_0u_2^{-1}$$
Without loss of generality we assume $m>2$. Again we write $u_1=\begin{bmatrix}
    v_1 & {^t\delta_1}\\
    & 1\\
\end{bmatrix}$, $u_2=\begin{bmatrix}
    1 & \delta_2\\
    & v_2'\\
\end{bmatrix}=\begin{bmatrix}
    1 & \gamma_2 & x\\
    & v_2 & {^t\gamma'_2}\\
    & & 1\\
\end{bmatrix}$, and $X_0=\begin{bmatrix}
    0 & X'_1\\
    x_{r,1} & 0\\
\end{bmatrix}=\begin{bmatrix}
    0 & X_1 & {^t\gamma}\\
    x_{r,1} & 0 & 0\\
\end{bmatrix}$. By the previous argument we have $\gamma'_2=-\gamma_2v_2^{-1}J_{2m-2}'$.
Then $$u_1X_0u_2^{-1}=\begin{bmatrix}
    {^t\delta_1}x_{r,1} & -x_{r,1}{^t\delta_1}\delta_2(v'_2)^{-1}+ v_1 X'_1(v'_2)^{-1}\\
    x_{r,1} & -x_{r,1}\delta_2(v'_2)^{-1}\\
\end{bmatrix}$$$$=\begin{bmatrix}
g    {^t\delta_1}x_{r,1} & -x_{r,1}{^t\delta_1}\gamma_2v_2^{-1}+v_1 X_1v_2^{-1} & x_{r,1}{^t\delta_1}\gamma_2J'_{2m-2}{^t\gamma_2}-x_{r,1}{^t\delta_1}x-v_1X_1 J'_{2m-2}{^t\gamma_2}+v_1{^t\gamma}\\
    x_{r,1} & -x_{r,1}\gamma_2v_2^{-1} & x_{r,1}\gamma_2J'_{2m-2}{^t\gamma_2}-x_{r,1}x
\end{bmatrix}.$$
We write this map as
$$(\delta_1, \gamma_2, x, x_{r,1}, \gamma, v_1, v_2, X_1)$$$$\mapsto (x_{r,1}\delta_1, -x_{r,1}\gamma_2v_2^{-1}, x_{r,1}\gamma_2J'_{2m-2}{^t\gamma_2}-x_{r,1}x, x_{r,1},$$$$ x_{r,1}{^t\delta_1}\gamma_2J'_{2m-2}{^t\gamma_2}-x_{r,1}{^t\delta_1}x-v_1X_1J'_{2m-2}{^t\gamma_2}+v_1{^t\gamma}, -x_{r,1}{^t\delta_1\gamma_2v_2^{-1}+v_1X_1v_2^{-1}}).$$

The Jacobian matrix of this map is given by
$$\begin{bmatrix}
    x_{r,1}I_{r-1} & 0 & 0 & \delta_1 & 0 & 0 & 0 & 0\\
    0 & -x_{r,1}v_2^{-1} & 0 & \gamma_2v_2^{-1} & 0 & 0 & \frac{\partial(-x_{r,1}\gamma_2v_2^{-1})}{\partial v_2} & 0\\
    0 & * & -x_{r,1} & * & 0 & 0 & 0 & 0\\
    0 & 0 & 0 & 1 & 0 & 0 & 0 & 0\\
    * & * & * & * & \frac{\partial (v_1{^t\gamma})}{\partial\gamma} & * & 0 & *\\
    * & * & 0 & * & 0 & \frac{\partial(v_1 X_1v_2^{-1})}{\partial v_1}& \frac{\partial(-x_{r,1}{^t\delta_1}\gamma_2v_2^{-1})}{v_2}+\frac{\partial(v_1 X_1v_2^{-1})}{v_2} & \frac{\partial (v_1X_1v_2^{-1})}{X_1}\\
\end{bmatrix}.$$ Multiplying the second row by $-\delta_1$ and adding it to the last row, one computes that the determinant of the Jacobian matrix is equal to
$$\vert x_{r_1}\vert^{r-2m-2} \cdot \vert \frac{\partial(v_1{^t\gamma})}{\partial \gamma}\vert\cdot \vert\det\begin{bmatrix}
   \frac{\partial(v_1 X_1v_2^{-1})}{\partial v_1} & \frac{\partial(v_1 X_1v_2^{-1})}{\partial v_2} & \frac{\partial (v_1X_1v_2^{-1})}{\partial X_1}\\
\end{bmatrix}\vert$$ From the form of $R_X$ we observe that
$$\vert \frac{\partial(v_1{^t\gamma})}{\partial \gamma}\vert=\begin{cases}
    1, & \mathrm{\  if\ } r\le m,\\
    \vert x_{r-2m+1, 2m}\vert^{r-2m}, &\mathrm{\ if \ } r\ge 2m.\\
\end{cases}$$ 
   And the last term $ \vert\det\begin{bmatrix}
   \frac{\partial(v_1 X_1v_2^{-1})}{\partial v_1} & \frac{\partial(v_1 X_1v_2^{-1})}{\partial v_2} & \frac{\partial (v_1X_1v_2^{-1})}{\partial X_1}\\
\end{bmatrix}\vert$ gives the Jacobian of the same type of action but with rank drop by 1. Proceed by induction we obtain the invariant measure on $R_X$. When $r<m<2m$, $\frac{\partial(v_1{^t\gamma})}{\partial \gamma}=1$ until our induction procedure goes to the $(2m-r+1)$-th step, where $\gamma=(x_{2,r-1},*,\cdots,*)$, in which case we have $\vert\frac{\partial (v_1{^t\gamma})}{\partial \gamma}\vert=\vert x_{2,r-1}\vert$. Then the induction continues as case (2). Consequently, we obtain the invariant measure $d\mu_X$ on $R_X$ as
$$d\mu_X=\begin{cases}

\vert x_{r,1}^{r+2m-2}x_{r-1, 2}^{r+2m-5}x_{r-2, 3}^{r+2m-8}\cdots x_{1,r}^{2m-2r+1}\vert \prod dx_{ij}, & r\le m,\\

\vert x_{r,1}^{r+2m-2}x_{r-1, 2}^{r+2m-5}x_{r-2, 3}^{r+2m-8}\cdots x_{r-m+1,m}^{r-m+1}x_{r-m, m+1}^{r-m-1} x_{r-m-1,m+2}^{r-m-2}\cdots x_{2,r-1}\vert \prod dx_{ij}, & m<r<2m,\\

\vert x_{r,1}^{r+2m-2}x_{r-1, 2}^{r+2m-5}x_{r-2, 3}^{r+2m-8}\cdots x_{r-m+1,m}^{r-m+1}x_{r-m, m+1}^{r-m-1} x_{r-m-1,m+2}^{r-m-2}\cdots x_{r-2m+1, 2m}^{r-2m}\vert \prod dx_{ij}. & r\ge 2m,\\
\end{cases}$$
where the product runs over all $(i,j)$ such that $x_{i,j}\neq 0$.

   Recall that when $r\ge 2m$, we will also need to consider the action of $U_{\mathrm{GL}_{r-2m}}$ on $\mathrm{Sym}^r$. The corresponding invariant measure on $R_Z$ can also be proved inductively. Write $u_1=\begin{bmatrix}
       u'_1 & w\\
       & u_1''\\
   \end{bmatrix}$, where $u_1'\in U_{\mathrm{GL}_{r-2m}}$, $u_1''\in U_{\mathrm{GL}_{2m}}$, $w\in\mathrm{Mat}_{(r-2m)\times 2m }$, and $Z_0=\begin{bmatrix}
       Z_1 & h\\
       {^th} & Z_2\\
   \end{bmatrix}$ with $Z_1=\mathrm{diag}\{z_{1,1},\cdots, z_{r-2m, r-2m}\}$, $Z_2\in \mathrm{Sym}^{2m}$, $h\in\mathrm{Mat}_{(r-2m)\times 2m}$. As we see from the proof of $R_Z$ above, it suffices to only consider the action $(U_1', Z_1) \mapsto u_1' Z_1{^tu_1'}$. Without loss of generality, we assume $r-2m>0$. Then we can write $u_1'=\begin{bmatrix}
       \tilde{u}_1 & \alpha_1\\
       & 1\\
   \end{bmatrix}$, with $\tilde{u}_1\in U_{\mathrm{GL}_{r-2m-1}}$, $\alpha_1\in \mathbb{A}^{r-2m-1}$, and $Z_1=\begin{bmatrix}
       Z_1' & \\
       & z_{r-2m,r-2m}\\
   \end{bmatrix}$ where $Z_1'=\mathrm{diag}\{z_{1,1},\cdots, z_{r-2m-1,r-2m-1}\}$. Then the map $(u_1', Z_1)\mapsto u_1'Z_1{^tu_1'}$ can be written as $$(z_{r-2m,r-2m}, \alpha_1, \tilde{u}_1, Z_1')\mapsto (z_{r-2m,r-2m}, z_{r-2m,r-2m}\alpha_1, \tilde{u}_1Z_1'{^t\tilde{u}_1}+{^t\alpha_1}z_{r-2m,r-2m}\alpha_1),$$ whose Jacobian matrix if of the form
   $$\begin{bmatrix}
       1 & 0 & 0 & 0\\
       0 & z_{r-2m,r-2m}I_{r-2m-1} & 0 & 0\\
       * & * & \frac{\partial (\tilde{u}_1Z_1'{^t\tilde{u}_1})}{\partial \tilde{u}_1} & \frac{\partial (\tilde{u}_1Z_1'{^t\tilde{u}_1})}{\partial Z'_1}
   \end{bmatrix}.$$ So the absolute value of its determinant is equal to 
   $\vert x_{r-2m,r-2m}^{r-2m-1}\vert\cdot \vert\det\begin{bmatrix}
       \frac{\partial (\tilde{u}_1Z_1'{^t\tilde{u}_1})}{\partial \tilde{u}_1} & \frac{\partial (\tilde{u}_1Z_1'{^t\tilde{u}_1})}{\partial Z'_1}
   \end{bmatrix}\vert$, where the second term is the absolute value of the Jacobian matrix of the same type of action with rank drop by 1. Hence we can proceed by induction. As a result, we obtain the invariant measure $d\mu_Z$ on $R_Z$ in this case given by
   $$d\mu_Z=\vert z_{r-2m, r-2m}^{r-2m-1}z_{r-2m-1, r-2m-1}^{r-2m-2}\cdots z_{2,2}\vert\prod dz_{ij}$$ where the product runs over all $(i,j)$ such that $z_{i,j}\neq 0$. For the other two cases, $d\mu_Z=\prod_{i,j}dz_{i,j}$ where the product runs over all $(i,j)$.
   
Finally, since $Z=J_r{^tY}+\frac{(-1)^rXJ'_{2m}{^tX}}{2}$, and ${^tZ}=Z$, we see that
$$Y={^t(J_r^{-1}(Z-\frac{(-1)^rXJ'_{2m}{^tX}}{2}))}={^tZJ_r}+\frac{(-1)^rXJ'_{2m}{^tX}J_r}{2}=ZJ_r+\frac{X\theta_{r,m}(X)}{2}.$$ 
Let $R_Y=R_ZJ_r+\frac{R_X\theta_{r,m}(R_X)}{2}$. Since the map $(X,Y)\mapsto (X,Z)$ is bijective and $U_M$-equivariant, $(R_X,R_Y)$ gives a set of orbit space representatives for the action $(X,Y)\mapsto (u_1Xu_2^{-1}, u_1Y\theta_r(u_1^{-1}))$. If we denote the corresponding measure on $R_Y$ as $d\mu_Y$, since $\theta_{r,m}$ is linear, we have
$$d\mu_Y=d\mu_Z\cdot J_r+\frac{d\mu_X\cdot \theta_{r,m}(X)+X\cdot \theta_{r,m}(d\mu_X)}{2}$$ Therefore
$$d\mu_X\wedge d\mu_Y=d\mu_X\wedge (d\mu_Z\cdot J_r).$$
This shows that we replace integrals over $R_X\times R_Y$ by the ones over $R_X\times R_Z$, without changing the corresponding invariant Haar measures.
\end{proof}

In order to apply Theorem 6.2 of \cite{Sha02}, we need the following lemma:
\begin{lem}\label{lem4.2}
    There exists an injection $\alpha^\vee: F^\times \longrightarrow Z_G\backslash Z_M$, s.t., $\alpha(\alpha^\vee(t))=t$, for $\forall t\in F^\times$, where $\alpha=\alpha_r$.
\end{lem}

\begin{proof}
    Since $Z_M=\cap_{i\neq r}\ker(\alpha_i)=\{\mathrm{diag}(\underbrace{t,\cdots, t}_{r}, \pm 1,\cdots, \pm 1, \underbrace{t^{-1},\cdots, t^{-1}}_{r})\} $, and $Z_G=\cap_{i=1}^n\ker(\alpha_i)=\{\pm 1\}$, the quotient $Z_G\backslash Z_M$ can be identified as $\{\mathrm{diag}(\underbrace{t,\cdots, t}_{r}, 1, \cdots, 1, \underbrace{t^{-1},\cdots, t^{-1}}_{r}):t\in\mathbb{G}_m\} $.
    Define $$\alpha^\vee(t)=\mathrm{diag}(\underbrace{t,\cdots, t}_{r}, 1,\cdots, 1, \underbrace{t^{-1}, \cdots,t^{-1}}_{r}),$$ then clearly we have $\alpha(\alpha^\vee(t))=t$ for $\forall t\in F^\times$.
\end{proof} 

Let $Z_M^0$ denote the image of $\alpha^\vee$. Notice that $Z_M^0$ is the connected component of the center $Z_M$ of $M$. As we will see later in section 5, the test functions in the local coefficient formula are compactly supported modulo $Z_M$, so it is necessary to consider the action of $Z_M^0$ on $U_M\backslash N$. As a result, the local coefficient is an integral of certain partial Bessel function over the space $Z_M^0U_M\backslash N$. The corresponding orbit space representatives and invariant measures can be given by the following proposition: 
\begin{prop}\label{prop4.3} The action of $Z_M^0$ on $R=R_X\times R_Z=U_M\backslash N$ admits a set of orbit space representatives of the form $R':=R_{X'}\times R_{Z'}$ by setting $x'_{i,j}=\frac{x_{i,j}}{x_{r,1}}$, and $z'_{i,j}=\frac{z_{i,j}}{x^2_{r,1}}$ for elements in $R_X$ and $R_Z$ respectively. The corresponding invariant measure is given by $d\mu'=d\mu_{X'}\wedge d\mu_{Z'}$, where
$$d\mu_{X'}=\begin{cases}

\vert {x'}_{r-1, 2}^{r+2m-5}{x'}_{r-2, 3}^{r+2m-8}\cdots {x'}_{1,r}^{2m-2r+1}\vert \prod dx'_{ij}, & r\le m,\\

\vert {x'}_{r-1, 2}^{r+2m-5}{x'}_{r-2, 3}^{r+2m-8}\cdots {x'}_{r-m+1,m}^{r-m+1}{x'}_{r-m, m+1}^{r-m-1} {x'}_{r-m-1,m+2}^{r-m-2}\cdots {x'}_{2,r-1}\vert \prod dx'_{ij}, & m<r<2m,\\

\vert {x'}_{r-1, 2}^{r+2m-5}{x'}_{r-2, 3}^{r+2m-8}\cdots {x'}_{r-m+1,m}^{r-m+1}{x'}_{r-m, m+1}^{r-m-1} {x'}_{r-m-1,m+2}^{r-m-2}\cdots {x'}_{r-2m+1, 2m}^{r-2m}\vert \prod d{x'}_{ij}. & r\ge 2m,\\
\end{cases}$$
 and $$d\mu_{Z'}=\begin{cases}
     \prod_{i,j}dz'_{i,j} & \mathrm{\ if \ } r\leq 2m,\\
     \vert {z'}_{r-2m, r-2m}^{r-2m-1}{z'}_{r-2m-1, r-2m-1}^{r-2m-2}\cdots {z'}_{2,2}\vert\prod dz'_{ij} & \mathrm{\ if \ } r>2m.
 \end{cases}$$   
\end{prop}

\begin{proof}
For $z=\alpha^\vee(t)$, we have $zn(R_X,R_Y)z^{-1}=n(tR_X,t^2R_Y)$. This is equivalent to say that the action of $Z_M^0$ on $U_M\backslash N'\simeq R_X\times R_Z$ is given by
$$R_X\mapsto tR_X, R_Z\mapsto t^2R_Z$$ From this we identify $Z_M^0U_M\backslash N$ with $R_{X'}\times R_{Z'}$ where $X'$ is given by setting $x_{r,1}=1$ in $R_X$, and $Z'$ is of the same form with $Z$. We will construct the invariant measure $d\mu'$ on the orbit space $R':=Z_M^0U_M
\backslash N'\simeq F^\times\backslash (R_X\times R_Z)=R_{X'}\times R_{Z'}$ such that it is compatible with our invariant measure $d\mu=d\mu_X\wedge d\mu_Z$, in the sense that
$$\int_{R'(F)}\int_{F^\times}f(tX', t^2Z')q^{\langle 2\rho, H_M(\alpha^\vee(t))\rangle}\frac{dt}{\vert t\vert} d\mu'=\int_{R(F)}f(X,Z)d\mu_X\wedge d\mu_Z$$ for any integrable function $f$ on $R=R_X\times R_Z$, where $\rho$ is half sum of the positive roots in $N$.
Observe that $dx_{i,j}=\vert t\vert dx'_{i,j}$, $dz_{i,j}=\vert t\vert^2 dz_{i,j}$, and $$2\rho=\sum_{1\le i\le r,r+1\le j\le n}(e_i\pm e_j)+\sum_{1\le i\le j\le r}(e_i+e_j)=(2m+2)\sum_{i=1}^r e_i+\sum_{1\le i<j\le r}(e_i+e_j)=(2m+r+1)\sum_{i}^r e_i$$ Hence $q^{\langle 2\rho, H_M(\alpha^\vee(t))\rangle}=\vert t\vert^{r(2m+r+1)}$. The above observation implies that the measures on $R_{X'}$ and $R_{Z'}$ must be of the form $d\mu_{X'}=\prod \vert x'_{ij}\vert^{k_{ij}}dx'_{ij}$, and $d\mu_{Z'}=\prod \vert z'_{st}\vert^{l_{st}} dz'_{st}$. Therefore 
$$\int_{R'(F)}\int_{F^\times}f(tX', t^2Z')q^{\langle 2\rho, H_M(\alpha^\vee(t))\rangle}\frac{dt}{\vert t\vert} d\mu'=\int_{F^{\times}\times R'(F)}f(tX',t^2Z')\vert t\vert^{r(2m+r+1)}\frac{dt}{\vert t\vert}(d\mu_{X'}\wedge d\mu_{Z'}) $$
$$=\int_{R(F)}f(X,Z)(\vert x_{r,1}\vert^{r(2m+r+1)-1}\prod_{(i,j)\neq (r,1)} \vert x_{r,1}\vert^{-k_{ij}-1} \vert x_{ij}\vert^{k_{ij}}dx_{ij})\wedge(\prod_{(s,t)}\vert x_{r,1}\vert^{-l_{st}-2}\vert z_{st}\vert^{l_{st}}dz_{st})$$ Compare this with the formulas for $R_X$ and $R_Z$ we obtained, if forces that all the powers of $x'_{i,j}$ and $z'_{i,j}$ remain the same as those of $x_{i,j}$ and $z_{i,j}$ respectively, except $x_{r,1}$, i.e. the invariant measure on $Z_M^0U_M\backslash N$ can be given as $d\mu'=d\mu_{X'}\wedge d\mu_{Z'}$, where 
 $$d\mu_{X'}=\begin{cases}

\vert {x'}_{r-1, 2}^{r+2m-5}{x'}_{r-2, 3}^{r+2m-8}\cdots {x'}_{1,r}^{2m-2r+1}\vert \prod dx'_{ij}, & r\le m,\\

\vert {x'}_{r-1, 2}^{r+2m-5}{x'}_{r-2, 3}^{r+2m-8}\cdots {x'}_{r-m+1,m}^{r-m+1}{x'}_{r-m, m+1}^{r-m-1} {x'}_{r-m-1,m+2}^{r-m-2}\cdots {x'}_{2,r-1}\vert \prod dx'_{ij}, & m<r<2m,\\

\vert {x'}_{r-1, 2}^{r+2m-5}{x'}_{r-2, 3}^{r+2m-8}\cdots {x'}_{r-m+1,m}^{r-m+1}{x'}_{r-m, m+1}^{r-m-1} {x'}_{r-m-1,m+2}^{r-m-2}\cdots {x'}_{r-2m+1, 2m}^{r-2m}\vert \prod d{x'}_{ij}. & r\ge 2m,\\
\end{cases}$$
 and $$d\mu_{Z'}=\begin{cases}
     \prod_{i,j}dz'_{i,j} & \mathrm{\ if \ } r\leq 2m,\\
     \vert {z'}_{r-2m, r-2m}^{r-2m-1}{z'}_{r-2m-1, r-2m-1}^{r-2m-2}\cdots {z'}_{2,2}\vert\prod dz'_{ij} & \mathrm{\ if \ } r>2m.
 \end{cases}$$
Moreover, to match the power of $x_{r,1}$ in the expressions of $d\mu_X$, one also needs to verify that the following identities hold in each cases:

If $r\le m$, then 
   $r(r+2m+1)-1-\sum_{k=1}^{r-1}((r+2m-2-3k)+1)-\sum_{k=1}^{r-1} k -2\sum_{k=1}^rk=r+2m-2$ ;

If $m<r<2m$, then
   
  $ r(r+2m+1)-1-\sum_{k=1}^{m-1}((r+2m-2-3k)+1)-\sum_{k=1}^{r-m-1}((r-m-k)+1)-((r-m)^2-(r-m-1))-\frac{(3r-2m-1)(2m-r)}{2}-2\sum_{k=1}^r k=r+2m-2$;

  If $r\ge 2m$, then
   $r(r+2m+1)-1-\sum_{k=1}^{m-1}((r+2m-2-3k)+1)-\sum_{k=1}^m((r-m-k)+1)-m(m-1)-2\sum_{k=1}^{r-2m}((k-1)+2)-2(\frac{r(r+1)}{2}-\frac{(r-2m)(r-2m-1)}{2})=r+2m-2$. 
These identities follow by straightforward calculations.
\end{proof}

\section{Local coefficient as Mellin transforms of partial Bessel functions}
In this section we apply Theorem 6.2 of \cite{Sha02} to represent the local coefficient in our case as the Mellin transform of certain partial Bessel integrals.
Recall that given irreducible $\psi$-generic representations $\sigma$ and $\tau$ of $\mathrm{GL}_n(F)$ and $\mathrm{Sp}_{2m}(F)$ respectively, then $\pi=\sigma\otimes \tau$ is a $\psi$-generic representation of $M(F)\simeq \mathrm{GL}_n(F)\times \mathrm{Sp}_{2m}(F)$. The corresponding local coefficient is a product of two $\gamma$-factors:
$$C_\psi(s,\pi)=\gamma(s,\sigma\times\tau, \psi)\gamma(s, \sigma, \wedge^2,\psi)$$ As the stability of $\gamma(s,\sigma, \wedge^2,\psi)$ is proved in \cite{CST17}, the stability of $\gamma(s,\sigma\times\tau,\psi)$ is equivalent to the stability of the local coefficient $C_\psi(s,\pi)$.  
Let $n\in N(F)$ such that the Bruhat decomposition $\dot{w}_0^{-1}n=mn'\bar{n}$ holds. Denote by $$U_{M,n}:=\{u\in U_M: unu^{-1}=n\}, \mathrm{\ and \ } 
U_{M,m}':=\{u\in U_{M}: mum^{-1}\in U_M\ \ \& \ \ \psi(mum^{-1})=\psi(u)\}.$$ By \cite{Sun07}, except for a set of measure zero on $N(F)$, we have 
$U_{M,n}=U_{M,m}'$. Together with Lemma \ref{lem4.2}, this implies that the assumptions for Theorem 6.2 in \cite{Sha02} are satisfied. Bessel functions are defined on Bruhat cells. Suppose in the decomposition $\dot{w}_0^{-1}n=mn'\bar{n}$, $m$ lies in some Bruhat cell $C_M(w)=B_M w B_M$, where $B_M=B\cap M$. We define the Bessel function associated to $\pi$ and $w$ as $$j_{\pi,w}(m)=\int_{U_M} W(mu)\psi^{-1}(u)du=\int_{U_M}W(mu^{-1})\psi(u)du$$ where $W\in W(\pi,\psi)$ is a Whittaker model defined by a Whittaker functional $\lambda$ on $\pi$ as $W(m)=\lambda(\pi(m)v)$ for certain $v\in\pi $, normalized so that $W(e)=1$, where $e$ is the identity element of $M(F)$. Then It is immediate that $$j_{\pi,w}(u_1m u_2)=\psi(u_1u_2)j_{\pi,w}(m)$$ for any $u_1,u_2\in U_M(F)$. It follows by the Bruhat decomposition that the Bessel functions are essentially functions on $wT$, where $T$ is the maximal (split) torus of $G$. We refer the readers to \cite{CpS05} for partial Bessel funcitons for quasi-split groups.

However, the local coefficient is represented by an integral of certain partial Bessel function, not the Bessel function itself. To define that, we need an exhaustive sequence of open compact subgroups $\overline{N}_{0,\kappa}\subset \overline{N}(F)$ such that $\alpha^\vee(t)\overline{N}_{0,\kappa}\alpha^\vee(t)^{-1}$ depends only on $\vert t\vert$. For a matrix $X$ of size $l\times k$, let
$$\varphi_\kappa(X)=\begin{cases}
    1 & \vert X_{i,j}\vert\leq q^{((k-i)+(l-j)+1)\kappa}, \\
    0 & \mathrm{\ otherwise \ }.\\
\end{cases}$$
Set $$\overline{N}_{0,\kappa}=\{\bar{n}(X,Y): \varphi_\kappa(\varpi^{-(d+g)}X)\cdot\varphi_\kappa(\varpi^{-2(d+g)}Y)=1\}$$ where $d=\mathrm{Cond}(\psi)$ and $g=\mathrm{Cond}(\omega_\pi^{-1}(w_0(\omega_\pi)))$. 
 As we have seen before that in our case $\alpha^\vee(t)\bar{n}(X,Y)\alpha^\vee(t)^{-1}=\bar{n}(tX, t^2Y)$, it follows that our definition of $\overline{N}_{0,\kappa}$ makes $\alpha^\vee(t)\overline{N}_{0,\kappa}\alpha^\vee(t)^{-1}$ depends only on $\vert t\vert$. Denote by $\varphi_{\overline{N}_{0,\kappa}}$ the characteristic function of $\overline{N}_{0,\kappa}$.
Then the partial Bessel function is a function on $M(F)\times Z_M^0(F)$ defined as
$$j_{\overline{N}_{0,\kappa}, \pi, w}(m,z):=\int_{U_{M,n}\backslash U_M}W(mu)\varphi_{\overline{N}_{0,\kappa}}(zu^{-1}\overline{n}uz^{-1})\psi^{-1}(u)du.$$ Given that $\dot{w}_0^{-1}n=mn'\bar{n}$, let  
$$j_{\overline{N}_{0,\kappa}\pi,w}(n):=j_{\overline{N}_{0,\kappa}, \pi, w}(m,\alpha^{\vee}(\varpi^{d+g}u_{\alpha_r}(\dot{w}_0\bar{n}\dot{w}_0^{-1}))).$$ Suppose $\sigma$ and $\tau$ are irreducible $\psi$-generic representations of $\mathrm{GL}_r(F)$ and $\mathrm{Sp}_{2m}(F)$ respectively,  and $L=\prod_{k=1}^d\mathrm{GL}_{n_k}\times \prod_{j=1}^l\mathrm{GL}_{t_j}\times \mathrm{Sp}_{2m'}$ with $\sum_{k=1}^d n_k=r$ and $\sum_{j=1}^l n_j+m'=m$. Then we can find $\psi$-generic supercuspidal representations $\sigma'_k,(1\le k\le d)$, $\sigma''_{j},(1\le j\le l)$, and $\tau'$ of $\mathrm{GL}_{n_k}(F)$, $\mathrm{GL}_{t_j}(F)$, and $\mathrm{Sp}_{2m'}(F)$ respectively, so that $\pi=\sigma\boxtimes \tau\hookrightarrow \mathrm{Ind}_{L(F)N(F)}^{M(F)}\sigma'\otimes \textbf{1}_{N(F)}$, where
$\sigma'=\otimes_{k=1}^d\sigma'_{k}\otimes\otimes_{j=1}^l \sigma''_j\otimes\tau'$. By the multiplicativity of Langlands-Shahidi local $\gamma$-factors in our cases, see Theorem 3.1 of \cite{sha90} or Theorem 8.3.2 of \cite{sha10}, we obtain that
$$\gamma(s,\sigma\times \tau,\psi)=\prod_{k=1}^d\prod_{j=1}^l\gamma(s,\sigma'_k\times\sigma''_j,\psi)\gamma(s,\sigma_k\times\tilde{\sigma}''_j,\psi)\prod_{k=1}^d\gamma(s, \sigma'_k\times \tau',\psi)$$ where $\tilde{\sigma}''_j$ is the contragredient of $\sigma''_j$. Consequently, it suffices to prove stability for $\psi$-generic supercuspidal representations. Suppose $\sigma$ and $\tau$ are $\psi$-generic supercuspidal representations of $\mathrm{GL}_n(F)$ and $\mathrm{Sp}_{2m}(F)$ respectively. Set $\pi=\sigma\boxtimes \tau$, then its central character $\omega_\pi=\omega_\sigma\boxtimes \omega_\tau$. In this case, one can find a function $f\in C^\infty_c(M(F);\omega_\pi)$, the space of smooth functions of compact support modulus the center on $M(F)$, such that $f(zm)=\omega_\pi(m)f(z)$ for all $z\in Z_M(F)$, and $m\in M(F)$. Then $W(m)=W_f(m):=\int_{U_M(F)}f(xm)\psi^{-1}(x)dx$ defines a non-zero Whittaker model. We normalize it so that $W_f(1)=1$. For simplicity, we omit $w$ and denote
$j_{\pi,\kappa, f}(n):=j_{\overline{N}_{0,\kappa},\pi, w}(n)$ where we replace $W$ by $W_f$. Let $\alpha=\alpha_r$ and $\tilde{\alpha}=\langle \rho, \alpha\rangle ^{-1}\rho$, where $\rho$ is half sum of the positive roots in $N$. Then by Theorem 6.2 of \cite{Sha02}, we have
\begin{prop}
     Let $\sigma$ and $\tau$ be $\psi$-generic supercuspidal representations of $\mathrm{GL}_r(F)$ and $\mathrm{Sp}_{2m}(F)$ respectively, set $\pi=\sigma\boxtimes \tau$, then
$$C_\psi(s,\pi)^{-1}=\gamma(2\langle\tilde{\alpha},\alpha^\vee\rangle s, \omega_\pi(\dot{w}_0\omega_\pi^{-1}), \psi^{-1})\int_{Z_M^0(F)U_M(F)\backslash N(F)}j_{\pi,\kappa, f}(\dot{n})\omega_{\pi_s}^{-1}(\dot{w}_0\omega_{\pi_s})(x_\alpha)q_F^{\langle s\tilde{\alpha}+\rho, H_M(m)\rangle}d\dot{n}$$ for sufficiently large $\kappa$,
where $m$ is the image of $\dot{n}\mapsto m$ with $\dot{n}$ a representative of the orbit of $n \in Z_M^0(F)U_M(F)\backslash N(F)$ in $N(F)$ via $\dot{w}_0^{-1}\dot{n}=mn'\bar{n}$, which holds off a subset of measure zero on $N(F)$. Here $d\dot{n}$ is the invariant measure on the orbit space $Z_M^0(F)U_M(F)\backslash N(F)$, $\pi_{s}=\pi\otimes q^{\langle s\tilde{\alpha}, H_M(\cdot)\rangle}$,  $x_\alpha=u_{\alpha_r}(\dot{w}_0\bar{n}\dot{w}_0^{-1})$, and $\gamma(2\langle\tilde{\alpha},\alpha^\vee\rangle s, \omega_\pi(\dot{w}_0\omega_\pi^{-1}), \psi^{-1})$ is the Abelian $\gamma$-factor depending only on $\omega_\pi$.
\end{prop}

Let us simplify this formula in our cases. First, recall that we have 
$2\rho=(r+2m+1)\sum_{i=1}^re_i$,
$\alpha=\alpha_r=e_r-e_{r+1}$. So
$$\langle\rho,\alpha\rangle=\frac{2(\rho,\alpha)}{(\alpha,\alpha)}=\frac{((r+2m+1)\sum_{i=1}^r e_i, e_r-e_{r+1})}{(e_r-e_{r+1}, e_r-e_{r+1})}=\frac{r+2m+1}{2}$$ and therefore $\tilde{\alpha}=\langle \rho, \alpha\rangle^{-1}\rho=\sum_{i=1}^re_i$. It follows that
$\langle\tilde{\alpha},\alpha^\vee\rangle=\langle \sum_{i=1}^re_i, \sum_{i=1}^re_i^*\rangle=r$, and if we write $m=\mathrm{diag}\{m_1, m_2, \theta_r(m_1)\}$, then $$q^{\langle s\tilde{\alpha}, H_M(m)\rangle}=\vert \det(m_1)\vert^{s}, q^{\langle s\tilde{\alpha}+\rho, H_M(m)\rangle}=\vert \det(m_1)\vert^{s+\frac{r+2m+1}{2}}.$$ 

Moreover,
$$\omega_\pi(\dot{w}_0\omega_\pi^{-1})(\alpha^\vee(t))=\omega_\pi(\alpha^\vee(t))\omega_\pi^{-1}(\dot{w}_0^{-1}\alpha^\vee(t)\dot{w}_0)$$ 
$$=\omega_\pi(\mathrm{diag}(\underbrace{t,\cdots,t}_r,1,\cdots,1,\underbrace{t^{-1},\cdots,t^{-1}}_r))\omega_\pi^{-1}(\mathrm{diag}(\underbrace{t^{-1}, \cdots,t^{-1}}_r, 1,\cdots, 1, \underbrace{t,\cdots, t}_r))=\omega_\pi^2(\alpha^\vee(t)).$$ Hence $\omega_\pi(\dot{w}_0\omega_\pi^{-1})\circ\alpha^\vee=\omega_\pi^2\circ\alpha^\vee$. Similarly
$\omega_{\pi_s}^{-1}(\dot{w}_0\omega_{\pi_s})(\alpha^\vee(t))=\omega^{-2}_{\pi_s}(\alpha^\vee(t))=\omega_{\pi}^{-2}(\alpha^\vee(t))\vert t\vert^{-rs}$. From our previous calculation of the bruhat decomposition $\dot{w}_0^{-1}n=mn'\bar{n}$ with $n=n(X,Y)$, we see that $$\dot{w}_0\bar{n}\dot{w}_0^{-1}=\begin{bmatrix}
    & & I_r \\
    & I_{2m} & \\
    (-1)^rI_r & & \\
\end{bmatrix}\begin{bmatrix}
    I_r & & \\
    J_{2m}' {^tX}{^tY^{-1}}J_r & I_{2m} & \\
    Y^{-1} & Y^{-1}X & I_r \\ 
\end{bmatrix}\begin{bmatrix}
     & & (-1)^r I_r\\
     & I_{2m} & \\
     I_r & & \\
\end{bmatrix}$$$$=\begin{bmatrix}
    I_r & Y^{-1}X & Y^{-1}\\
    & I_{2m} & (-1)^rJ_{2m}'{^tX}{^tYY^{-1}}J_r\\
    & & I_r\\
\end{bmatrix}=n(Y^{-1}X, Y^{-1}).$$ Hence $u_{\alpha_r}(\dot{w}_0\bar{n}\dot{w}_0^{-1})$ is the lower-left entry of $Y^{-1}X$. When we take our orbit space representative $R=R_X\times R_Y\simeq R_X\times R_Z\simeq U_M\backslash N$, and $n=n(X, Y)$ with $(X,Y)\in R_X\times R_Y$, this is just $\frac{y^*_{rr}x_{r,1}}{\det Y}$, where $y^*_{rr}$ is the $(r,r)$-th entry of the adjoint matrix of $Y=ZJ_r+\frac{X\theta_{r,m}(X)}{2}$, which is a polynomial function of $(X, Z)$, we denote it as $P(X,Z)$. When passing to the orbit of $Z_M^0(F)$, we obtain that $u_{\alpha_r}(\dot{w}_0\bar{n}(X',Y')\dot{w}_0^{-1})=\frac{P(X', Z')}{\det Y'}$. 

From now on, we will also use $n(X, Z)$(resp. $m(X,Z)$) to denote $n(X,Y)$(resp. $m(X,Y)$) if we emphasis that it is parameterized by $(X,Z)$ instead of $(X,Y)$. With the discussion of the orbit space structure of $Z_M^0U_M\backslash N$ and its invariant measure in Proposition \ref{prop4.3}, we obtain the following result:

\begin{prop}\label{prop5.2}
     Let $\sigma$ and $\tau$ be $\psi$-generic supercuspidal representations of $\mathrm{GL}_r(F)$ and $\mathrm{Sp}_{2m}(F)$ respectively, set $\pi=\sigma\boxtimes \tau$, then
$$C_\psi(s,\pi)^{-1}=\gamma(2rs, \omega_\pi^2, \psi^{-1})\int_{R_{X'}\times R_{Z'}}j_{\pi,\kappa, f}(n(X',Z'))\omega_{\pi}^{-2}(\frac{P(X',Z')}{\det(Z'J_r+\frac{X' \theta_{r,m}(X')}{2})})$$$$\cdot \vert \frac{P(X',Z')}{\det(Z'J_r+\frac{X'\theta_{r,m}(X')}{2})}\vert^{-rs} \vert \det(m_1(X',Z'))\vert^{s+\frac{r+2m+1}{2}} d\mu_{X'}\wedge d\mu_{Z'}$$ for sufficiently large $\kappa$,
where $m(X',Z')$ is the image of $n(X',Z')\mapsto m(X', Z')$ via the Bruhat decomposition $\dot{w}_0^{-1}\dot{n}=mn'\bar{n}$,  which holds off a subset of measure zero on $N(F)$, $m=\mathrm{diag}\{m_1, m_2,\theta_r(m_1)\}$ with $m_1\in\mathrm{GL}_r(F)$, $m_2\in\mathrm{Sp}_{2m}(F)$, and $\gamma(2r s, \omega_\pi^2, \psi^{-1})$ is the Abelian $\gamma$-factor depending only on $\omega_\pi$.
\end{prop}
To study the stability of local coefficient, we need to consider $C_\psi(s,\pi\otimes \chi)$ with $\chi$ a highly ramified character of $F^\times$, regarded as a character of $M(F)\simeq \mathrm{GL}_r(F)\times \mathrm{Sp}_{2m}(F)$ via $\chi(m_1, m_2):=\chi(\det(m_1))$ in our case. Therefore it is necessary to choose the open compact subgroups $\{\overline{N}_{0,\kappa}\}_{\kappa}$ of $\overline{N}(F)$ to be independent of $\chi$. As in the proof of Theorem 6.2 of \cite{Sha02}, given an irreducible $\psi$-generic representation $\pi$ of $M(F)$ with ramified central character $\omega_\pi$, we choose a section $h\in I(s, \pi)=\mathrm{Ind}_{P(F)}^{G(F)}(\pi\otimes q^{\langle s\tilde{\alpha}+\rho, H_M(\cdot)\rangle})\otimes \textbf{1}_{N(F)}$, such that $h$ is compactly supported modulus $P(F)$, and use it to obtain the integral representation of the local coefficient formula as above. So we choose a sufficiently large open compact subgroup $\overline{N}_0$ of $\overline{N}(F)$ such that $\mathrm{Supp}(h)\subset P(F)\overline{N}_0$. In our situation, we fix a character $\chi_0$ of $F^\times$ such that $\omega_{\pi\otimes\chi_0}=\omega_\pi\chi_0^r$ is ramified, and choose a $\kappa_0$ and $h_0\in I(s,\pi)$ such that $\mathrm{Supp}(h_0)\subset P(F)\overline{N}_{0,\kappa_0}$. Since $\overline{N}_{0,\kappa}\subset\overline{N}_{0,\kappa_0}$ for all $\kappa\ge\kappa_0$, $\overline{N}_{0,\kappa}$  is independent of $\chi_0$ if $\kappa\ge\kappa_0$. Suppose $\chi$ is any character of $F^\times$ such that $\omega_{\pi\otimes \chi}$ is ramified, choose $h_\chi\in I(s,\pi\otimes \chi)$ such that $\mathrm{Supp}(h_\chi)\subset P(F)\overline{N}_{0,\chi}$ for some open compact subgroup $\overline{N}_{0,\chi}$ of $\overline{N}(F)$. If $\overline{N}_{0,\chi}\subset \overline{N}_{0,\kappa_0}$, we are fine as we just discussed. If not, we replace $h_\chi$ by a right shift $R(\alpha^\vee(t))h_\chi: g\mapsto h_\chi(g\alpha^\vee(t))$, then $R(\alpha^\vee(t)h_\chi)$ is supported in $\alpha^\vee(t)\overline{N}_{0,\chi}\alpha^\vee(t)^{-1}$ mod $P(F)$, as $\alpha^\vee(t)\in M(F)$. Since $\alpha^\vee(t)\bar{n}(X,Y)\alpha^\vee(t)^{-1}=\bar{n}(tX, t^2Y)$, we choose $\vert t\vert$ to be small enough so that $\alpha^\vee(t)\overline{N}_{0,\chi}\alpha^\vee(t)^{-1}\subset \overline{N}_{0,\kappa_0}$, then we are done. So we obtain the following stronger version of the above proposition:

\begin{prop}\label{prop5.3}
     Let $\sigma$ and $\tau$ be $\psi$-generic supercuspidal representations of $\mathrm{GL}_r(F)$ and $\mathrm{Sp}_{2m}(F)$ respectively, such that $\omega_\pi$ is ramified where $\pi=\sigma\boxtimes \tau$, then there exists a $\kappa_0$, such that for all $\kappa\ge\kappa_0$ and all characters $\chi$ of $F^\times$ such that $\omega_{\pi\otimes\chi}=\omega_\pi\chi^r$ is ramified, we have
$$C_\psi(s,\pi\otimes \chi)^{-1}=\gamma(2rs, \omega_\pi^2\chi^{2r}, \psi^{-1})\int_{R_{X'}\times R_{Z'}}j_{\pi\otimes\chi,\kappa, f}(n(X',Z'))(\omega_{\pi}^{-2}\chi^{-2r})(\frac{P(X',Z')}{\det(Z'J_r+\frac{X' \theta_{r,m}(X')}{2})})$$$$\cdot \vert \frac{P(X',Z')}{\det(Z'J_r+\frac{X'\theta_{r,m}(X')}{2})}\vert^{-rs} \vert \det(m_1(X',Z'))\vert^{s+\frac{r+2m+1}{2}} d\mu_{X'}\wedge d\mu_{Z'}$$
where $m(X',Z')$ is the image of $n(X',Z')\mapsto m(X', Z')$ via the Bruhat decomposition $\dot{w}_0^{-1}\dot{n}=mn'\bar{n}$,  which holds off a subset of measure zero on $N(F)$, $m=\mathrm{diag}\{m_1, m_2,\theta_r(m_1)\}$ with $m_1\in\mathrm{GL}_r(F)$, $m_2\in\mathrm{Sp}_{2m}(F)$, and $\gamma(2r s, \omega_\pi^2\chi^{2r}, \psi^{-1})$ is the Abelian $\gamma$-factor depending only on $\omega_\pi$ and $\chi$.
\end{prop}
\section{Asymptotics of partial Bessel integrals}
So far we obtained our local coefficient as the Mellin transform of certain partial Bessel functions. The partial Bessel functions appeared in our discussions can be reformulated using partial Bessel integrals, which admit nice asymptotic expansion formulas, as well as some uniform smooth properties. Theses properties are crucial for the proof of stability using our method. In this section, we discuss some general properties of partial Bessel integrals and obtain the asymptotic formula and uniform smoothness in our case. In particular, different from the known results using this method, we observe that the orbit space $Z_M^0U_M\backslash N$ is no longer isomorphic to a torus, and we will separate its 'toric'-part out, which plays the same role as the torus over which the integral representing the local coefficient is taken in the known cases. This is a new phenomenon and we believe it can be generalized in our future work. 
\subsection{Some general properties of partial Bessel integrals}
We begin by introducing partial Bessel integrals and some of its important properties. The structure and results in this subsection is the same as section 6 of \cite{She23} but are under a more general setting, so we will reprove some of the important results. Let $M$ be a connected reductive group defined over a local field $F$. Fix a Borel subgroup $B_M=AU_M$, where $A$ is the maximal (split) torus and $U_M$ is the unipotent radical of $B_M$. Suppose $\pi$ is a $\psi$-generic supercuspidal representation of $M(F)$ with central character $\omega_\pi$. Take a matrix coefficient $f$ of $\pi$, then $f\in C^\infty_c(M(F);\omega_\pi)$, the space of smooth functions of compact support on $M(F)$ modulo the center, such that $f(zm)=\omega_\pi(z)f(m)$ for all $z\in Z_M(F)$, $m\in M(F)$. Then the integral $W_f(m)=\int_{U_M(F)}f(xm)\psi^{-1}(x)dx$ converges since $Z_M(F)U_M(F)m$ is closed in $M(F)$, and hence defines a non-zero Whittaker model attached to $\pi$. We normalize it so that $W_f(e)=1$, where $e$ is the identity element of $M(F)$. 

From now on, we will not distinguish algebraic group and the group of its $F$-points unless it is necessary in order to simplify the notations, we hope this will not cause any troubles for the readers. Given an $F$-involution $\Theta_M: M\rightarrow M$, we define the partial integral as
$$B^M_\varphi(m,f):=\int_{U_{M,m}^{\Theta_M}\backslash U_M}W_f(mu)\varphi(\Theta_M(u^{-1})m'u)\psi^{-1}(xu)du=\int_{U_{M,m}^{\Theta_M}\backslash U_M}\int_{U_M}f(xmu)\varphi(\Theta_M(u^{-1})m'u)\psi^{-1}(xu)dxdu$$
where $U_{M,m}^{\Theta_M}=\{m\in U_M: \Theta_M(u^{-1})mu=m\}$ is the twisted centralizer of $m$ in $U_M$, and $\varphi$ is certain characteristic function of a subset of $M(F)$. And $m'$ is obtained from $m$ by stripping off the center, i.e. $m=zm'$ with $z\in Z_M$. Similarly one can define partial Bessel integrals on any Levi subgroup $L$ of $M$ by an involution $\Theta_L$ of $L$.

Let $\Delta_M$ denote the set of simple roots in $M$, $W(M,A)$ the Weyl group, and $w_M\in W(M,A)$ the long Weyl group element. The following objects needed for our study of the asymptotic expansion of partial Bessel integrals.

\begin{itemize}
    \item {\textbf{\textrm{B(M)}.}} Following \cite{CST17}, the subset of Weyl group elements that supports Bessel functions is given by $B(M)=\{w\in W(M, A): \alpha\in \Delta_M \ \ s.t. \ \ w\alpha>0 \Rightarrow w\alpha\in \Delta_M \}$, there is a bijection $$B(M)\leftrightarrow \{L: \textrm{ Levi\ of\ standard\ parabolic\ subgroups \ of\  } M\}$$ by
    $w\mapsto L=Z_M(\cap_{\alpha\in\theta^+_{M,w}} \ker\alpha)$, where  $\theta_{M,w}^+=\{\alpha\in \Delta_M: w\alpha>0\}$, and conversely
$L \mapsto w=w_M w_L^{-1}$.

\item
$\textbf{U}_{M,w}^+, \textbf{U}_{M,w}^-.$ For each $w\in W(M,A)$, define
$U_{M,w}^+=\{u\in U_M: wuw^{-1}\in U_M\}$
and $U_{M,w}^-=\{u\in U_M:wuw^{-1}\in U_M^-\}.$ Then $U_M=U^+_{M,w}U^-_{M,w}$ and $U_{M,w}^+$ normalizes $U_{M,w}^-$.
In particular, if $w\in B(M)$ with $\dot{w}=\dot{w}_M\dot{w}_L^{-1}$, then $U_{M,w}^+=U_L:=U_M\cap L, U_{M,w}^-=N_L$, the unipotent radical of the standard parabolic subgroup of $M$ with Levi component $L$. In the extreme cases, if $w=w_L$, then
$U_{M,w_L}^+=N_L, U_{M,w_L}^-=U_L$. If $w=w_M$, then
$U_{M,w_M}^+=\{e\}, U_{M,w_M}^-=U_M$.

\item\textbf{Bessel distance} For $w,w'\in B(M)$ with $w>w'$ define
$$d_B(w,w')=max\{m: \exists w_i\in B(M)\ \ s.t \ \ w=w_m>w_{m-1}>\cdots>w_0=w' \}.$$

\item \textbf{Bruhat order} For $w,w'\in W(M,A)$, $w\leq w'\Longleftrightarrow C(w)\subset \overline{C(w')}.$

\item\textbf{The relevant torus $\textbf{A}_w$}.
For $w\in B(M)$, define $A_w=\{a\in A: a \in \cap_{\alpha\in \theta_{M,w}^+}\ker \alpha\}^\circ\subset A$. Note that it is also the connected center of $L_w=Z_M(\cap_{\alpha\in \theta_{M,w}^+}\ker \alpha)$.

\item \textbf{The relevant Bruhat cell} $\textbf{C}_{r}(\dot{w})$. We call $C_r(\dot{w})=U_M\dot{w}A_wU_{M,w}^-$ the relevant part of the Bruhat cell $C(w)$. Note that $C_r(\dot{w})$ depends on the choice of the representative $\dot{w}$ of $w$. We choose the representative $\dot{w}$ of $w$ so that it is compatible with $\psi$ in the sense that $\psi(\dot{w}u\dot{w}^{-1})=\psi(u)$ for all $u\in U_{M,w}^+$, as in Proposition 5.1 of \cite{CST17}.
\item \textbf{Transverse tori} For $w,w'\in B(M)$, let $L=L_w$ and $L'=L_{w'}$ be their associated Levi subgroups respectively. Suppose $w'\leq w$. Then $L\subset L'$ and $A_{w'}\supset A_w$. Let $A^{w'}_w=A_w\cap L'_{der}=Z_L\cap L'_{der}$. In particular  $A^{w'}_w\cap A_{w'}=A^{w'}_{w'}=Z_{L'}\cap L'_{der}$ is finite and the subgroup $A^{w'}_w A_{w'}\subset A_w$ is open and of finite index. So this decomposition is essentially a "transfer principal" for relevant tori.
\end{itemize}
Let $\{x_\alpha\}_{\alpha\in \Delta_M}$ be a set of 1-parameter subgroups, i.e. $x_\alpha:\mathbb{G}_a\simeq U_\alpha$, Then $(B_M, A, \{x_{\alpha}\}_{\alpha\in \Delta_M})$ defines an $F$-splitting of $M$. We make the following assumptions from now on:

1. $\Theta_M$ preserves the $F$-splitting and is compatible with Levi subgroups and $\psi$. To be precise, if $L\subset M$ is a Levi subgroup, we require that $\tau^M_L\circ \Theta_M\vert_L=\Theta_L$, where $\tau^M_L=\mathrm{Int}((\dot{w}^{M}_L)^{-1})$. And $\Theta_L$ also preserves the $F$-splitting $\{B_L:=B_M\cap L, A, \{x_\alpha\}_{\alpha\in \Delta_L}\}$. Here $\mathrm{Int}(w)(m)=wmw^{-1}$, and we denote $\dot{w}^M_L:=\dot{w}_M\dot{w}_L^{-1}$. Consequently $\Theta_M(U_M)= U_M$, $\Theta_L(U_L)=U_L$. The representatives $\dot{w}^M_L$ are chosen to be compatible with the generic character $\psi$ as in Proposition 5.1 of \cite{CST17}, i.e., $\psi(\tau^M_L(u))=\psi(\dot{w}^M_L u(\dot{w}^M_L)^{-1})=\psi(u)$ for all $u\in U_L(F)$. By $\Theta_M$ is compatible with $\psi$, we mean that $\psi(\Theta_M(u))=\psi(u)$ for all $u\in U_M(F)$. As a result, if $u\in U_L(F)$, then $\psi(\Theta_L(u))=\psi(\tau^M_L\circ \Theta_M(u))=\psi(\Theta_M(u))=\psi(u)$.

2. $\varphi$ is the characteristic function of some subset of $M(F)$ so that it is invariant under the $\Theta_M$-twisted conjugate action by some open compact subgroup $U_0$ of $U_M(F)$.i.e. $\varphi(\Theta_M(u^{-1})mu)=\varphi(m)$ if $u\in U_0$.

These assumptions guarantees some nice behavior of partial Bessel functions, and are satisfied in all the known cases. We will see later that they also hold in our cases.

We would like to first understand the structure of twisted centralizer in the definition of partial Bessel integrals, in order to get a descent formula for partial Bessel integrals to those defined over Levi subgropus of $M$. This serves as an important step to prove the general asymptotic expansion formula for partial Bessel integrals. 
\begin{lem}\label{lem6.1}
    For $m=u_1\dot{w}au_2\in C(w^M_L)$, then $U^{\Theta_M}_{M,m}\subset u_2^{-1}U_Lu_2=u_2^{-1}U_{M,w^M_L}^+u_2$.
\end{lem}

\begin{proof}
We have $$u\in U^{\Theta_M}_{M,m}\Leftrightarrow \Theta_M(u^{-1})u_1w_Mw_L^{-1}au_2u=u_1w_Mw_L^{-1}au_2\Leftrightarrow w_M^{-1}u_1^{-1}\Theta_M(u^{-1})u_1w_M=w_L^{-1}au_2u^{-1}u_2^{-1}a^{-1}w_L.$$
    Since $\Theta_M(U_M)=U_M$, the left hand side of the last equality lies in $U_M^{-}$. Hence
    $au_2u^{-1}u_2^{-1}a^{-1}\in U_{w_L}^-=U_L$, therefore $u_2u^{-1}u_2^{-1}\in a^{-1}U_L a=U_L$, so $u\in u_2^{-1}U_Lu_2$.  
\end{proof}

\begin{lem}\label{lem6.2}
    Let $H\subset L\subset M$ be Levi subgroups of $M$, then for $m=u_1\dot{w}^{M}_H au_2\in C_r(\dot{w}^M_H)$, we have
    $$U^{\Theta_{M}}_{M,m}=(\tilde{u}_1^-)^{-1}U^{\Theta_L}_{L,l}\tilde{u}_1^-\cap (u_2^-)^{-1} U^{\Theta_L}_{L,l}u_2^-$$
    where $u_1=u_1^-u_1^+$ with $u_1^-\in U_{M,w^{-1}}^-$ and $u_1^+\in U_{M,w^{-1}}^+$, $u_2=u_2^+u_2^-$ with $u_2^+\in U_{M,w}^+$ and $u_2^-\in U_{M,w}^{-}$, $\tilde{u}_1^-=\tau^M_L\circ \Theta_M\circ \tau^M_L((u_1^-)^{-1})$, $l=w^{-1}u_1^+w w^L_Ha u_2^+\in L$, where $w=w^{M}_L\in B(M)$.
\end{lem}
\begin{proof} By Lemma 6.1 of \cite{She23}, $m=u_1\dot{w}^{M}_Ha u_2\in C_r(\dot{w}^M_H)\subset \Omega_{w}=U_{M,w^{-1}}^-\times wL\times U_{M, w}^-$, where $\Omega_w:=\bigsqcup_{w\leq w'} C(w')$, and the decomposition is unique. Decompose $u_1=u_1^-u_1^+$, $u_2=u_2^+u_2^-$ as stated in the Lemma. Then we can write $m=u_1^- \dot{w}(\dot{w}^{-1}u_1^+\dot{w}\dot{w}^L_H a u_2^+)u_2^-$ and it is easy to see that $l=\dot{w}^{-1}u_1^+\dot{w}\dot{w}^L_H a u_2^+\in L$. Then
$$u\in U^{\Theta_M}_{M,m}\Leftrightarrow \Theta_M(u^{-1})mu=m\Leftrightarrow (\Theta_M(u^{-1})u_1^-\dot{w}\Theta_L(u^+)\dot{w}^{-1})\dot{w} \Theta_L((u^+)^{-1})lu^+((u^+)^{-1}u_2^-u^+u^-)=u_1^-\dot{w}lu_2^-$$
where we decompose $u=u^+u^-$ with $u^+\in U_{M,w}^+=U_L$ and $u^-\in U_{M,w}^-=N_L$. We claim that $$\Theta_M(u^{-1})u_1^-\dot{w}\Theta_L(u^+)\dot{w}^{-1}\in U_{L,w^{-1}}^-.$$ To see this, note that it is equivalent to  $\dot{w}^{-1}\Theta_M(u^{-1})u_1^-\dot{w}\Theta_L(u^+)\in U_{M}^-$. But $\dot{w}^{-1}\Theta_M(u^{-1})u_1^-\dot{w}\Theta_L(u^+)=\dot{w}^{-1}\Theta_M((u^-)^{-1})\dot{w} (\dot{w}^{-1}\Theta_M((u^+)^{-1})\dot{w})(\dot{w}^{-1}u_1^-\dot{w})\Theta_L(u^+)$. Since $(u^-)^{-1}\in N_L$, and by our assumption,
$\tau^M_L\circ \Theta_M$ preserves the $F$-splitting when restricted to $L$, so $\Theta_M(N_L)=N_L$. In addition, since $\dot{w}=\dot{w}^M_L$, we have $\dot{w}^{-1}N_L{\dot{w}}\cap U_L=\{1\}$. As a result, $\dot{w}^{-1}\Theta_M((u^-)^{-1})\dot{w}\in U_M^-$. This also implies that $\tau^M_L\circ \Theta_M\vert_{N_L}\subset U_M^-$. Moreover, $\dot{w}^{-1}\Theta_M((u^+)^{-1})\dot{w}=\tau^M_L\circ\Theta_L((u^+)^{-1})=\Theta_L((u^+)^{-1})\in U_L$ by our assumption on $\Theta_L$. We also have $\dot{w}^{-1}U_{M,w^{-1}}^-\dot{w}=(U_{M,w}^-)^-=N_L^-$, so $\dot{w}^{-1}u_1^-\dot{w}\in N_L^-$. Since $U_L$ normalizes $N_L^-$, $(\dot{w}^{-1}\Theta_M((u^+)^{-1})\dot{w})(\dot{w}^{-1}u_1^-\dot{w})\Theta_L(u^+)\in N_L^-\subset U_M^-$. This proves the claim. It follows from the claim and the uniqueness of the decomposition 
$\Omega_{w}=U_{L,w^{-1}}^-\times \dot{w}L\times U_{L, w}^-$ that
$$\textrm{(1).}\  \Theta_M(u^{-1})u_1^-\dot{w}\Theta_L(u^+)\dot{w}^{-1}=u_1^-,\ \textrm{(2).}\  \Theta_L((u^+)^{-1})lu^+=l,\ \textrm{(3).}\ (u^+)^{-1}u_2^-u^+u^-=u_2^-.$$
Note that (1)$\Leftrightarrow \dot{w}^{-1}(u_1^-)\Theta_M(u^{-1})u_1^-\dot{w}=\Theta_L(u^+)^{-1}\Leftrightarrow \tau^M_L(u_1^-)\tau^M_L\circ\Theta_M(u^{-1})\tau^M_L(u_1^-)=\Theta_L(u^+)^{-1}$. Apply $\tau^M_L\circ \Theta_M$ on both sides, and use $\tau^M_L\circ \Theta_M\vert_L=\Theta_L$, then take inverse on both sides, we obtain that $\tilde{u}_1^- u(\tilde{u}_1^-)^{-1}=u^+$, where $\tilde{u}_1^-=\tau^M_L\circ \Theta_M\circ \tau^M_L((u_1^-)^{-1})$. (2) is saying that $u^+\in U^{\Theta_L}_{L,l}$. Therefore (1) and (2) imply that $U^{\Theta_M}_{M,m}\subset (\tilde{u}_1^-)^{-1}U^{\Theta_L}_{L,l}\tilde{u}_1^-$. On the other hand, (3)$\Leftrightarrow u_2^-u(u_2^-)^{-1}=u^+$, then (2) and (3) imply that $U^{\Theta_M}_{M,m}\subset (u_2^-)^{-1}U^{\Theta_L}_{L,l}u_2^-$. Hence $$U^{\Theta_{M}}_{M,m}\subset(\tilde{u}_1^-)^{-1}U^{\Theta_L}_{L,l}\tilde{u}_1^-\cap (u_2^-)^{-1} U^{\Theta_L}_{L,l}u_2^-.$$
Conversely, if $u=(\tilde{u}_1^-)^{-1}u'\tilde{u}_1^-=(u_2^-)^{-1}u'' u_2^- $ with $u',u''\in U^{\Theta_L}_{L,l}$, then
$u=u^+u^-=u' (u'^{-1}(\tilde{u}_1^-)^{-1}u'\tilde{u}_1^-)=u''((u'')^{-1}(u_2^-)^{-1}u'' u_2')$. One can see from the above argument that $\tilde{u}_1^-, u_2^-\in N_L$. Since $U_L\cap N_L=\{1\}$, and $U_L$ normalizes $N_L$, this implies that $u^+=u'=u''$, this is equivalent to (2). Replace $u'$ and $u''$ by $u^+$ we obtain that $\tilde{u}_1^- u(\tilde{u}_1^-)^{-1}=u^+$ and $u_2^-u(u_2^-)^{-1}=u^+$, which are equivalent to (1) and (3) respectively. Reverse the above process we obtain that (1), (2), (3) hold simultaneously, hence the reversed inclusion. 
\end{proof}

By Lemma 5.9 of \cite{CST17}, we have a surjective map $C^{\infty}_c(\Omega_{w};\omega_{\pi})\rightarrowdbl  C^{\infty}_c(L;\omega_{\pi})$ given by
$f\mapsto h_f$
where $$h_f(l)=\int_{U^-_{M,w}}\int_{U^-_{M,w^{-1}}}f(x^-\dot{w}lu^-)\psi^{-1}(x^-u^-)dx^-du^-.$$
We have the following descent formula for partial Bessel integrals:

\begin{prop}\label{prop6.3} Suppose $H\subset L\subset M$ are Levi subgroups of $M$, and $\varphi$ is invariant under the $\Theta_M$-twisted conjugate action of some open compact subgroups of $U_M(F)$. Let $w=w^M_L\in B(M)$. Then for any fixed $m\in C_r(\dot{w}^M_H)$, decompose $m$ as $m=u_1^- \dot{w} l u_2^-$, with $u_1^-\in U_{M,w^{-1}}^-$, $u_2^-\in U_{M,w}^-$, $l\in L$, we have 
$$B^M_\varphi(m,f)=\psi(u_1^-u_2^-)B^L_\varphi(u_1^-, u_2^-, l.h_f)$$
where $$B_\varphi^L(u_1^-,u_2^-, m,h_f)=\int_{U_{L,l}\cap n_0U_{L,l}n_0^{-1}\backslash U_L}\int_{U_L}h_f(x'lu')\varphi(\Theta_M((n_0{u'})^{-1}\dot{w}l'u')\psi^{-1}(x')\psi^{-1}(u')dx'du'$$
in which $n_0=\tilde{u}_1^-(u_2^-)^{-1}\in N_L$, $\tilde{u}_1^-=\tau^M_L\circ \Theta_M\circ \tau^M_L((u_1^-)^{-1})$, and $l=zl'$ with $z\in Z_M$.
\end{prop}

\begin{proof} By Lemma \ref{lem6.2}, $U^{\Theta_{M}}_{M,m}=(\tilde{u}_1^-)^{-1}U^{\Theta_L}_{L,l}\tilde{u}_1^-\cap (u_2^-)^{-1} U^{\Theta_L}_{L,l}u_2^-$. Denote $ n_1=(\tilde{u}_1^-)^{-1}, n_2=(u_2^-)^{-1},$ and $n_0=\tilde{u}_1^-(u_2^-)^{-1}$, then $n_1, n_2\in U_{M,w}^-=N_L\subset U_M$, so $U_M=n_1^{-1}U_Mn_1$. Hence if we decompose $x=x^-x^+$ with $x^-\in U_{M,w^{-1}}^-$, $x^+\in U_{M,w^{-1}}^+$, make a change of variable $u\mapsto n_1 un_1^{-1}$, and then decompose $u=u^+u^-$ with $u^+\in U_{M,w}^+=U_L$, $u^-\in U_{M,w}^-=N_L$ in the integral
$$B^M_\varphi(m,f)=\int_{U_{M,m}^{\Theta_M}\backslash U_M}\int_{U_M}f(xmu)\varphi(\Theta_M(u^{-1})m'u)\psi^{-1}(xu)dxdu,$$ 
we obtain
$$B^M_\varphi(m,f)=\int_{n_1(U^{\Theta_L}_{L,l}\cap n_0U^{\Theta_L}_{L,l}n_0^{-1})n_1^{-1}\backslash n_1U_L n_1^{-1}}\int_{U_{M, w^{-1}}^+}\int_{U_{M,w}^-}\int_{U_{M,w^{-1}}^-}f(x^-x^+u_1^-\dot{w}lu_2^-n_1 u^+u^-n_1^{-1})$$
$$\varphi(\Theta_M(n_1(u^-)^{-1}(u^+)^{-1}n_1^{-1})u_1^-\dot{w}l'u_2^-n_1u^+u^-n_1^{-1})\psi^{-1}(x^+x^-)\psi^{-1}(n_1u^+u^-n_1^{-1})dx^-du^-dx^+du^+.$$
Since $n_1^{-1}=\tilde{u}_1^-=\tau^M_L\circ \Theta_M\circ \tau^M_L((u_1^-)^{-1})$, we see that $\Theta_M(n_1^{-1})=(\Theta_M\circ\tau^M_L)^2((u_1^-)^{-1})=(u_1^-)^{-1}$. Let $y^-=x^-x^+u_1^-(x^+)^{-1}$, $x'=\dot{w}^{-1}x^+\dot{w}$, $v^-=(u^+)^{-1}u_2^-n_1u^+u^-n_1^{-1}$, and $u'=u^+$, then $y^-\in U_{M,w^{-1}}^-$, $x'\in U_{M,w}^+=U_L$, and $v^-\in U_{M,w}^-=N_L$. By the compatibility of $\psi$ with $\dot{w}$ by \cite[Proposition 5.1]{CST17}, we have $\psi(x')=\psi(x^+)$. Therefore we can rewrite the above integral as
$$B^M_{\varphi}(m,f)=\psi(u_1^-u_2^-)\int_{U^{\Theta_L}_{L,l}\cap n_0U^{\Theta_L}_{L,l}n_0^{-1}\backslash U_L}\int_{U_L}\int_{U_{M,w^{-1}}^-}\int_{U_{M,w}^-}f(y^-\dot{w}x'lu'v^-)$$$$\varphi(\Theta_M((n_0u'v^-)^{-1})\dot{w}l'u'v^-)\psi^{-1}(v^-y^-x'u')dv^-dy^-dx'du'.$$
Since $f\in C^\infty_c(\Omega_{w};\omega_\pi)$, $\Omega_w=U_{M,w^{-1}}^-\times \dot{w}L\times U_{M,w}^-$, $U_{M,w^{-1}}^-$ and $U_{M,w}^-$ are closed in $\Omega_w$, so there exists compact subgroups $U_1$ of $U_{M,w}^-$ and $U_2$ of $U_{M,w^{-1}}^-$ such that $f(y^-\dot{w}x'lu'v^-)\neq 0$ implies that $y^-\in U_1$ and $v^-\in U_2$. On the other hand, by our assumption on $\varphi$, $\varphi$ is invariant under the $\Theta_M$-twisted conjugate action by some open compact subgroup $U_0$ of $U_M(F)$. We shrink $U_2$ if necessary so that $v^-$ lies in $U_0$. It follows that
$$\varphi(\Theta_M(n_0u'v^-)^{-1}\dot{w}l'u'v^-)=\varphi(\Theta_M(n_0u')^{-1}\dot{w}l'u').$$ 
So we obtain that
$$B^M_\varphi(m,f)=\psi(u_1^-u_2^-)\int_{U^{\Theta_L}_{L,l}\cap n_0U^{\Theta_L}_{L,l}n_0^{-1}\backslash U_L}\int_{U_L}\int_{U_{M,w^{-1}}^-}\int_{U_{M,w}^-}f(y^-\dot{w}x'lu'v^-)\psi^{-1}(v^-y^-)dv^-dy^-$$$$\varphi(\Theta_M((n_0u')^{-1})\dot{w}l'u')\psi^{-1}(x'u')dx'du'$$
$$=\psi(u_1^-u_2^-)\int_{U^{\Theta_L}_{L,l}\cap n_0U^{\Theta_L}_{L,l}n_0^{-1}\backslash U_L}\int_{U_L}h_f(x'lu')\varphi(\Theta_M((n_0u')^{-1})\dot{w}l'u')\psi^{-1}(x'u')dx'du'=\psi(u_1^-u_2^-)B^L_\varphi(u_1^-,u_2^-,l, h_f).$$ To show the last integral is well-defined, note that if $v''=n_0u''n_0^{-1}$, with $u'',v''\in U^{\Theta_L}_{L,l}$, then $$l'=\Theta_L((v'')^{-1})l'v''=\Theta_L(({u''}^{-1}n_0u''n_0^{-1})^{-1})\Theta_L((u'')^{-1})l'u''({u''}^{-1}n_0u''n_0^{-1})$$$$=\tau^M_L\circ\Theta_M(({u''}^{-1}n_0u''n_0^{-1})^{-1})l'({u''}^{-1}n_0u''n_0^{-1}),$$
so $${l'}^{-1}\tau^M_L\circ\Theta_M({u''}^{-1}n_0u''n_0^{-1})l'={u''}^{-1}n_0u''n_0^{-1}.$$
Since $\Theta_M(U_M)=U_M$, $U_M=U_LN_L$, and $\tau^M_L(N_L)=\overline{N_L}$, the opposite of $N_L$.
The left hand side of the above equality lies in $\overline{N_L}$, while the right hand side lies in $N_L$, as $N_L\cap\overline{N_L}=\{1\}$ it forces that ${u''}^{-1}n_0u''n_0^{-1}=1$, i.e. $u''n_0=n_0u''$, and $v''=u''$. It follows that 
$$\Theta_M((n_0v''u')^{-1})\dot{w}l'v''u'=\Theta_M((u'' n_0u')^{-1})\dot{w}\Theta_L(u'')(\Theta_L(u'')l'u'')u'$$
$$=\Theta_M((n_0u')^{-1})\Theta_M((u'')^{-1})\dot{w}\Theta_L(u'')l'u'=\Theta_M((n_0u')^{-1})\dot{w}l'u'$$ where we used $\dot{w}^{-1}\Theta_M(u'')\dot{w}=\tau^M_L\circ\Theta_M(u'')=\Theta_L(u'')$.
\end{proof}
\begin{remark}\label{rem6.4}
    If we impose a stronger condition on the cut-off function $\varphi$, i.e, we require that $\varphi(\Theta_M(x_1^{-1})mx_2)=\varphi(m)$ for $x_1,x_2\in U_0$, then we can show that $$B^M_\varphi(m,f)=\psi(u_1^-u_2^-)\mathrm{Vol}(U^{\Theta_L}_{L,l}\cap n_0U^{\Theta_L}_{L,l}n_0^{-1}\backslash U^{\Theta_L}_{L,l})B^L_{\varphi_{\dot{w}}}(l,h_f)$$ for all $m\in C_r(\dot{w})$ where $\varphi_{\dot{w}}(\cdot)=\varphi(\dot{w}\cdot)$ is the left translation of $\varphi$ by $\dot{w}$. But notice that in general the subset consists of the image $n\mapsto m$ via $\dot{w}_0^{-1}n=mn'\bar{n}$ is not invariant under the action $m\mapsto \Theta_{M}(x_1^{-1})mx_2$ for $x_1,x_2\in U_M$. In particular, one checks easily from our calculation of the Bruhat decomposition in section 3 that in our case the set of elements of the form $m(X,Y)$ is not invariant under this action. So we do not use this formula, although it is of a nicer form. 
\end{remark}
\begin{remark}\label{rem6.5} We observe easily that
    in particular, if $n_0=1$, i.e. $\tilde{u}_1^-=u_2^-$, we have
    $$B^M_\varphi(m, f)=B^L_{\varphi_{\dot{w}}}(l, h_f).$$
\end{remark}

\subsection{Asymptotic expansion and uniform smoothness}
The main idea to prove stability is that once we relate our local coefficient formula with the Mellin transform of the partial Bessel integrals, there are two parts in the asymptotic expansion formula of the partial Bessel integrals, one depends only on the central character of $\pi$, the other is certain uniform smooth function. Therefore under highly ramified twists, the second part vanishes.

To obtain the asymptotic expansion of partial Bessel integrals, we study its boundary behavior on $\Omega_w$ with $w\in B(M)$, i.e. on Bruhat cells of smallest possible dimension, and do induction on Bessel distance. This is a standard procedure as in \cite{CST17} or \cite{She23}. We begin with the small cell $C(e)=B_M$ of $M$. First note that 
$L_e=M, A_e=Z_M$, $U_{M,e}^+=U_M$, and 
$\Omega_e=M$. Take $\dot{e}=I$, the identity matrix in $M$. Recall that $\pi$ is $\psi$-generic supercuspidal and we pick an matrix coefficient $f\in C^{\infty}_c(M;\omega_{\pi})$, normalized with $W_f(e)=1.$ Then we have

\begin{prop}\label{prop6.6}
Fix an auxiliary function $f_0\in C^{\infty}_c(M;\omega_{\pi})$ with $W_{f_0}(e)=1$. Then for each $f\in C_c^{\infty}(M;\omega_{\pi})$ with $W_f(e)=1$ and $w'\in B(M)$ with $d_B(e,w')=1$, there exists a function $f_{w'}\in C^{\infty}_c(\Omega_{w'};\omega_{\pi})$ such that for any $w\in B(M)$ and $m\in M$, we have 
$$B^M_{\varphi}(m,f)=B^M_\varphi(m,f_1)+\sum_{w'\in B(M), d_B(w',e)=1}B^M_{\varphi}(m,f_{w'})$$ where $f_1(m):=\sum_{m=m_1c}f_0(m_1)B^M(\dot{e}c, f)=\sum_{m=m_1c}f_0(m_1)\omega_{\pi}(c)$, the sum runs over all possible decompositions $m=m_1c$ with $m_1\in M_{der}$, and $c\in A_e=Z_M$.
\end{prop}
\begin{proof}

 Decompose $M=M_{der}A_e=M_{der}Z_M$, where $M_{der}$ is the derived group of $M$. Write $m=m_1c$ with $m_1\in M_{der}$ and $c\in Z_M$, then there are only finitely many such decompositions which are indexed by elements in the transverse torus $A^e_e=M_{der}\cap Z_M$. 
Define
$$f_1(m):=\sum_{m=m_1c}f_0(m_1)B^M(\dot{e}c, f)=\sum_{m=m_1c}f_0(m_1)\omega_{\pi}(c)$$ Then $f_1(m)\in C_c^{\infty}(M;\omega_{\pi})$. For $a\in A_e=Z_M$, $a'=e$, therefore
$$B^M_\varphi(\dot{e}a, f_1)=\omega_\pi(a)\int_{U_{M,e}^{\Theta_M}\backslash U_M}\int_{U_M} f_1(xu)\varphi(\Theta_M(u^{-1}) u)\psi^{-1}(xu)dx du=\omega_\pi(a)W_{f_1}(e)\int_{U_{M,e}^{\Theta_M}\backslash U_M} \varphi(\Theta_M(u^{-1})u)du.$$ 
On the other hand, $W_{f_1}(e)=\int_U f_1(x)\psi^{-1}(x)dx$, while $x\in U_M\subset M_{der}$ so $f_1(x)=f_0(x)$, thus $W_{f_1}(e)=W_{f_0}(e)=1$, hence we obtain that
$B^M_\varphi(\dot{e}a, f_1)=B^M_\varphi(\dot{e}a, f)$ for all $a\in A_e$. Therefore $B^M_\varphi(\dot{e}a, f-f_1)=0$ for all $a\in A_e$. Apply \cite[Lemma 5.13]{CST17}, there exists $f_2'\in C_c^\infty(\Omega_e^\circ,;\omega_\pi)$, such that $B^M_\varphi(m,f-f_1)=B^M_\varphi(m, f_2')$ for all $m\in M$, where $\Omega_e^\circ=\Omega_e-C(e)=M-B_M$. Let $\Omega_1=\cup_{w\in B(M), w\neq e}\Omega_w=\cup_{w\in B(M), d_B(w,e)=1}\Omega_w$, and $\Omega_0=M-C(e)=\Omega^\circ_e$. By \cite[Lemma 5.14]{CST17}, there exists $f_2\in C^\infty_c(\Omega_1;\omega_\pi)$, so that $B^M_\varphi(m, f_2)=B^M_\varphi(m, f_2')=B^N_\varphi(m, f-f_1)$, for a sufficiently large $\varphi$ depending only on $f_1$. Then by a partition of unity argument, for each $w'\in B(M)$, $d_B(w',e)=1$, we can find $f_{w'}\in C^\infty_c(\Omega_{w'};\omega_\pi)$, such that $f_2=\sum_{w'\in B(M), d_B(w',e)=1}f_{w'}$. 
\end{proof}

We would like to do the same type of analysis on each $B^M_\varphi(m,f_{w'})$ and obtain an asymptotic expansion formula for $B^M_\varphi(m,f)$ indexed by Weyl group elements that supports Bessel functions, and for all $m$ lying in the relevant part of Bruhat cells. To be precise we will show that
\begin{prop}\label{prop6.7}
    Fix an auxiliary function $f_0\in C^{\infty}_c(M;\omega_{\pi})$ with $W_{f_0}(e)=1$. Then for each $f\in C_c^{\infty}(M;\omega_{\pi})$ with $W_f(e)=1$ and $w'\in B(M)$ with $d_B(e,w')\ge 1$, there exists a function $f_{w'}\in C^{\infty}_c(\Omega_{w'};\omega_{\pi})$ such that for any $m=u_1\dot{w}au_2\in C_r(\dot{w})$ with $\dot{w}=\dot{w}^M_L\in B(M)$, we have 
$$B^M_{\varphi}(m,f)=\sum_{a=bc}\omega_\pi(c)B^M_\varphi(u_1 \dot{w}bu_2, f_0)+\sum_{w'\in B(M), d_B(w',e)\ge 1}B^M_{\varphi}(m,f_{w'})$$ 
where $a=bc$ runs over the possible decompositions of $a\in A_w$ with $b\in A^e_w$ and $c\in A_e=Z_M.$
\end{prop}
\begin{proof}

Note that in particular
if $m=u_1\dot{w}au_2\in C_r(\dot{w})$, by Lemma \ref{lem6.1}, $U^{\Theta_M}_{M,m}\subset u_2^{-1}U_{M.w}^+u_2$, we can write $U_M=u_2^{-1}U_{M,w}^+ U_{M,w}^- u_2^{-1}=u_2^{-1}U_{M,w}^+ u_2(u_2^{-1}U_{M,w}^-u_2)$. Let $u=u'(u_2^{-1}u^-u_2)$ where $u'=u_2^{-1}u^+u_2$ with $u^+\in U_{M,w}^+$, $u^-\in U_{M,w}^-$, we have
 $$B^M_{\varphi}(m,f_1)=\int_{U_{M,m}^{\Theta_M}\backslash u_2^{-1}U_{M,w}^+u_2}\int_{U_{M,w}^-}\int_{U_M}f_1(xu_1\dot{w}au_2u'u_2^{-1}u^-u_2)$$$$\cdot\varphi(\Theta_M(u_2^{-1}(u^-)^{-1}u_2 {u'}^{-1})u_1\dot{w} a' u_2 u' u_2^{-1}u^-u_2)\psi^{-1}(xu' u_2^{-1}u^-u_2)dxdu^-du'.$$
 As $u^+=u_2u'u_2^{-1}\in U_{M,w}^+$ and $a\in A_{w}$, we have $xu_1\dot{w}au_2^{-1}u'u_2^{-1}u^-u_2=xu_1(\dot{w}u_2u'u_2^{-1}\dot{w}^{-1})\dot{w}u^-u_2$. Let $x'=xu_1(\dot{w}u_2u'u_2^{-1}\dot{w}^{-1})\in U_{M,w}^+$, $v^-=u^-u_2\in U_{M,w}^-= U_{M,w}^-$, and by compatibility of $\psi$ with $\dot{w}$, we obtain that
 $$B^M_\varphi(m,f_1)=\psi(u_1u_2)\int_{U_{M,m}^{\Theta_M}\backslash u_2^{-1}U_{M,w}^+u_2}\int_{U_{M,w}^-}\int_{U_M}f_1(x'\dot{w} av^-)$$$$\cdot\varphi(\Theta_M((v^-)^{-1}u_2 {u'}^{-1})u_1\dot{w} a' u_2 u' u_2^{-1}v^-)\psi^{-1}(x' v^-)dx'dv^-du'$$
 So by the construction of $f_1$, we need to decompose $x'\dot{w}av^-=m_1c$ with $m_1\in M_{der}$ and $c\in A_e=Z_M$, which is equivalent to $ac^{-1}=\dot{w}^{-1}(x')^{-1}m_1(v^-)^{-1}$. Since we pick $\dot{w}\in M_{der}$, and $x', v^-\in U_M\subset M_{der}$, this is saying that $b:=ac^{-1}\in M_{der}\cap A_{w}=A_{w}^e$. It follows that $f_1(x'\dot{w}av^-)=\sum_{a=bc}f_0(x'\dot{w}bv^-)\omega_\pi(c)$, and consequently
$B^M_\varphi(m,f_1)=\sum_{a=bc}\omega_\pi(c)B^M_\varphi(u_1 \dot{w}bu_2, f_0)$ where $a=bc$ runs over the possible decompositions of $a\in A_w$ with $b\in A^e_w$ and $c\in A_e=Z_M.$

We continue to obtain expansions for $B^M_\varphi(m,f_{w'})$. Suppose $w'=w^M_L\in B(M)$. Let $h=h_{f_{w'}}\in C^\infty_c(L;\omega_\pi)$ be the image of $f_{w'}$ under $C^\infty_c(\Omega_{w'}; \omega_\pi)\twoheadrightarrow C^{\infty}_c(L;\omega_\pi)$. Pick $h_0\in C^\infty_c(L;\omega_\pi)$ normalized so that $B^L_{\varphi_{\dot{w}'}}(\dot{e}, h_0)=\frac{1}{\vert Z_L\cap A^{w'}_{w'}\vert}$, and $B^M(b, h_0)=0$ for $b\in A^{w'}_{w'}$ but $b\notin Z_M\cap A^{w'}_{w'}$. Similar to the construction of $f_1$ from $f$, let $h_1(l)=\sum_{l=l_1c}h_0(l_1)B^L(c,h)$, where the decomposition is taken over all possible $l=l_1c$ where $l_1\in L_{der}$ and $c\in A_{w'}=Z_L$. By \cite[Proposition 5.4]{CST17}, 
$B^L_{\varphi_{\dot{w}'}}(a, h_1)=B^L_{\varphi_{\dot{w}'}}(a, h)$, for all $a\in A_{w'}=Z_L$. Choose $f_1$ which maps to $h_1$ under $C^{\infty}_c(\Omega_{w'}, \omega_{\pi}) \twoheadrightarrow C^{\infty}_c(L;\omega_{\pi})$, apply Proposition \ref{prop6.3} to case when $H=L$, and $m=w'a$, $a\in A_{w'}=Z_L$, we get
$B^M_{\varphi}(w'a,f_1)=B^L_{\varphi_{\dot{w}'}}(a, h_1)=B^L_{\varphi_{\dot{w}'}}(a,h)=B^M_\varphi(w'a,f_{w'})$. So $B^M_\varphi(w'a, f_{w'}-f_1)=0$, for all $a\in A_{w}$. Apply Apply \cite[Lemma 5.13, 5.14]{CST17} again, together with a partition of unity argument, there exists $f_{w',w''}\in C^\infty_c(\Omega_{w''};\omega_\pi)$ for $w'<w''\in B(M)$ such that $d_B(w',w'')=1$, and
$$B^M_{\varphi}(m, f_{w'})=B^M_{\varphi}(m,f_1)+\sum_{w''\in B(M), w''>w', d_B(w'',w')=1}B^M_{\varphi}(m,f_{w''}).$$ Proceed by repeating this process and do induction on the Bessel distance, and combine with the above argument on $B^M_\varphi(m,f_1)$, we obtain the asymptotic expansion formula as desired.
\end{proof}

Note that the first sum in Proposition \ref{prop6.7} depends only on $\omega_\pi$ and a fixed auxiliary function $f_0$. The next important step is to show that each $B^M_\varphi(m,f_{w'})$ satisfy some uniform smooth property. Recall that for $h\in C^\infty_c(\Omega_{w'};\omega_\pi)$ with $w'=w^M_L\in B(M)$, we constructed $h_1(l)=\sum_{l=l_1c}h_0(l_1)B^L(c,h)$. Suppose $m=u_1^-\dot{w}' l u_2^-$ and write $l=zl'$ with $z\in Z_M$, $n_0=\tilde{u}_1^-(u_2^-)^{-1}$, where $\tilde{u}_1^-=\tau^M_L\circ \Theta_M\circ \tau^M_L((u_1^-)^{-1})$, we have $B^L_\varphi(u_1^-,u_2^-, l,h_1)=\omega_\pi(z)B^L_\varphi(u_1^-, u_2^-, l', h_1)$. Since 
$$B^L_\varphi(u_1^-,u_2^-,l',h_1)=\int_{U^{\Theta_L}_{L,l}\cap n_0U^{\Theta_L}_{L,l}n_0^{-1}\backslash U_L}\int_{U_L}h_1(xl'u)\varphi(\Theta_M((n_0u)^{-1})\dot{w}'l'u)\psi^{-1}(xu)dxdu,$$ we need to write $xl'u=l_1c$ with $l_1\in L_{der}$, $c\in Z_L$. Suppose $l=v_1\dot{w}av_2\in C^L(\dot{w})$ a Bruhat cell in $L$ with $\dot{w}\in L_{der}$, then $l'=v_1\dot{w}a'v_2$ where $a=za'$. As $v_1, v_2,\dot{w}\in L_{der}$, so to write $xl'u=xv_1\dot{w}a'v_2=l_1c$ is equivalent to write $a'=bc$ with $b\in A\cap L_{der}$ and $c\in Z_L$.
Let $b=z_bb'$, $c=z_cc'$ with  $z_b,z_c\in Z_M$, $b'\in A'$, $c'\in Z_L'$, then $a'=z_bz_cb'c'\Rightarrow z_bz_c=1$, and $a'=b'c'$. Let $l_{b}=v_1\dot{w} b v_2$,  and similarly $l_{b'}= v_1\dot{w} b' v_2$, then
$h_1(xl'u)=\sum_{a'=b'c'}h_0(xl_{b'}u)B^L(c' , h)$. Since $U_{L,l'}^{\Theta_L}=\{u\in U_L: \Theta_L(u^{-1})l'u=l'\}=\{u\in U_L: \Theta_L(u^{-1})l_{b'}c'u=l_{b'}c'\}=\{u\in U_L:\Theta_L(u^{-1})l_{b'}uc'=l_{b'}c'\}=\{u\in U_L:\Theta_L(u^{-1})l_{b'}u=l_{b'}\}=U_{L, l_{b'}}^{\Theta_L}$, we obtain that
$$B^L_\varphi(u_1^-,u_2^-,l', h_1)=\sum_{a'=b'c'}B^L(c', h)\int_{U_{L,l'}^{\Theta_L}\cap n_0 U^{\Theta_L}_{L,l'}n_0^{-1}\backslash U_L}\int_{U_L}h_0(xl_{b'}u)\varphi(\Theta_M((n_0u)^{-1})l_{b'}uc')\psi^{-1}(xu)dxdu$$$$=\sum_{a'=b'c'}B^L(c',h)B^L_{\varphi^{c'}}(u_1^-,u_2^-, l_{b'}, h_0)$$ where $\varphi^{c'}(\cdot):=\varphi(\cdot c')$ is the right shift of $\varphi$ by $c'$.

A smooth function on a p-adic group is \textbf{uniform smooth} if it is uniformly locally constant, i.e,, there exists an open compact subgroup $K_0$ such that the function is constant on $aK_0$ for any $a$ in the group. Recall that $A$ is the maximal split torus of $M$. We say a unipotent element $u\in U_M$ is \textbf{rationally parameterized} by $(a,d)\in A\times \mathbb{A}^k$ if all the entries in the 1-parameter subgroups corresponding to $u$ are rational functions of $(a,d)$, here $\mathbb{A}^k$ is the affine space of dimension $k$, i.e., if one has $u=\prod_{\alpha}u_\alpha(x_\alpha)$ with $x_\alpha\in \mathbb{G}_a$, then each $x_\alpha$ is a rational function of $(a,d)$.
\begin{prop} \label{prop6.8}Let $H\subset L\subset M$ be Levi subgroups of $M$, suppose $m=u_1^-\dot{w}'lu_2^-=m(a,d)=u_1(a,d)\dot{w}au_2(a,d)\in C_r(\dot{w})$ is rationally parameterized by $(a,d)\in A^{w'}_w A_{w'}\times \mathbb{A}^k\subset A\times \mathbb{A}^k$ for some $k\ge 0$ with $w=w^M_H, w'=w^M_L\in B(M)$. Assume that the rational functions that parameterize $u_1$ and $u_2$ have no singularities on $A\times \mathbb{A}^k$. Then if one writes $a=bc$ with $b\in A^{w'}_{w}$, $c\in A_{w'}$ and $c=zc'$ with $z\in Z_M$, $c\in A'_{w'}$, where $A_{w'}=Z_MA'_{w'}$, we have that
$$B^M_\varphi(m,f_1)=\omega_\pi(z)\psi(u_1^-(bc'z,d),u_2^-(bc'z,d))\sum_{a'=b'c'}B^L(c', h)B^L_{\varphi^{c'}}(u_1^-(bc'z), u_2^-(bc'z),l_{b'}(bc'z,d), h_0)$$ is uniformly smooth as a function of $c'\in A'_{w'}$, where $f_1\mapsto h_1$ via $C^{\infty}_c(\Omega_{w'}; \omega_{\pi}) \twoheadrightarrow C^{\infty}_c(L;\omega_{\pi})$.
\end{prop}
\begin{proof}
    We have
    $B^M_\varphi(m,f_1)=\psi(u_1^-, u_2^-)B^L_\varphi(u_1^-,u_2^-, l, h_1)=\psi(u_1^-,u_2^-)\sum_{a=bc}B^L(c, h)B^L_{\varphi^{c}}(u_1^-(bc,d),u_2^-(bc,d), l_{b}, h_0)$ by Proposition 6.3 and the above argument. By fixing a decomposition $a=bc$ with $b\in A_{w_M}^{w'}$, $c\in A_{w'}$, all such decompositions are of the form $a=(b\xi^{-1})(\xi c)$ with $\xi\in A_{w}^{w'}\cap A_{w'}=A_{w'}^{w'}=Z_L\cap L_{der}$, a finite set, so $\vert\xi\vert=1$. Therefore $$B^L_\varphi(u_1^-,u_2^-, l, h_1)=\sum_{\xi\in A_{w'}^{w'}}B^L(\xi c, h)B^L_{\varphi^{\xi c}}(u_1^-(bc,d), u_2^-(bc,d), l_{b\xi^{-1}}, h_0)$$ where $l_{b\xi^{-1}}=\dot{w}'^{-1}u_1^+\dot{w}'\dot{w}^L_H b\xi^{-1}u_2^+$. As $\varphi$ depends only on the absolute value, $\varphi^{\xi c}=\varphi^{c}$. On the other hand,
    $$B^L(\xi c, h)=\int_{U_M(F)}h(x\xi c)\psi^{-1}(x)dx=\omega_{\pi}(\xi_1 z)\int_{U_M(F)}h(x\xi' c')\psi^{-1}(x)dx$$ where we write $\xi=\xi_1\xi'$ with $\xi_1\in Z_M$ and $\xi'=\xi \xi_1^{-1}$. As $h\in C^\infty_c(L;\omega_\pi)$ is smooth of compact support modulo $Z_M$, and the small Bruhat cell $C^L(e_L)=AU_L=Z_M A' U_L$ of $L$ is closed in $L$, and $Z_L'\subset A'$ is closed, there exists compact subsets $K_1\subset Z_L'$ and $U_1\subset U_M$ such that $h(x\xi'c')\neq 0\Rightarrow \xi'c'\in K_1$ and $x\in U_1$. Writing $c=zc'$, we see that 
    $B^L_\varphi(u_1^-, u_2^-, l, h_1)=\omega_\pi(z)\sum_{\xi\in A_{w'}^{w'}}B^L(\xi c', h)B^L_{\varphi^{ c}}( u_1^-(bc'z,d), u_2^-(bc'z,d), l_{b\xi^{-1}}, h_0)$ does not vanish unless $c'\in K:=\cup_{\xi'}\xi'^{-1}K_1$, a compact subset. From this we also see that the support of $h$ in $c'$ is compact and independent of the decomposition $a=bc$. 

    Moreover, since the support of $h$ in $c'$ is compact, and $h$ is smooth, there exists an open compact subgroup $C_0\subset Z_L'$ such that $h(x\xi'c'c_0)=h(x\xi'c')$ for all $x\in U_1$, $c'\in Z_L'$, and $c_0\in C_0$. By our assumption the entries in the 1-parameter subgroups corresponding to $u_i(a,d)$, and hence $u_i^\pm(a, d)=u_i^\pm(bc'z, d)$ $(i=1,2)$ are rational functions in $(b,c',z,d)$ without singularities, it follows that $u_i^\pm$ and hence also $\psi(u_1^-, u_2^-)$ are smooth functions of $(b,c',z,d)$ with no singularities, and the quotient space $U^{\Theta_L}_{L,l}\cap n_0U^{\Theta_L}_{L,l}n_0^{-1}\backslash U_L$ is also smoothly parameterized by $(b,c', z,d)$. Therefore $B^M_\varphi(m(bc'z,d), f_1)$ is zero when $c'\notin K$. So for each fixed $z, b, d$, we simultaneously choose $C_0$, depending on $z, b, d$, so that $u_i^{\pm}(zbc'c_0, d)=u_i^{\pm}(zbc', d)$ for all $c_0\in C_0, c'\in K$. Finally, We shrink $C_0$ so that $C_0\subset Z_L'(\mathcal{O}_F)$ if necessary, then $\varphi^{cc_0}=\varphi^{c}$ for all $c_0\in C_0$ and $c\in A_{w'}=Z_L$. Consequently, there exists an open compact subgroup $C_0\subset Z_L'$ such that
    $$B^M_\varphi(m(z bc'c_0, d), f_1)=B^M_\varphi(m(zbc', d), f_1)$$ for all $a=bc$ and $c_0\in C_0\subset Z_L'=A_{w'}'$, i.e. $B^M_\varphi(m, f_1)$ is uniformly smooth in $c'\in A_{w'}'$.
\end{proof}
\subsection{The final local coefficient formula and the separation of the toric part} We will relate the partial Bessel functions defined in section 5 with partial Bessel integrals introduced in section 6, so that the local coefficient formula in Proposition \ref{prop5.2} can be restated in terms of Mellin transforms of partial Bessel integrals. With nice asymptotic expansion formulas in Proposition \ref{prop6.7} and uniform smoothness in Proposition \ref{prop6.8}, we will be able to show our desired analytic stability. 

For our applications, we set the involution $\Theta_M: m\mapsto \dot{w}_0^{-1}m\dot{w}$, and it is easily verified that it satisfies the assumption that it preserves the $F$-splitting and is compatible with Levi subgroups and the generic character $\psi$. Recall that for a matrix $X$ of size $k\times l$, we defined
$\varphi_\kappa(X)=\begin{cases}
    1 & \vert X_{i,j}\vert\leq q^{((k-i)+(l-j)+1)\kappa} \\
    0 & \mathrm{\ otherwise \ }\\
\end{cases}$, and $\overline{N}_{0,\kappa}=\{\bar{n}(X,Y): \varphi_\kappa(\varpi^{-(d+g)}X)\cdot\varphi_\kappa(\varpi^{-2(d+g)}Y)=1\}$, where $d=\mathrm{Cond}(\psi)$ and $g=\mathrm{Cond}(\omega_\pi^{-1}(w_0(\omega_\pi)))$. We also denoted by $\varphi_{\overline{N}_{0,\kappa}}$ the characteristic function of $\overline{N}_{0,\kappa}$. Suppose $n(X,Y)\mapsto m(X,Y)$ via the Bruhat decomposition $\dot{w}_0^{-1}n=mn'\bar{n}$. Define $$\varphi(m(X,Y)):=\varphi_\kappa((Y^{-1}X)_{r1}X)\varphi_\kappa(((Y^{-1}X)_{r1})^2 Y)$$ where $(Y^{-1}X)_{r1}$ is the $(r,1)$-th entry of $Y^{-1}X$. We will show that $\varphi$ satisfies the assumption that it is invariant under the $\Theta_M$-twisted conjugate action by some open compact subgroup $U_0$ of $U_M(F)$. Consequently, the assumptions on $\Theta_M$ and $\varphi$ in the beginning of section 6 are both satisfied and we have

\begin{prop}\label{prop6.9} Let $m(X', Z')$ be the image of $n(X',Z')$ under the map $n\mapsto m$ via $\dot{w}_0^{-1}n=mn'\bar{n}$. Then
$$j_{\overline{N}_{0,\kappa},\pi,w}(n(X',Z'))=B^M_\varphi(m(X',Z'), f)$$ where $\varphi_{\overline{N}_{0,\kappa}}$ and $\varphi$ are defined as above and $f\in C^\infty_c(M(F);\omega_{\pi})$ such that $W=W_f$, normalized so that $W_f(e)=1$. 
    
\end{prop}

\begin{proof}
Given the Bruhat decomposition $\dot{w}_0^{-1}n=mn'\bar{n}$, for $z=\alpha^\vee(t)\in Z_M^0$, if $u\mapsto zu^{-1}nuz^{-1}$, then $m\mapsto \Theta_M(zu^{-1})muz^{-1}=\Theta_M(u)^{-1}\Theta_M(z)z^{-1}m u=\Theta_M(u^{-1})\alpha^\vee(t^{-2})mu$, as $\Theta_M(z), z\in Z_M$. 

In our cases, where $n=n(X,Y)$, then  $m=m(X,Y)=\mathrm{diag}\{m_1(X,Y), m_2(X,Y), \theta_r(m_1(X,Y))\}$. So $$\Theta_M(m(X,Y))=\mathrm{diag}(\theta_r(m_1(X,Y)),m_2(X,Y), m_1(X,Y))$$ It shows that $\Theta_M$ induces $\theta_r$. Moreover, for $z=\alpha^\vee(t)\in Z_M^0$, $u=\mathrm{diag}\{u_1, u_2, \theta_r(u_1)\}\in U_M$, we have
$$zu^{-1}\bar{n}(X,Y)u z^{-1}=\bar{n}(u_1^{-1}(tX)u_2, u_1^{-1}(t^2Y)\theta_r(u_1))=u^{-1}\bar{n}(tX, t^2Y)u.$$ Hence the action $\bar{n}(X,Y)\mapsto zu^{-1}\bar{n}(X,Y)uz^{-1}$ is given by $X\mapsto u_1^{-1}tXu_2, Y\mapsto u_1^{-1}t^2Y\theta_r(u_1)$. On the other hand, $$\Theta_M(zu^{-1})muz^{-1}=\Theta_M(u)^{-1}\alpha^\vee(t^{-2})m(X,Y)=\begin{bmatrix}
    \theta_r(u_1^{-1}) & & \\
    & u_2^{-1} &\\
    & & u_1^{-1}\\
\end{bmatrix} \begin{bmatrix}
    t^{-2}\theta_r(Y) & & \\
    & (I_{2m}-J'_{2m}{^tX}{^tY}^{-1}J_r X)^{-1} & \\
    & & t^2Y\\
\end{bmatrix}$$$$\cdot \begin{bmatrix}
    u_1 & & \\
    & u_2& \\
    & & \theta_r(u_1)\\
\end{bmatrix}=\begin{bmatrix}
    \theta_r(u_1^{-1})\theta_r(t^2Y)u_1 &  & \\
        &  (I_{2m}-J'_{2m}{^t(u_1^{-1}tXu_2)}{^t(u_1^{-1}t^2Y\theta_r(u_1))}J_r(u_1^{-1}tXu_2))^{-1} &  \\
        &    & u_1^{-1}t^2Y\theta_r(u_1) \\
\end{bmatrix}.$$ This verifies that for $z=\alpha^\vee(t)$, $u=\mathrm{diag}(u_1,u_2,\theta_r(u_1))\in U_M$, we have
$$\Theta_M(zu^{-1})m(X,Y)uz^{-1}=m(u_1^{-1}tXu_2, u_1^{-1}t^2Y\theta_r(u_1)).$$
 Let $U_{0,\kappa}=U_{1,\kappa}\times U_{2,\kappa}$ with $U_{1,\kappa}\subset U_{\mathrm{GL}_r}(F)$ and $U_{2,\kappa}\subset U_{\mathrm{Sp}_{2m}}(F)$ given by
$$U_{1,\kappa}:=\{(u_{1,ij})\in U_{\mathrm{GL}_r(F)}: \vert u_{1, ij}\vert\leq q^{(j-i)\kappa}\}, \quad U_{2,\kappa}:=\{(u_{2,ij})\in U_{\mathrm{Sp}_{2m}(F)}: \vert u_{2,ij}\vert\leq q^{(j-i)\kappa}\}.$$
Then $U_{0,\kappa}$ are open compact subgroups of $U_M(F)$. For $u=(u_1,u_2)\in U_{0,\kappa}$, one checks easily that $\Theta_M(u^{-1})\overline{N}_{0,\kappa}u\subset \overline{N}_{0,\kappa}$, i.e,  $\overline{N}_{0,\kappa}$ is invariant under the twisted conjugate action by $U_{0,\kappa}$.

 Let $\varphi(m(X,Y)):=\varphi_\kappa((Y^{-1}X)_{r1}X)\varphi_\kappa(((Y^{-1}X)_{r1})^2 Y)$, so $\varphi$ is defined on the subset of $M$ which lie in the image of $n(X,Y)$ under the map $n\mapsto m$ via $\dot{w}_0^{-1}n=mn'\bar{n}$, where $(Y^{-1}X)_{r1}$ is the $(r,1)$-th entry of $Y^{-1}X$. We will show that $\varphi$ is invariant under the $\Theta_M$-twisted conjugate action $U_{0,\kappa}$. Suppose $u=(u_1, u_2)\in U_{0,\kappa}$, then
 $$\varphi(\Theta_M(u^{-1})m(X,Y)u)=\varphi(m(u_1^{-1}Xu_2, u_1^{-1}Y\theta_r(u_1)))$$$$=\varphi_{\kappa}((\theta_r(u_1^{-1})Y^{-1}Xu_2)_{r1})u_1^{-1}Xu_2)\varphi_{\kappa}((\theta_r(u_1^{-1})Y^{-1}Xu_2)_{r1})^2u_1^{-1}Y\theta_r(u_1))$$
 $$=\varphi_{\kappa}((Y^{-1}X)_{r1})u_1^{-1}Xu_2)\varphi_{\kappa}((Y^{-1}X)_{r1})^2u_1^{-1}Y\theta_r(u_1)).$$ The last equality holds since $(Y^{-1}X)_{r1}$ is clearly invariant under the left or right translation by an upper triangular unipotent matrix. On the other hand, since $u_1(Y^{-1}X)_{r1}u_1^{-1}=(Y^{-1}X)_{r1}$,
 \begin{align*}
\varphi_\kappa((Y^{-1}X)_{r1}u_1^{-1}Xu_2)=1&\Leftrightarrow (Y^{-1}X)_{r1}u_1^{-1}Xu_2\in X(\kappa)&\Leftrightarrow X\in u_1(Y^{-1}X)_{r1}u_1^{-1} (u_1X(\kappa)u_2^{-1})=(Y^{-1}X)_{r1}X(\kappa)\\
&\Leftrightarrow \varphi_{\kappa}((Y^{-1}X)_{r1}X)=1,&
 \end{align*} where $X(\kappa)=\{X=(x_{ij})\in\mathrm{Mat}_{r\times 2m}: \vert x_{ij}\vert \leq q^{((r-i)+(2m-j)+1)\kappa}\}$, which can be checked by direct calculation that it
   is invariant under the left-action by $U_{1,\kappa}$ and right-action by $U_{2,\kappa}$. Similarly we have
 $\varphi_{\kappa}(((Y^{-1}X)_{r1})^2 u_1^{-1}Y\theta_r(u_1))=\varphi_{\kappa}(((Y^{-1}X)_{r1})^2Y)$. It follows that
 $$\varphi(\Theta_M(u^{-1})m(X,Y)u)=\varphi(m(X,Y))$$ for all $u\in U_{0,\kappa}$.

Then for $(X',Z')\in R_{X'}\times R_{Z'}$, $z=\alpha^\vee(\varpi^{d+g}u_{\alpha_r}(\dot{w}_0\bar{n}(X',Z')\dot{w}_0^{-1}))=\alpha^\vee(\varpi^{d+g}\frac{{y'}^*_{rr}}{\det(Y'))})$, where ${y'}_{rr}^*$ is the adjoint matrix of the $(r,r)$-th entry of $Y'$, we have
$$\varphi_{\overline{N}_{0,\kappa}}(zu^{-1}\bar{n}(X',Z')uz^{-1})=\varphi_\kappa(u_1^{-1}\frac{{y'}^*_{rr}}{\det(Y')}X' u_2)\varphi_\kappa(u_1^{-1}(\frac{{y'}^*_{rr}}{\det Y')})^2Y'\theta_r(u_1))$$
$$=\varphi_{\kappa}(u_1^{-1}({X'}^{-1}Y')_{r1}X'u_2)\varphi_\kappa(u_1^{-1}(({Y'}^{-1}X')_{r1})^2 Y'\theta_r(u_1))=\varphi(\Theta_M(u^{-1})m(X', Z')u).$$ 

Since the uniqueness of the Bruhat decomposition $\dot{w}_0^{-1}n=mn'\bar{n}$ implies that when $n\mapsto u^{-1}nu$, then $m\mapsto \Theta_M(u)^{-1}mu$. By our calculation, for $u=(u_1,u_2)\in U_M$, $n=n(X,Y)\mapsto m=m(X,Y)$, we have $u^{-1}n(X,Y)u=n(u_1^{-1}Xu_2, u_1^{-1}Y\theta(u_1))$, $\Theta_M(u)^{-1}m(X,Y)u=m(u_1^{-1}Xu_2, u_1^{-1}Y\theta(u_2))$. So both the actions are given by $X\mapsto u_1^{-1}Xu_2$, $Y\mapsto u_1^{-1}Y\theta(u_1)$ Therefore, $$U_{M, m(X,Y)}^{\Theta_M}=U_{M,n(X,Y)}=\{u=(u_1,u_2)\in U_M: u_1^{-1}Xu_2=X, u_1^{-1}Y\theta(u_1)=Y\}$$ This shows that the centralizer of $n(X,Y)$ in $U_M$ agrees with the $\Theta_M$-twisted centralizer of $m(X,Y)$ in $U_M$. Compare the definitions of $j_{\overline{N}_{0,\kappa},\pi,w}(n(X',Z'))$ and $B^M_\varphi(m(X',Z'), f)$, we obtain the desired equality.
\end{proof}
Finally, by Proposition \ref{prop5.2} and Proposition \ref{prop6.8}, we can restate Proposition \ref{prop5.2} as:
\begin{prop}\label{prop6.10}
     Let $\sigma$ and $\tau$ be $\psi$-generic supercuspidal representations of $\mathrm{GL}_r(F)$ and $\mathrm{Sp}_{2m}(F)$ respectively, such that $\omega_\pi$ is ramified where $\pi=\sigma\boxtimes \tau$, then there exists a $\kappa_0$, such that for all $\kappa\ge\kappa_0$ and all characters $\chi$ of $F^\times$ such that $\omega_{\pi\otimes\chi}=\omega_\pi\chi^r$ is ramified, we have
$$C_\psi(s,\pi\otimes \chi)^{-1}=\gamma(2rs, \omega_\pi^2\chi^{2r}, \psi^{-1})\int_{R_{X'}\times R_{Z'}}B^M_\varphi(m(X',Z'), f)(\omega_{\pi}^{-2}\chi^{-2r})(\frac{P(X',Z')}{\det(Z'J_r+\frac{X' \theta_{r,m}(X')}{2})})$$$$\cdot \vert \frac{P(X',Z')}{\det(Z'J_r+\frac{X'\theta_{r,m}(X')}{2})}\vert^{-rs} \vert \det(m_1(X',Z'))\vert^{s+\frac{r+2m+1}{2}} d\mu_{X'}\wedge d\mu_{Z'}$$
where $m(X',Z')$ is the image of $n(X',Z')\mapsto m(X', Z')$ via the Bruhat decomposition $\dot{w}_0^{-1}\dot{n}=mn'\bar{n}$,  which holds off a subset of measure zero on $N(F)$, $m=\mathrm{diag}\{m_1, m_2,\theta_r(m_1)\}$ with $m_1\in\mathrm{GL}_r(F)$, $m_2\in\mathrm{Sp}_{2m}(F)$, and $\gamma(2r s, \omega_\pi^2\chi^{2r}, \psi^{-1})$ is the Abelian $\gamma$-factor depending only on $\omega_\pi$ and $\chi$.
\end{prop}

One important observation in our case is that our orbit space $R'=R_{X'}\times R_{Z'}$ is no longer isomorphic to a torus, and attempts to reparameterize $m(X,Z)$ by the maximal torus together with non-torus parts fail since the map $(X,Z)\mapsto (m_1, m_2)$ is not surjective in general. So in order to apply the asymptotic analysis of partial Bessel integrals, we will need to study the toric action on the orbit space and separate its toric part out from the integral. First we would like to understand the action of the maximal torus $A\simeq A_1\times A_2$ on $R=R_X\times R_Z=U_M\backslash N$,  where $A_1$ is the maximal torus of $\mathrm{GL}_r$ and $A_2$ is the maximal torus of $\mathrm{Sp}_{2m}$. Recall that through the Bruhat decomposition $\dot{w}_0^{-1}n=mn'\bar{n}$, for $s\in A$, if $n\mapsto sns^{-1}$, then $m\mapsto \Theta_M(s)ms^{-1}$. In our case, for $m=(m_1,m_2)$ with $m_1\in \mathrm{GL}_r$ and $m_2\in\mathrm{Sp}_{2m}$, the action is given by $m_1\mapsto s'm_1\theta_r({s'}^{-1})$, and $m_2\mapsto s''m_2s''^{-1}$, which is equivalent to $X\mapsto s'X{s''}^{-1}$, $Y\mapsto s''Y\theta_r(s'')^{-1}$, or to $X\mapsto s'X{s''}^{-1}$, $Z\mapsto s'Zs'$, where $ s=(s',s'')$ with $s'\in A_1$ and $s''\in A_2$. Write $s'=\mathrm{diag}(s_1,s_2,\cdots, s_r)\in A_1, s'=\mathrm{diag}(s_{r+1},s_{r+2},\cdots, s_{r+m},{s_{r+m}}^{-1},\cdots, {s_{r+1}}^{-1})\in A_2$. Define
$$T_{X,Z}=\begin{cases}
    \{(x_{r,1},x_{r-1,2}\cdots, x_{r-m+1,m}, z_{1,1},\cdots, z_{r,r}): (X,Z)\in R_{X}\times R_Z\}, & \mathrm{ \ if \ } r\ge m,\\
    \{(x_{r,1},x_{r-1,2},\cdots, x_{1,r},z_{1,1},\cdots, z_{r,r}): (X,Z)\in R_X\times R_Z\}, &\mathrm{\ if\ } r<m.\\
\end{cases}$$
\begin{prop} \label{prop6.11}Let $\Phi: R_X\times R_Z\rightarrow M\rightarrow A=A_1\times A_2$ be the map given by $(X,Z)\mapsto m=m(X,Z)=(m_1, m_2)\mapsto (t_1,t_2)$, where $m(X,Z)$ is the image of $n(X,Z)$ under the map $n\mapsto m$ via $\dot{w}_0^{-1}n=mn'\bar{n}$, and $m_1=u_1\dot{w}_1 t_1 u_2$, $m_2=v_1\dot{w}_2t_2 v_2$ are the Bruhat decomposition of $m_1$ and $m_2$ respectively. Then $\Phi$ is compatible with the toric action, i.e., $\Phi(s'X{s''}^{-1}, s'Zs')=(m_1(s'X{s''}^{-1}, s'Zs'), m_2(s'X{s''}^{-1}, s'Zs'))$ for $s=(s',s'')\in A$. The toric action on $R_X\times R_Z$ is completely determined by its restriction on $T_{X,Z}$, up to the square of each $s_i(1\le i\le r+\mathrm{min}\{r, m\})$. Moreover, the restriction $\Phi\vert_{T_{X,Z}}$ is a finite $\acute{e}tale$ map onto its image in $A$. 
\end{prop}

\begin{proof} The compatibility of $\Phi$ with the toric action follows directly from $\Theta_M(s)m(X,Z)s^{-1}=m(s'X{s''}^{-1}, s'Zs') $. Recall that we also write $m=m(X,Z)=m(X,Y)$ where $Y=ZJ_r+\frac{X\theta_{r,m}(X)}{2}$. We claim that if we inductively write $Y=Y_r=\begin{bmatrix}
    \alpha & y_{1,r}\\
    Y_{r-1} & {^t\beta}\\
\end{bmatrix}$ and $Y=Y^{(r)}=\begin{bmatrix}
    {^t\alpha'} & Y^{(r-1)}\\
    y_{r,1} & \beta'\\
\end{bmatrix}$, and assume that $\det(Y_i)$'s and $\det(Y^{(j)})$'s are all non-zero, then the image $m(X,Y)$ of $n(X,Y)$ lies in the big cell of $M$. If we decompose $m_1=u_1\dot{w}_1t_1 u_2$, $m_2=v_1\dot{w}_2t_2 v_2$ with $\dot{w}_1=J_r$, $\dot{w}_2=J'_{2m}$, the long Weyl group elements of $\mathrm{GL}_r$ and $\mathrm{Sp}_{2m}$ respectively, we will show that
    \[
t_1=\mathrm{diag}((-1)^{r-1}\frac{\det(Y_{r-1})}{\det(Y_r)}, (-1)^{r-2}\frac{\det(Y_{r-2})}{\det(Y_{r-1})}, \cdots, \frac{1}{\det(Y_1)} )
\] and $t_2=\mathrm{diag}(\bar{t}_2,{\bar{t}_2}^{-1})$ with
$$\bar{t}_2=\begin{cases}
    (-\frac{\det(Y^{(r-1)})}{\det(Y^{(r)})}x_{r,1}^2, (-1)^2\frac{\det(Y^{(r-2)})}{\det(Y^{(r-1)})}x_{r-1,2}^2, \cdots, (-1)^{r-1}\frac{1}{\det(Y^{(1)})}x_{1,r}^2,1,\cdots, 1 )  & \mathrm{\  if \ } r<m\\
    (-\frac{\det(Y^{(r-1)})}{\det(Y^{(r)})}x_{r,1}^2, (-1)^2\frac{\det(Y^{(r-2)})}{\det(Y^{(r-1)})}x_{r-1,2}^2, \cdots, (-1)^{m}\frac{\det(Y^{(r-m)})}{\det(Y^{(r-m+1)})}x_{r-m+1,m}^2 ) & \mathrm{\ if \ }r\ge m\\
\end{cases}
$$
Let us first prove this for $t_1$ by induction on $r$. If $r=1$ then there is nothing to show. Assume $r>1$ and write $Y_r=\begin{bmatrix}
    \alpha & y_{1,r}\\
    Y_{r-1} & {^t\beta}\\
\end{bmatrix}$. Take the representative of the long Weyl group element in $\mathrm{GL}_r$ as $J_r$. Then $m_1=\theta_{r}(Y)=u_1 J_r t_1 u_2\Leftrightarrow \tilde{u}_1Y\tilde{u}_2=t_1^{-1}J_r$, where $\tilde{u}_i=\theta_r(u_i)^{-1}(i=1,2)$.
Write $\tilde{u}_1=\begin{bmatrix}
    1 & \delta_1\\
    & u'_1\\
\end{bmatrix}$, and $\tilde{u}_2=\begin{bmatrix}
    u'_2 & {^t\delta_2}\\
     & 1\\
\end{bmatrix}$, with $u_i'\in U_{\mathrm{GL}_{r-1}}$, $\delta_i\in \mathbb{A}^{r-1}$. Then
$$\tilde{u}_1 Y\tilde{u}_2=\begin{bmatrix}
    (\alpha+\delta_1Y_{r-1})u_2' & (\alpha+\delta_1Y_{r-1}){^t\delta_2}+(y_{1,r}+\delta_1{^t\beta})\\
    u_1' Y_{r-1} u_2' & u_1'Y_{r-1}{^t\delta}_2+u_1'\beta\\
\end{bmatrix}.$$
Assume $\det(Y_{r-1})\neq 0$, $y_{1,r}\neq 0$, choose $\delta_1=-\alpha Y^{-1}_{r-1}$, ${^t\delta_2}=-Y_{r-1}^{-1}\beta$, write $t_1=\mathrm{diag}\{a_1,\cdots, a_r\}$ then $\tilde{u}_1Y\tilde{u}_2=t_1^{-1}J_r\Leftrightarrow$
$\begin{bmatrix}
    0 & y_{1,r}-\alpha Y_{r-1}^{-1}{^t\beta}\\
    u'_1Y_{r-1}u'_2 & 0
\end{bmatrix}=\begin{bmatrix}
    & & a_1^{-1}\\
    & \iddots & \\
    (-1)^{r-1}a_r & & \\
\end{bmatrix}$. By induction hypothesis, it suffices to show that $y_{1,r}-\alpha Y_{r-1}^{-1}{^t\beta}=a_1^{-1}$. Note that
$\det(Y_r)=\det\begin{bmatrix}
    \alpha & y_{1,r}\\
    Y_{r-1} & {^t\beta}\\
    \end{bmatrix}=(-1)^{r-1}\det\begin{bmatrix}
    Y_{r-1} & {^t\beta}\\
     \alpha & y_{1,r}\\
    \end{bmatrix}$. On the other hand, $\begin{bmatrix}
        Y^{-1}_{r-1} & \\
        & 1\\
    \end{bmatrix}\begin{bmatrix}
    Y_{r-1} & {^t\beta}\\
     \alpha & y_{1,r}\\
    \end{bmatrix}=\begin{bmatrix}
        I_{r-1} & Y_{r-1}^{-1}{^t\beta}\\
        \alpha & y_{1,r}\\
    \end{bmatrix}$. Taking determinant on both sides, we obtain that $$(-1)^{r-1}\det(Y_r)\det(Y_{r-1})^{-1}=\det\begin{bmatrix}
        I_{r-1} & Y_{r-1}^{-1}{^t\beta}\\
        \alpha & y_{1,r}\\
    \end{bmatrix}=\det\begin{bmatrix}
        I_{r-1} & Y_{r-1}^{-1}{^t\beta}\\
        0 & -\alpha Y^{-1}_{r-1}{^t\beta}+y_{1,r}\\
    \end{bmatrix}= y_{1,r}-\alpha Y^{-1}_{r-1}{^t\beta},$$
    i.e., $a_1=(-1)^{r-1}\frac{\det(Y_{r-1})}{\det(Y_r)}.$

    Next, we work on the formula for $t_2$. Let us first assume that $r\ge m$ and we prove the formula by induction on $m$.
When $m=1$, the result follows from direct calculation.
Assume $m>1$ and suppose
$m_2=v_1J'_{2m}t_2 v_2$, then
$m_2^{-1}=I_{2m}-J'_{2m}{^tX}{^tY^{-1}}J_rX=-v_2^{-1}t_2^{-1}J'_{2m} v_1^{-1}$, which is equivalent to $v_2(I_{2m}-J'_{2m}{^tX}{^tY^{-1}}J_rX)v_1=-t_2^{-1}J'_{2m}$.
Write $v_i=\begin{bmatrix}
    1 & \gamma_i & b_i\\
    & v_i' & \gamma_i^*\\
    & & 1
\end{bmatrix}, (i=1,2)$, $X=\begin{bmatrix}
    0 & X_1 & {^t\gamma}\\
    x_{r,1} & 0 & 0\\
\end{bmatrix}$, $Y=Y^{(r)}=\begin{bmatrix}
    {^t\alpha'} & Y^{(r-1)}\\
    y_{r,1} & \beta'\\
\end{bmatrix}$. Assume that each $Y^{(j)},(1\le j\le r)$ is invertible, $y_{r,1}\neq 0$, and set ${^t(Y^{(r)})^{-1}}=\begin{bmatrix}
    {^t\delta'_1} & H\\
    x & \delta_2'\\
\end{bmatrix}$. Then one computes that
$$I_{2m}-J_{2m}'{^tX}{^tY^{-1}}J_rX=I_{2m}-\begin{bmatrix}
     & & 1 \\
     & -J'_{2m-2} & \\
     -1 & & \\
\end{bmatrix}\begin{bmatrix}
    0 & x_{r,1}\\
    {^tX_1} & 0\\
    \gamma & 0\\
\end{bmatrix}\begin{bmatrix}
    {^t\delta'_1} & H \\
    x & \delta'_2\\
\end{bmatrix}\begin{bmatrix}
    & 1\\
    -J_{r-1} & \\
\end{bmatrix}\begin{bmatrix}
    0 & X_1 & {^t\gamma}\\
    x_{r,1} & 0 & 0\\
\end{bmatrix}$$
$$=\begin{bmatrix}
    1-\gamma{^t\delta_1'}x_{r,1} & \gamma HJ_{r-1}X_1 & \gamma H J_{r-1}{^t\gamma}\\
    J'_{2m-2} {^tX_1}{^t\delta_1'}x_{r,1} & I_{2m-2}-J'_{2m-2}{^tX_1}HJ_{r-1}X_1 & -J'_{2m-2}{^tX_1}HJ_{r-1} {^t\gamma}\\
    x^2_{r,1}x & -x_{r,1}\delta'_2J_{r-1}X_1 & 1-x_{r,1}\delta'_2 J_{r-1}{^t\gamma}\\
\end{bmatrix}$$

By definition of ${^t(Y^{(r)})^{-1}}$, we have
$$\alpha'{^t\delta_1'}+y_{r,1}x=1, \alpha' H+y_{r,1}\delta'_2=0,{^tY^{(r-1)}}{^t\delta_1'}+{^t\beta'x}=0, {^tY^{(r-1)}}H+{^t\beta'\delta_2'}=I_{r-1}.$$
Multiply by $v_2$ on the left and $v_1$ on the right, and compare with $-t_2^{-1}J'_{2m}=\begin{bmatrix}
    & & a_{r+1}^{-1}\\
    & -t_2'^{-1}J'_{2m-2} & \\
    -a_{r+1} & & \\
\end{bmatrix}$, where we write $t_2=\mathrm{diag}\{a_{r+1}, \cdots, a_n, a_n^{-1},\cdots, a_{r+1}^{-1}\}=\mathrm{diag}\{a_{r+1},t_2', a_{r+1}^{-1}\}$. The middle entry of $v_2(I_{2m}-J'_{2m}{^tX}{^tY^{-1}}J_rX)v_1$ gives
$$v_2'(I_{2m-2}-J'_{2m-2}{^tX_1}{^t(Y^{(r-1)})}^{-1}J_{r-1}X_1)v_1'$$
$$+v_2'J'_{2m-2}{^tX_1}{^t(Y^{(r-1)})^{-1}}{^t\beta'\delta_2'}J_{r-1}X_1v_1'+\gamma_2^*(-x_{r,1}\delta_2'J_{r-1}X_1)v_1'=-t_2'^{-1}J'_{2m-2}.$$ Therefore we have to show that $$v_2'J'_{2m-2}{^tX_1}{^t(Y^{(r-1)})^{-1}}{^t\beta'\delta_2'}J_{r-1}X_1v_1'+\gamma_2^*(-x_{r,1}\delta_2'J_{r-1}X_1)v_1'=0\ \ \ \ \ \ \ \ \ (*)$$ in order to continue with our induction procedure.
Since ${^tY^{(r-1)}}{^t\delta_1'}+{^t\beta'x}=0,$ we have $ {^t\delta_1'}+{^t(Y^{(r-1)})^{-1}}{^t\beta'}x=0$, thus ${^t\delta_1'}\delta_2'+{^t(Y^{(r-1)})^{-1}}{^t\beta'}\delta_2'x=0$. We also note that $\gamma_2^*=v_2' J'_{2m-2}{^t\gamma_2}$ by the structure of $\mathrm{Sp}_{2m-2}$.
So the left hand side of (*) is equal to $$v_2'J'_{2m-2}(-{^tX}_1{^t\delta_1'}\delta_2'-{^t\gamma_2}x_{r,1}\delta_2')J_{r-1}X_1v_1'.$$
On the other hand, one computes that the $(2,1)$-th entry of $v_2(I_{2m}-J'_{2m}{^tX}{^tY^{-1}}J_rX)v_1$ gives
$v_2'J'_{2m-2}({^tX}_1{^t\delta_1'}x_{r,1}+{^t\gamma_2}x^2_{r,1}x)=0$,
hence ${^tX}_1{^t\delta_1'}x_{r,1}+{^t\gamma_2}x^2_{r,1}x=0$. As a result,
$${^tX}_1{^t\delta_1'}\delta_2'x_{r,1}+{^t\gamma_2}x^2_{r,1}\delta_2'x=0.$$ Since $x_{r,1}\neq 0$, this shows that the left hand side of (*) is equal to 0. 

By induction hypothesis, it suffices to show that $a_{r+1}=-x_{r,1}^2\frac{\det(Y^{(r-1)})}{\det(Y^{(r)})}$. Note that the $(3,1)$-th entry of $v_2(I_{2m}-J'_{2m}{^tX}{^tY^{-1}}J_rX)v_1$ implies that
$x_{r,1}^2 x=-a_{r+1}$. One also checks by direct calculation that other identities in the comparison of entries in $v_2(I_{2m}-J'_{2m}{^tX}{^tY^{-1}}J_rX)v_1$ and $-t_2^{-1}J'_{2m}$ are automatically satisfied. On the other hand, $\alpha'{^t\delta'_1}+y_{r,1}x=1$ together with ${^t\delta_1}=-{^t(Y^{(r-1)})^{-1}}{^t\beta'}x$ imply that
$x(y_{r,1}-\alpha'{^t(Y^{(r-1)})^{-1}}{^t\beta'})=1$.
We have $$\frac{\det(Y^{(r)})}{\det(Y^{(r-1)})}=\det(Y^{(r)}\begin{bmatrix}
    (Y^{(r-1)})^{-1} & \\
    & 1
\end{bmatrix})=\det(\begin{bmatrix}
    Y^{(r-1)} & {^t\alpha'}\\
    \beta' & y_{r,1}\\
\end{bmatrix}\begin{bmatrix}
    (Y^{(r-1)})^{-1} & \\
    & 1
\end{bmatrix})$$
$$=\det\begin{bmatrix}
    I_{r-1} & {^t\alpha_1}\\
    \beta'(Y^{(r-1)})^{-1} & y_{r,1}\\
\end{bmatrix}=\det\begin{bmatrix}
    I_{r-1} & {^t\alpha'}\\
    0 & -\beta'(Y^{(r-1)})^{-1}{^t\alpha'}+y_{r,1}\\
\end{bmatrix}=y_{r,1}-\alpha'{^t(Y^{(r-1)})^{-1}}{^t\beta'}.$$
Hence by assuming that $y_{r,1}-\alpha'{^t(Y^{(r-1)})^{-1}}{^t\beta'}\neq 0$, we obtain that
$a_{r+1}=-x_{r,1}^2\frac{\det(Y^{(r-1)})}{\det(Y^{(r)})}$ as desired. 

When $r<m$, recall that the toric action has a stabilizer isomorphic to the maximal torus of $\mathrm{Sp}_{2(m-r)}$, given by the entries $(s_{2r+1}, \cdots, s_{r+m})$. So the above induction process terminates at the $r$-th step. As a result, the formula directly follows from the same procedure. 

Given the explicit formula for $t_1$ and $t_2$ in the Bruhat decomposition, when restricted to $T_{X,Z}$,  
the map $\Phi$ is completely determined by $Y\mapsto \det(Y_i)(1\le i\le r )$, $Y\mapsto \det(Y^{(i)})(1\le i\le r)$, $X\mapsto x_{r-j+1,j}^2(1\le j\le \mathrm{min}(r,m))$, which are clearly finite $\acute{e}tale$ maps since we assumed that $\det(Y^i)$'s, $\det(Y^{(i)})$'s, and $x_{r-j+1,j}$'s are all non-zero. 
The restriction of the toric action $X\mapsto s'X{s''}^{-1}$, $Z\mapsto s' Zs'$, or equivalently, $X\mapsto s'X{s''}^{-1}$, $Y\mapsto s'Y\theta_r(s')^{-1}$, on $T_{X,Z}$ is given by
$x_{r-i+1,i}\mapsto s_{r-i+1}s_{r+i}^{-1}x_{r-i+1}$, $z_{j,j}\mapsto s_j^2z_{j,j}$ for $1\le i\le m$ if $r\ge m$, $1\le i\le r$ if $r<m$, and $1\le j\le r$. Moreover, when $r<m$, the action has a non-trivial stabilizer isomorphic to the maximal torus of $\mathrm{Sp}_{2(m-r)}$, given by the variables $(s_{2r+1}, \cdots, s_{r+m})$. A simple calculation shows that it takes $\det(Y_i)$ to $\prod_{k=1}^rs^2_{r-k+1}\det(Y_i)$, $\det(Y^{i})$ to $\prod_{k=1}^i s_i^2\det(Y^{(i)})$, and $x_{r-k+1,k}$ to $\frac{s_{r-k+1}}{s_{r+k}}x_{r-k+1}$. Together with the compatibility of the toric action, it implies that $$t_1\mapsto {s'}^{-2}t_1, t'_2\mapsto \begin{cases}
    \mathrm{diag}\{s_{r+1}^{-2},\cdots, s_{2r}^{-2}, 1,\cdots, 1,\}t_2',& \mathrm{\ if \ }r<m,\\
    {s''}^{-2}t_2' ,& \mathrm{\ if \ } r\ge m.\\
\end{cases}$$ Consequently, the toric action on $R_X\times R_Z$ is completely determined by its action on $T_{X,Z}$ up to the square of each $s_i(1\leq i\leq r+\mathrm{min}\{r,m\})$. It follows that the covering group is isomorphic to finitely many copies of $\mathbb{Z}/2\mathbb{Z}$.
\end{proof}  
Finally, in order to apply our uniform smooth results for partial Bessel integrals, we need to verify the non-singular condition in Proposition \ref{prop6.8}:

\begin{prop}\label{prop6.12}
In the Bruhat decompositions    
$m_1=u_1\dot{w}_1t_1u_2$, and $m_2=v_1\dot{w}_2t_2 v_2$, the entries in the 1-parameter subgroups for $u_1, u_2, v_1, v_2$ are rational functions of the image of $T_{X,Z}$ under $\Phi$ without singularities.
\end{prop}

\begin{proof} Based on the calculation in Proposition \ref{prop6.11}, it suffices to show that the denominator of entries in $u_1, u_2, v_1, v_2$ are polynomial functions of $\det(Y^{(i)})$, $\det(Y_i)$ and $x_{r-i+1, i}$ without singularities. In fact, we show that they are monomials of these factors. We start with $m_1=u_1w_1 t_1 u_2$. Denote
$u_1=(u^1_{ij}), u_2=(u^2_{ij})
$, $Y=(y_{ij})$. Let $$Y_{i,j}=Y[r-j+1, r-j+2, \cdots,\widehat{r-i+1}\cdots, r;1,\cdots,  j-1],$$ 
the matrix of size $(j-1)\times (j-1)$ with elements from the indicated columns and rows of $Y$. We will show by induction on $r$ that
\[
u_{ij}^{1}=\frac{\det(Y_{i,j})}{\det(Y_{j-1})}.
\] When $r=1$ there is nothing to show. Suppose $r>1$, write $u_1=\begin{bmatrix}
    u_1 & {^t\delta_1}\\
    & 1\\
\end{bmatrix}$, $u_2=\begin{bmatrix}
    1 & \delta_2\\
    & u_2'\\
\end{bmatrix}$, $t=\mathrm{diag}\{a_1,\cdots, a_r\}=\mathrm{diag}\{a_1, t_1'\}$, and $Y=\begin{bmatrix}
    \alpha & y_{1,r}\\
    Y_{r-1} & {^t\beta}
\end{bmatrix}$. Note that $m_1=\theta_r(Y)=u_1 J_r t_1 u_2\Leftrightarrow Y=Y_{r}=J_r{^tu_1^{-1}}J_r^{-1}t_1^{-1}{^tu_2^{-1}}J_r$.
We compute that
$$\begin{bmatrix}
    \alpha & y_{1,r}\\
    Y_{r-1} & {^t\beta}
\end{bmatrix}=\begin{bmatrix}
    & 1\\
    -J_{r-1} & \\
\end{bmatrix}\begin{bmatrix}
    {^tu_1'^{-1}} & \\
    -\delta_1{^tu'}_1^{-1} & 1\\
\end{bmatrix}\begin{bmatrix}
    & -J_{r-1}^{-1}\\
    1 & \\
\end{bmatrix}\begin{bmatrix}
    a_1^{-1} & \\
    & t_1'^{-1}\\
\end{bmatrix}\begin{bmatrix}
    1  & \\
    -{^tu'}_2^{-1}{^t\delta_2} & {^tu'}_2^{-1}
\end{bmatrix}\begin{bmatrix}
    & 1 \\
    -J_{r-1} & \\
\end{bmatrix}$$
$$=\begin{bmatrix}
    -\delta_1{^tu'}_1^{-1}J_{r-1}^{-1}t_1'^{-1}{^tu'}_2^{-1} J_{r-1} & a_1^{-1}-\delta_1{^tu'}_1^{-1}J_{r-1}^{-1}t_1'^{-1}{^tu'}_2^{-1} {^t\delta_2}\\
    -J_{r-1}{^tu'}_1^{-1}J_{r-1}^{-1}t_1'^{-1}{^tu'}_2^{-1} J_{r-1} & -J_{r-1}{^tu'}_1^{-1}J_{r-1}^{-1}t_1'^{-1}{^tu'}_2^{-1} {^t\delta_2}\\
\end{bmatrix}.$$
Hence
$Y_{r-1}=-J_{r-1}{^tu'}_1^{-1}J_{r-1}^{-1}t_1'^{-1}{^tu'}_2^{-1} J_{r-1}$. Compare this with the expression for $Y_r$, we replace $Y=Y_r$ by $-Y_{r-1}$ and continue by induction. Suppose we show that $u^1_{i,r}=\frac{\det(Y_{i,r})}{\det(Y_{r-1})}, (1\le i\le r-1)$, then the induction hypothesis would imply that $$u^1_{i,r-1}={u'}^1_{i,r-1}=\frac{\det((-Y_{r-1})_{i,r-1})}{\det((-Y_{r-1})_{r-2})}=\frac{(-1)^{r-2}\det(Y_{i,r-1})}{(-1)^{r-2}\det(Y_{r-2})}=\frac{\det(Y_{i,r-1})}{\det(Y_{r-2})}, (1\le i\le r-2),$$ since $$(Y_{r-1})_{i,r-1}=Y_{r-1}[1,\cdots, \widehat{(r-1)-i+1},\cdots, r-1; 1,\cdots, r-2]$$$$=Y_r[2,\cdots,\widehat{r-i+1},\cdots, r; 1,\cdots, r-2]=Y_{i,r-1}.$$ Therefore it suffices to show that $u_{i,r}^1=\frac{\det(Y_{i,r})}{\det(Y_{r-1})}$ for $1\le i\le r-1$. From the above calculation we also have $\alpha=-\delta_1J_{r-1}^{-1}Y_{r-1}$, ${^t\beta}=-Y_{r-1}J_{r-1}^{-1}{^t\delta_2}$. So
$$\delta_1=-\alpha Y_{r-1}^{-1}J_{r-1}=(u^1_{1,r},u^1_{2,r},\cdots, u^1_{r-1,r}),$$$$ \delta_2=-\beta{^tY_{r-1}^{-1}}J_{r-1}^{-1}=(u^2_{1,2},u^2_{1,3},\cdots, u^2_{1,r}).$$
We will show the formula for $u_1$ first.
By the construction of $Y_{i,j}$, we expand the determinant of $Y$ along the last column to get  $$\det(Y_r)=(-1)^{r+1}y_{1,r}\det(Y_{r-1})+(-1)^{r+2}y_{2,r}\det(Y_{r-1,r})+\cdots +(-1)^{2r}y_{r,r}\det(Y_{1,r}).$$
Since $\beta=(y_{2,r},\cdots, y_{r-1,r})$, if we set $$\delta=((-1)^{r+2}\frac{\det(Y_{r-1,r})}{\det (Y_{r-1})},(-1)^{r+3}\frac{\det (Y_{r-2,r})}{\det (Y_{r-2})},\cdots, (-1)^{r+l+1}\frac{\det (Y_{r-l, r})}{\det (Y_{r-1})},\cdots, (-1)^{2r}\frac{\det (Y_{1,r})}{\det (Y_{r-1})}),$$ then the above formula is equivalent to 
$\det(Y_r)=(-1)^{r+1}y_{1,r}\det(Y_{r-1})+\det(Y_{r-1})\delta{^t\beta}$.
On the other hand, we have
$\begin{bmatrix}
        Y^{-1}_{r-1} & \\
        & 1\\
    \end{bmatrix}\begin{bmatrix}
    Y_{r-1} & {^t\beta}\\
     \alpha & y_{1,r}\\
    \end{bmatrix}=\begin{bmatrix}
        I_{r-1} & Y_{r-1}^{-1}{^t\beta}\\
        \alpha & y_{1,r}\\
    \end{bmatrix}$. Taking determinant on both sides, we obtain that $$(-1)^{r-1}\det(Y_r)\det(Y_{r-1})^{-1}=\det\begin{bmatrix}
        I_{r-1} & Y_{r-1}^{-1}{^t\beta}\\
        \alpha & y_{1,r}\\
    \end{bmatrix}=\det\begin{bmatrix}
        I_{r-1} & Y_{r-1}^{-1}{^t\beta}\\
        0 & -\alpha Y^{-1}_{r-1}{^t\beta}+y_{1,r}\\
    \end{bmatrix}= y_{1,r}-\alpha Y^{-1}_{r-1}{^t\beta}.$$
Hence
$\det(Y_r)=(-1)^{r-1}y_{1,r}\det(Y_{r-1})+(-1)^r\det(Y_{r-1})\alpha Y_{r-1}^{-1}{^t\beta}$.
Consequently, $((-1)^r\alpha Y_{r-1}^{-1}-\delta){^t\beta}=0$. Note that our formula works for any $\beta$, therefore  $\delta=(-1)^r\alpha Y_{r-1}^{-1}$. As $\delta_1=-\alpha Y_{r-1}^{-1}J_{r-1}$, we have
$\delta=(-1)^{r-1}\delta_1 J_{r-1}^{-1}=\delta_1 J_{r-1}$.
It follows that
$$(-1)^{r+l+1}\frac{\det(Y_{r-l,r})}{\det(Y_{r-1})}=(-1)^{r-l-1}u^1_{r-l,r}$$ for each $1\le l\le r-1$. i,e., $u_{i, r}^1=\frac{\det(Y_{i,r})}{\det(Y_{r-1})}$ for $1\le i\le r-1$. This completes the proof for the formula for $u_1$.

Now we turn to the formula for $u_2$. Let
\[
Y'_{ij}=Y[(i+1),\cdots, r; 1, \cdots, \widehat{r-j+1}, \cdots, r-i+1]
\]
which is of size $(r-i)\times (r-i)$.
We will show that
\[
u_{ij}^{2}=\frac{\det(Y'_{ij})}{\det(Y_{r-i})}.
\]
Similar to the proof for $u_1$, since $Y'_{1,1}=Y_{r-1}$, if we expand $\det(Y)$ along the first row, we get $$\det(Y)=\det(Y_r)=\sum_{j=1}^r(-1)^{j+1}y_{1,j}\det(Y'_{1,r-j+1})=(=1)^{r+1}y_{1,r}\det(Y_{r-1})+\sum_{j=1}^{r-1}y_{1,j}\det(Y'_{1,r-j+1}).$$
Set $\delta'=((-1)^{j+1}\frac{\det(Y'_{1,r-j+1})}{\det(Y_{r-1})})_{j=1}^{r-1}$,
 then $\frac{\det(Y_r)}{\det(Y_{r-1})}=(-1)^{r+1}y_{1,r}+\alpha{^t\delta'}$.

On the other hand, by the above argument we still have $(-1)^{r-1}\det(Y_r)\det(Y_{r-1})^{-1}=y_{1,r}-\alpha Y_{r-1}^{-1}{^t\beta}$.
 So $\alpha((-1)^r{^t\delta'-Y_{r-1}^{-1}{^t\beta}})=0$. Since this works for any $\alpha$, we obtain that
 $(-1)^r{^t\delta'}=Y_{r-1}^{-1}\beta$. Note that the formula $\delta_2=-\beta{^tY_{r-1}^{-1}}J_{r-1}^{-1}$ is equivalent to $ \delta'=(-1)^{r-1}\delta_2J_{r-1}$. As $\delta_2=(u_{1,2},\cdots, u_{1,r})$, this means that
 $(-1)^{j+1}\frac{\det(Y'_{1,r-j+1})}{\det(Y_{r-1})}=(-1)^ju^2_{1,r-j+1}, \  (1\le j\le r-1)$ and therefore
 $u^2_{1,j}=\frac{\det(Y'_{1,j})}{\det(Y_{r-1})},\ (1\le j\le r-1)$. Then a similar induction argument as the proof of formulas for entries of $u_1$ implies that $u^2_{ij}=\frac{\det(Y'_{i,j})}{\det(Y_{r-1})}$ for all $(i,j)$.

 For entries of $v_1$ and $v_2$, recall that in the proof of Proposition \ref{prop6.11}, we wrote $v_i=\begin{bmatrix}
    1 & \gamma_i & b_i\\
    & v_i' & \gamma_i^*\\
    & & 1
\end{bmatrix}, (i=1,2)$, $X=\begin{bmatrix}
    0 & X_1 & {^t\gamma}\\
    x_{r,1} & 0 & 0\\
\end{bmatrix}$, $Y=Y^{(r)}=\begin{bmatrix}
    {^t\alpha'} & Y^{(r-1)}\\
    y_{r,1} & \beta'\\
\end{bmatrix}$. We assumed that $\det(Y^{(j)})\neq 0,(1\le j\le r)$, $y_{r,1}\neq 0$, and wrote ${^t(Y^{(r)})^{-1}}=\begin{bmatrix}
    {^t\delta'_1} & H\\
    x & \delta_2'\\
\end{bmatrix}$. It follows immediately from the computation of $v_2(I_{2m}-J'_{2m}{^tX}{^tY^{-1}}J_rX)v_1=-t_2^{-1}J'_{2m}$ that
\begin{equation}\label{b2-formula-1}
    1-\gamma {^t\delta_1'}x_{r,1}+\gamma_2J'_{2m-2}{^tX_1}{^t\delta_1'} x_{r,1}+b_2x_{r,1}^2x=0,
\end{equation} 
\begin{equation}\label{b2-formula-2}
    \gamma HJ_{r-1}X_1+\gamma_2(I_{2m-2}-J'_{2m-2}{^tX_1}HJ_{r-1}X_1)-b_2x_{r,1}\delta'_2J_{r-1}X_1=0.
\end{equation} Multiplying (\ref{b2-formula-1}) by $\delta_2' J_{r-1}X_1$, (\ref{b2-formula-2}) by $x_{r,1}x$, and taking their sum, we get
$$\gamma_2(J'_{2m-2}{^tX_1}({^t\delta_1'}\delta_2'-Hx)J_{r-1}X_1x_{r,1}+x_{r,1}xI_{2m-2})+\delta_2'J_{r-1}X_1-\gamma({^t\delta_1}'\delta_2'-Hx)J_{r-1}X_1x_{r,1}=0.$$ Since ${^tY^{(r-1)}}H+{^t\beta}\delta_2'=I_{r-1}$, it follows that ${^t\delta_1'}\delta_2'-Hx=-{^t(Y^{(r-1)})^{-1}}x$. By induction hypothesis, $I_{2m-2}-J'_{2m-2}{^tX_1}{^t{Y^{(r-1)}}^{-1}}J_{r-1}X_1=-t_2'J'_{2m-2}$, where we write $-t_2^{-1}J'_{2m}=\begin{bmatrix}
    & & a_{r+1}^{-1}\\
    & -t_2'^{-1}J'_{2m-2} & \\
    -a_{r+1} & & \\
\end{bmatrix}$. Therefore $$\gamma_2=-x_{r,1}^{-1}x^{-1}(\delta_2'J_{r-1}X_1+\gamma{^t(Y^{(r-1)})^{-1}}J_{r-1}X_1xx_{r,1})v_1't_2'J_{2m-2}v_2'.$$ On the other hand, since $x$ is the $(r,1)$-th entry of ${^tY}^{-1}$, and recall that ${^tY}=\begin{bmatrix}
    {^t\alpha} & {^tY_{r-1}}\\
    y_{1,r} & \beta\\
\end{bmatrix}$. By computing the inverse of a matrix using its adjoint matrix, we have $x=(-1)^{r+1}\frac{\det(Y_{r-1})}{\det(Y^{(r)})}$. Since ${^tY^{-1}}=\begin{bmatrix}
    {^t\delta_1'} & H \\
    x & \delta_2'\\
\end{bmatrix}$, the common denominator of entries in $\delta_2'$ is still $\det(Y^{(r)})$ for the same reason. Similarly $^t(Y^{(r-1)})^{-1}=\frac{({^tY^{(r-1)}})^*}{\det(Y^{(r-1)})}$, where $({^tY^{(r-1)}})^*$ is the adjoint matrix of ${^tY^{(r-1)}}$. In addition, by induction hypothesis and the computation in Proposition \ref{prop6.11}, denominators of entries in $v_1'$, $v_2'$, and $t_2'$ are all monomials of $\det(Y^{(j)}), \det(Y_{j}), (j\le r-1)$, and $x_{r-i+1, i}, (2\le i\le \mathrm{min}\{r,m\})$. Moreover, since $\gamma_2^*=v_2' J'_{2m-2}{^t\gamma_2}$, and by (\ref{b2-formula-1}), with the fact that the common denominators of entries in $\delta_1'$ is still $\det(Y^{(r)})$, we conclude that the denominators of entries in $\gamma_2, \gamma_2^*$, and $b_2$ are monomials of $\det(Y_i)$, $\det(Y^{(i)}), (1\le i\le r)$, and $x_{r-i+1,i}, (1\le i\le \mathrm{min}\{r,m\})$. 

The same computation of $v_2(I_{2m}-J'_{2m}{^tX}{^tY^{-1}}J_rX)v_1=-t_2^{-1}J'_{2m}$ also implies that
\begin{equation}\label{gamma1-formula}
    x_{r,1}^2x\gamma_1-x_{r,1}\delta_2'J_{r-1}X_1v_1'=0,
\end{equation}
\begin{equation}\label{b1-formula}
    x_{r,1}^2xb_1-x_{r,1}\delta_2'J_{r-1}X_1\gamma_1^*+1-x_{r,1}\delta_2'J_{r-1}{^t\gamma}=0.
\end{equation}From the previous argument on $x$, $\delta_1', \delta_2'$, and the fact that $\gamma_1^*=v_1' J'_{2m-2}{^t\gamma_1}$, the formula (\ref{gamma1-formula}) and (\ref{b1-formula}) imply the same conclusion for denominators of entries in $\gamma_1$, $\gamma_1^*$, and $b_1$. 
\end{proof}
\section{Proof of Stability}
In the final section, we prove Theorem \ref{thm1.1} by proving the stability for the corresponding local coefficient. Recall that by the structure theorems it suffices to prove stability of local coefficient attached to $\psi$-generic supercuspidal representations. Given two $\psi$-generic supercuspidal representations $\sigma_i$, and $\tau_i$, $(i=1,2)$ of $\mathrm{GL}_r(F)$ and $\mathrm{Sp}_{2m}(F)$ respectively, with the same central character $\omega_{\sigma_1}=\omega_{\sigma_2}$, $\omega_{\tau_1}=\omega_{\tau_2}$. Denote $\pi_i=\sigma_i\boxtimes\tau_i$, then $\omega:=\omega_{\pi_1}=\omega_{\pi_2}$. Suppose $\chi$ is a continuous character of $F^\times$. Choose $f_1$ and $f_2$ to be matrix coefficients of $\pi_1$ and $\pi_2$ respectively, normalized so that $W_{f_1}(e)=W_{f_2}(e)=1$. We also choose $\kappa$ sufficiently large so that Proposition \ref{prop5.3} and Proposition \ref{prop6.7} hold for both $f_1$, $f_2$ and our fixed auxiliary function $f_0$ in Proposition \ref{prop6.7}, where $\varphi_{\overline{N}_0,\kappa}$ and $\varphi$ are related via Proposition \ref{prop6.9}.
By Proposition \ref{prop6.10}, 
$$C_\psi(s,\pi_1\otimes \chi)^{-1}-C_\psi(s,\pi_2\otimes \chi)^{-1}=\gamma(2rs, \omega^2\chi^{2r}, \psi^{-1})D_\chi(s)$$
where
$$D_\chi(s)=\int_{R_{X'}\times R_{Z'}}(B^M_\varphi(m(X',Z'), f_1)-B^M_\varphi(m(X',Z'), f_2))$$$$\cdot(\omega^{-2}\chi^{-2r})(\frac{P(X',Z')}{\det(Z'J_r+\frac{X' \theta_{r,m}(X')}{2})})\vert \frac{P(X',Z')}{\det(Z'J_r+\frac{X'\theta_{r,m}(X')}{2})}\vert^{-rs} \vert \det(m_1(X',Z'))\vert^{s+\frac{r+2m+1}{2}} d\mu_{X'}\wedge d\mu_{Z'}.$$
By Proposition \ref{prop6.7}, we can find $f_{1,w'}, f_{2,w'}\in C^\infty_c(\Omega_{w'}; \omega)$, for each $w'\in B(M)$ with $d_B(e,w')\ge 1$, such that $$B^M_{\varphi}(m(X',Z'),f_1)-B^M_\varphi(m(X',Z'),f_2)=\sum_{w'\in B(M), d_B(w',e)\ge 1}(B^M_{\varphi}(m(X',Z'),f_{1,w'})-B^M_{\varphi}(m(X',Z'),f_{2,w'})).$$

Since the toric action is determined by its restriction to $T_{X,Z}$ up to squares of torus entries according to Proposition \ref{prop6.11}, we call $T_{X,Z}$ the toric part of the orbit space $R_{X,Z}=R_{X}\times R_Z$. Denote by $N_{X,Z}$ the subset of $R_{X,Z}$ given by the variables not in $T_{X,Z}$, we call $N_{X,Z}$ the non-toric part of $R_{X,Z}$. Then by setting $x_{r,1}=1$ and pass to $R_{X'}\times R_{Z'}$, we can separate our integral over $R_{X',Z'}=R_{X'}\times R_{Z'}$ as a double integral over $N_{X',Z'}$ and $T_{X',Z'}$. 

Proposition \ref{prop6.12} implies that the assumption in Proposition \ref{prop6.8} is satisfied, since the denominators of rational functions in the 1-parameter subgroups for $u_1, u_2, v_1, v_2$ in the Bruhat decomposition $m_1=u_1\dot{w}_1t_1 u_2$, $m_2=v_1\dot{w}_2 t_2 u_2$ with $m=m(X',Z')=(m_1(X',Z'), m_2(X',Z'))$ are monomials of the image of $\phi$ on $T_{X',Z'}$.
The integral over $R_{X'}\times R_{Z'}$ breaks up into a double integral over $N_{X',Z'}$ and $T_{X',Z'}$. By Proposition \ref{prop6.11}, the restriction $\Phi\vert_{T_{X',Z'}}$ is finite $\acute{e}tale$ onto its image in $A'$, where $A=Z_MA'$, so the integral over $T_{X',Z'}$ can be written as a finite sum of integrals over $\Phi(T_{X',Z'})\subset A'$. In addition, since $A_{w_M}^{w'}A'_{w'}$ is open of finite index in $A'$, we can further write the integral over $\Phi(T_{X',Z'})$ as a double integral over $\Phi(T_{X',Z'})\cap A_{w_m}^{w'}$ and $\Phi(T_{X',Z'})\cap A'_{w'}$, so we have    
$$D_\chi(s)=\sum_{\mathrm{finite}}\sum_{1\le d_B(e,w')}\int_{N_{X',Z'}}\int_{\Phi(T_{X',Z'})\cap A_{w_M}^{w'}}\int_{\Phi(T_{X',Z'})\cap A'_{w'}}(B^M_\varphi(m(X',Z'), f_{1,w'})-B^M_\varphi(m(X',Z'), f_{2,w'}))$$$$\cdot(\omega^{-2}\chi^{-2r})(\frac{P(X',Z')}{\det(Z'J_r+\frac{X' \theta_{r,m}(X')}{2})})\vert \frac{P(X',Z')}{\det(Z'J_r+\frac{X'\theta_{r,m}(X')}{2})}\vert^{-rs} \vert \det(m_1(X',Z'))\vert^{s+\frac{r+2m+1}{2}} d\mu_{X'}\wedge d\mu_{Z'}.$$

By the proof of Proposition \ref{prop6.8}, the functions $B^M_\varphi(m(X',Z'), f_{i,w'})$, $(i=1,2)$, are compactly supported in the subtorus $A'_{w'}=Z_L'$ where $w'=w^M_L\in B(M)$. By the compatibility of the toric action as in Proposition \ref{prop6.11}, there exists an open compact subgroup $Z'_{L,0}$ of $Z'_L(F)$,  such that given $m=m(X',Z')$, for $s=(s',s'')\in Z'_{L,0}$, we have $$B^M(\Theta_M(s)m(X',Z')s^{-1}, f_{i,w'})=B^M_\varphi(m(s'X'{s''}^{-1}, s' Z's'), f_{i,w'})= B_\varphi^M(m(X',Z'), f_{i,w'})$$ for all $\kappa$ where $\varphi_{\overline{N}_{\kappa,0}}$ is related to $\varphi$ via Proposition \ref{prop6.9}.  
In our cases, 
$L\simeq \prod_{i=1}^k\mathrm{GL}_{r_i}\times \prod_{j=1}^t \mathrm{GL}_{m_i}\times \mathrm{Sp}_{2m'}$ with $\sum_{i=1}^kr_i=r$ and $\sum_{j=1}^t m_j+m'=m$. So we can pick $s=(s', s'')\in Z'_{L,0}$ with $$s'=\mathrm{diag}\{\underbrace{s_1,\cdots, s_1}_{r_1},\underbrace{s_2,\cdots,s_2}_{r_2},\cdots, \underbrace{s_k,\cdots, s_k}_{r_k}\},$$
$$s''=\mathrm{diag}\{\underbrace{s_{k+1},\cdots,s_{k+1}}_{m_1},\cdots,\underbrace{s_{k+l},\cdots, s_{k+1}}_{m_l}, \underbrace{1,\cdots, 1}_{2m'}, \underbrace{s^{-1}_{k+l},\cdots, s^{-1}_{k+1}}_{m_l},\cdots, \underbrace{s^{-1}_{k+1},\cdots,s^{-1}_{k+1}}_{m_1} \}.$$
Recall that the term $\frac{P(X',Z')}{\det(Z'J_r+\frac{X' \theta_{r,m}(X')}{2})}$ is just $\frac{{y'}^*_{rr}}{\det(Y')}$, where ${y'}^*_{rr}$ is the $(r,r)$-th entry of the adjoint matrix of $Y'$. Here $Y'=Z'J_r+\frac{X'\theta_{r,m}(X')}{2}$. One computes that the toric action takes $\frac{{y'}^*_{rr}}{\det(Y')}$ to $\frac{1}{s_k^2}\frac{{y'}^*_{rr}}{\det(Y')}$, and
takes $\det(m_1(X',Z'))$ to $\prod_{i=1}^ks_i^2\det(m_1(X',Z'))$, since $m_1(X',Z')=\theta_r(Y')$. The action of $Z_{L,0}$ preserves $\Phi(T_{X',Z'})\cap A'_{w'}$, therefore if we change variables under the toric action, the inner integral is equal to
$$(\omega\chi^{r})(s^4_k)\vert s_k\vert^{rs}\prod_{i=1}^k\vert s_i\vert^{2s+r+2m+1}\int_{\Phi(T_{X',Z'})\cap A'_{w'}}(B^M_\varphi(m(X',Z'), f_{1,w'})-B^M_\varphi(m(X',Z'), f_{2,w'}))$$$$\cdot(\omega^{-2}\chi^{-2r})(\frac{P(X',Z')}{\det(Z'J_r+\frac{X' \theta_{r,m}(X')}{2})})\vert \frac{P(X',Z')}{\det(Z'J_r+\frac{X'\theta_{r,m}(X')}{2})}\vert^{-rs} \vert \det(m_1(X',Z'))\vert^{s+\frac{r+2m+1}{2}}.$$ 
So $$D_\chi(s)=(\omega\chi^{r})(s^4_k)\vert s_k\vert^{rs}\prod_{i=1}^k\vert s_i\vert^{2s+r+2m+1} D_\chi(s).$$
Choose $\chi$ sufficiently ramified so that $(\omega\chi^{r})(s^4_k)\vert s_k\vert^{rs}\prod_{i=1}^k\vert s_i\vert^{2s+r+2m+1}\neq 1$, therefore $D_\chi(s)=0$. We finally conclude that
$$C_\psi(s,\pi_1\otimes \chi)=C_\psi(s,\pi_2\otimes \chi).$$
Consequently, $$\gamma(s,(\sigma_1\times \tau_1)\otimes \chi, \psi)=\gamma(s,(\sigma_2\times \tau_2)\otimes\chi, \psi)$$ for sufficiently ramified $\chi$.

\bibliographystyle{plain}
\bibliography{biblio}

\begin{thebibliography}{10}

\bibitem{Bou08}
Nicolas Bourbaki.
\newblock {\em Lie Groups and Lie Algebras Chapters 4-6}.
\newblock Springer Berlin, Heidelberg, 2008.

\bibitem{CpS05}
James~W Cogdell, Ilya~Iosifovich Piatetski-Shapiro, and Freydoon Shahidi.
\newblock Partial bessel functions for quasi-split groups.
\newblock {\em Automorphic Representations, L-functions and Applications:
  Progress and Prospects. Berlin: Walter de Gruyter}, pages 95--128, 2005.

\bibitem{Cog08}
James~W Cogdell, Ilya~Iosifovich Piatetski-Shapiro, and Freydoon Shahidi.
\newblock Stability of-factors for quasi-split groups.
\newblock {\em Journal of the Institute of Mathematics of Jussieu},
  7(1):27--66, 2008.

\bibitem{CST17}
James~W Cogdell, Freydoon Shahidi, and T-L Tsai.
\newblock Local langlands correspondence for {$GL_n$} and the exterior and
  symmetric square $\varepsilon$-factors.
\newblock {\em Duke Mathematical Journal}, 166(11):2053--2132, 2017.

\bibitem{Deligne73}
Pierre Deligne.
\newblock Les constantes des equations fonctionnelles des fonctions l.
\newblock {\em Modular Functions of One Variable, II(Antwerp,1972). Lecture
  Notes in Math.}, 349:501--597, 1973.

\bibitem{Jac85}
Herv{\'e} Jacquet and Joseph Shalika.
\newblock A lemma on highly ramified $\varepsilon$-factors.
\newblock {\em Math. Ann}, 271:319--332, 1985.

\bibitem{Ral05}
Stephen Rallis and David Soudry.
\newblock Stability of the local gamma factor arising from the doubling method.
\newblock {\em Mathematische Annalen}, 333(2):291--313, 2005.

\bibitem{sha90}
Freydoon Shahidi.
\newblock A proof of langlands's conjecture on plancharel measures,
  complementary series for p-adic groups.
\newblock {\em Annals of Mathematics}, 132:273--330, 1990.

\bibitem{Sha02}
Freydoon Shahidi.
\newblock Local coefficients as mellin transforms of bessel functions: Towards
  a general stability.
\newblock {\em International Mathematics Research Notices}, 39:2075--2119,
  2002.

\bibitem{sha10}
Freydoon Shahidi.
\newblock {\em Eisenstein Series and Automorphic L-Functions}, volume~58.
\newblock American Mathematical Society, 2010.

\bibitem{Shan18}
Daniel Shankman.
\newblock Local langlands correspondence for asai l-functions and epsilon
  factors.
\newblock {\em arXiv preprint arXiv:1810.11852}, 2018.

\bibitem{Shan19}
Daniel Shankman and Dongming She.
\newblock Stability of symmetric cube gamma factors for $gl (2)$.
\newblock {\em arXiv preprint arXiv:1911.03428}, 2019.

\bibitem{She23}
Dongming She.
\newblock Local langlands correspondence for the twisted exterior and symmetric
  square $\epsilon $-factors of $ \text {GL} \_n $.
\newblock {\em manuscripta mathematica}, pages 1--47, 2023.

\bibitem{Sun07}
Rajan Sundaravaradhan.
\newblock Some structural results for the stability of root numbers.
\newblock {\em International Mathematics Research Notices}, pages 1--22, 2007.

\end{thebibliography}

\end{document}